\documentclass[11pt,reqno]{amsart}
\usepackage[letterpaper,margin=1.05in]{geometry}
\usepackage[T1]{fontenc}
\usepackage{amsmath,amssymb,amsthm,mathtools,microtype}
\usepackage{lmodern}
\DeclareMicrotypeSet{T1text}{encoding=T1,family={rm*,sf*}}
\UseMicrotypeSet[expansion]{T1text}
\usepackage{xcolor,mathrsfs}
\usepackage[hyphens]{url}
\usepackage{hyperref}

\definecolor{annalsblue}{RGB}{18,53,110}
\hypersetup{
  colorlinks=true,
  linkcolor=annalsblue,
  citecolor=annalsblue,
  urlcolor=annalsblue,
  pdfencoding=auto, hypertexnames=false,
}

\numberwithin{equation}{section}

\newtheorem{theorem}{Theorem}[section]
\newtheorem{proposition}[theorem]{Proposition}
\newtheorem{lemma}[theorem]{Lemma}

\theoremstyle{remark}

\theoremstyle{definition}
\newtheorem{definition}[theorem]{Definition}

\newcommand{\ind}{\mathbf{1}}
\newcommand{\W}{\mathcal{W}}
\newcommand{\Pset}{\mathcal{P}}
\newcommand{\Spr}{\mathfrak{S}}
\newcommand{\dd}{\,\mathrm{d}}
\DeclareMathOperator{\supp}{supp}

\newcommand{\modf}{\,(\text{mod}\,4)}
\newcommand{\modq}{\,(\text{mod}\,q)}

\title[Primes of the form $u^2+v^2+1$]{An asymptotic formula for primes of the form $u^2+v^2+1$}

\author{Yujiao Jiang}

\address{Yujiao Jiang\\
	School of Mathematics and Statistics
	\\
	Shandong University
	\\
	Weihai
	\\
	Shandong 264209
	\\
	China}
\email{yujiaoj@sdu.edu.cn}

\author{Yuk-Kam Lau}

\address{Yuk-Kam Lau\\
	Department of Mathematics
	\\ 
	The University of Hong Kong
	\\ 
	Pokfulam Road
	\\ 
	Hong Kong\\
	China}
\email{yklau@maths.hku.hk}

\author{Guangshi L\"u}

\address{Guangshi L\"u\\
	School of Mathematics
	\\
	Shandong University
	\\
	Jinan
	\\
	Shandong 250100
	\\
	China}
\email{gslv@sdu.edu.cn}

\date{\today}


\newcommand{\Aset}{\mathcal A}
\newcommand{\Mmajor}{\mathfrak M}
\newcommand{\ip}[2]{\langle #1,#2\rangle}

\newtheorem{corollary}[theorem]{Corollary}

\begin{document}
\begin{abstract}
We prove Motohashi's conjecture on primes of the form $u^2+v^2+1$, with Iwaniec's corrected constant. Let $b$ denote the characteristic function of the positive integers representable as a sum of two squares. Then, as $x\to\infty$,
\[
 \sum_{p\le x}b(p-1)
 =\kappa\Spr\frac{x}{(\log x)^{3/2}}
 +O\!\left(\frac{x}{(\log x)^{3/2}\sqrt{\log\log x}}\right),
\]
where $\kappa$ is the Landau--Ramanujan constant and $\Spr=\prod_{p\equiv3\,({\rm mod}\,4)}\bigl(1-1/(p(p-1))\bigr)$. The proof combines auxiliary prime factors with a mean square estimate for bilinear sums. We obtain the minor arc bound through harmonic switching and high moments, replacing arithmetic divisibility conditions by independent projective lines. We remove the auxiliary weights using Selberg's sieve.
\end{abstract}
\maketitle
\setcounter{tocdepth}{2}
\tableofcontents

\section{Introduction}

Let $b$ denote the characteristic function of the positive integers of the form $u^2+v^2$, where $u,v\in\mathbb Z$. By the Fermat--Euler criterion \cite[Section~1.4]{IK},
\begin{equation}\label{mot:eq:b-def}
 b(p^j)=\begin{cases}
 1,&p=2\text{ or }p\equiv1\modf,\\
 \ind_{2\mid j},&p\equiv3\modf.
 \end{cases}
\end{equation}
Landau \cite{Landau} proved an asymptotic formula for the summatory function of $b$:
\[
 \sum_{n\le x}b(n)\sim\kappa\,\frac{x}{\sqrt{\log x}},
\]
where $\kappa=\frac{1}{\sqrt{2}}\prod_{p\equiv3\modf}\bigl(1-p^{-2}\bigr)^{-1/2}$ is the Landau--Ramanujan constant. We seek an asymptotic formula for $\sum_{p\le x}b(p-1)$, the number of primes $p\le x$ of the form $u^2+v^2+1$. Linnik \cite{Linnik} and Bredikhin \cite{Bredikhin} proved that there are infinitely many such primes. Motohashi \cite{Motohashi,MotohashiII} obtained a lower bound of order $x(\log x)^{-2}$ and formulated an asymptotic conjecture for this count (Conjecture~J$^*$ \cite{Motohashi}). Iwaniec \cite[p.~204]{Iwaniec72} corrected its leading constant and used the half-dimensional sieve to establish the true order of magnitude
\begin{equation}\label{eq:iwaniec-order}
	\sum_{p\le x} b(p-1) \asymp \frac{x}{(\log x)^{3/2}}.
\end{equation}

We prove the asymptotic formula conjectured by Motohashi, with Iwaniec's corrected constant and a quantitative error term.

\begin{theorem}\label{mot:thm:main}
As $x\to\infty$, we have
\[
 \sum_{p\le x}b(p-1)
 =\kappa\Spr\frac{x}{(\log x)^{3/2}}
 +O\!\left(\frac{x}{(\log x)^{3/2}\sqrt{\log\log x}}\right),
\]
where $\kappa$ is the Landau--Ramanujan constant and $\Spr=\prod_{p\equiv 3\, ({\rm mod}\, 4)}\bigl(1-\frac{1}{p(p-1)}\bigr)>0$.
\end{theorem}

Iwaniec \cite[Corollary~2, (2.1)]{Iwaniec76} related the error in this asymptotic formula to the distribution of primes in arithmetic progressions:
\[
\sum_{p<x}b(p-1)
=\kappa\Spr\frac{x}{(\log x)^{3/2}}
\bigl(1+O(\varepsilon^\gamma)\bigr)
+\theta_1 E(x,x^{1-\varepsilon}),
\]
where $\varepsilon>\log\log x/\log x$, $\gamma>0$, $|\theta_1|\le1$, and
\[
E(x,Q)=\sum_{q<Q}\max_{(a,q)=1}
\Bigl|\pi(x;q,a)-\frac{\operatorname{Li}(x)}{\varphi(q)}\Bigr|.
\]
The Elliott--Halberstam conjecture predicts $E(x,x^{1-\varepsilon})\ll_{\varepsilon,A}x(\log x)^{-A}$ for every fixed $0<\varepsilon<1$ and $A>0$. As Iwaniec observed after \cite[(2.3)]{Iwaniec76}, this bound would yield the conjectured asymptotic.

\subsection{Outline of the proof}
Write $L=\log x$. We first establish an asymptotic formula for the weighted sum $\Theta_{b,\W}(x,\Phi)=\sum_n\Lambda(n)b(n-1)\W(n-1)\Phi(n/x)$, where $\Phi$ is smooth and supported in $[1/2,1]$. We form $\W$ from the numbers of prime divisors in $K$ disjoint sets $\Pset_\nu$ of primes congruent to $1$ modulo $4$. These sets lie in the separated ranges defined in \eqref{common:eq:application-bands}.

Vaughan's identity, in Section~\ref{sec:common}, reduces the weighted sum to Type~I and Type~II sums. Section~\ref{sec:typeI} treats the Type~I sums through the distribution of $b(n)\W(n)$ in arithmetic progressions. We deduce the required Bombieri--Vinogradov estimate from Fouvry and Tenenbaum's theorem \cite[Theorem~1.8]{FT} by applying Cauchy's formula to a family of multiplicative functions with complex parameters. A Dirichlet convolution with $b$ then reduces the evaluation of the main term to Landau's asymptotic formula and an Euler product.

The Type~II sums have the form
\[
 \sum_{m\asymp M}\sum_{n\asymp N}
 \alpha_m\beta_n b(mn-1)\W(mn-1),
 \qquad MN\asymp x,\quad M,N\ge x^{1/5},
\]
where $\alpha_m=\mu(m)m^{iv}$ on its support. The support and coefficient conditions are stated in Proposition~\ref{common:prop:typeII}. Expanding $\W(mn-1)$ extracts a divisor $a=p_1\cdots p_K$, with $p_\nu\in\Pset_\nu$. Since every $p_\nu$ is congruent to $1$ modulo $4$, the relation $mn-1=ah$ implies $b(mn-1)=b(h)$ without a coprimality condition. After dividing $a$ into dyadic intervals, Cauchy's inequality over $h$ removes $b(h)$ and reduces the problem to the mean square $Q_Y$ in \eqref{core:eq:typeII-Q}. For each prescribed $R>0$, Theorem~\ref{core:quadratic-theorem} proves $Q_Y\ll XYL^{-R}$, where $X=MN$, once $K$ is sufficiently large.

Expanding this mean square pairs the equations $mn-1=ah$ and $rs-1=ch$, and hence introduces the determinant relation
\[
 cmn-ars=c-a.
\]
We compare $Q_Y$ with its major arc approximation. Sections~\ref{core:sec:minor-reduction} and~\ref{core:sec:typeII-minor} estimate the difference by harmonic switching and high moments of an operator on labelled primitive points. We follow closely the method of \cite{PD,LWW}, adapting it to the auxiliary weights $\W$ considered here. We establish the high moment estimate and the minor arc comparison in Proposition~\ref{core:lem:restricted-determinant}.

Section~\ref{core:sec:typeII-major} estimates the major arcs. We write $c=pc'$ with $p\in\Pset_1$ and, for fixed $c'$, obtain a prime Dirichlet polynomial of length $x^{o(1)}$. We collect the unit intervals containing a point where this polynomial is large; \cite[Lemma~6.5]{KMT} bounds their number. The contribution outside their union is negligible. On these intervals, we combine M\"obius cancellation from Lemma~\ref{core:lem:sifted-mobius} with the mean square estimate in Lemma~\ref{core:lem:prime-sampling}. The M\"obius estimate uses the bound for $1/L(s,\chi)$ in \cite[Theorem~24.20]{MVIII}, together with the exclusion of Landau--Siegel zeros \cite{BanksSiegel,GreenSiegel}. This completes the proof of Theorem~\ref{core:quadratic-theorem}.

Finally, Section~\ref{mot:sec:application} removes the auxiliary weights using Selberg's sieve. Smooth approximation then completes the proof of Theorem~\ref{mot:thm:main}.

\section{Preliminary reductions}\label{sec:common}
\subsection{Auxiliary primes and weights}
For the rest of the paper, write $L=\log x$ and $z=\exp(L^{0.24})$. Put $P(t)=\prod_{p\le t}p$. Write $P^-(n)$ and $P^+(n)$ for the least and greatest prime factors, with $P^-(1)=\infty$ and $P^+(1)=1$. Write $\omega(n)=\sum_{p\mid n}1$ for the number of distinct prime factors. We use $e(t)=\exp(2\pi it)$ and $e_q(t)=e(t/q)$. An integer $n$ is called $z$-rough if $P^-(n)>z$. Let $K\ge1$ be a fixed integer. Let $J=J(x)$ be an integer satisfying
\begin{equation}\label{common:eq:family-range}
 1\le J\le\left\lfloor\frac{\log L}{20K\log2}\right\rfloor.
\end{equation}
Choose
\begin{equation}\label{common:eq:family-exponents}
 \vartheta_{j,\nu}
 =\frac1{10}+\frac{(\nu-1)J+j}
                    {10(KJ+1)}
 \qquad(1\le j\le J,\ 1\le\nu\le K).
\end{equation}
All $KJ$ exponents form an equally spaced grid in $(1/10,1/5)$. Within each family,
\[
 \vartheta_{j,\nu+1}-\vartheta_{j,\nu}
 =\frac{J}{10(KJ+1)}
 \ge\frac1{10(K+1)}\ge\frac1{20K}.
\]
When working with one family we write $\vartheta_\nu$ for $\vartheta_{j,\nu}$, so that
\begin{equation}\label{common:eq:band-exponents}
 \frac1{10}<\vartheta_1<\cdots<\vartheta_K<\frac15,\qquad
 \vartheta_{\nu+1}-\vartheta_\nu\ge\frac1{20K}
 \quad(1\le\nu<K).
\end{equation}
We use
\begin{equation}\label{common:eq:application-bands}
 \Pset_\nu=\{p:\exp(L^{\vartheta_\nu})<p\le\exp(2L^{\vartheta_\nu}),
                      \ p\equiv1\modf\}.
\end{equation}
These sets, and hence the weights below, depend on $x$. They are held fixed when an auxiliary summation limit varies. Define
\begin{equation}\label{common:eq:weight}
 V_\nu=\sum_{p\in\Pset_\nu}\frac1p,\qquad
 \omega_\nu(n)=\sum_{p\in\Pset_\nu}\ind_{p\mid n},\qquad
 \W(n)=\prod_{\nu=1}^K\frac{\omega_\nu(n)}{V_\nu}.
\end{equation}
\begin{equation}\label{common:eq:products}
 \Pset_*=\bigcup_\nu\Pset_\nu,\qquad
 \mathcal A(\Pset)=\{p_1\cdots p_K:p_\nu\in\Pset_\nu\},\qquad
 \mathcal V(\Pset)=\prod_{\nu=1}^K V_\nu.
\end{equation}
For \eqref{common:eq:application-bands}, the prime number theorem in fixed progressions \cite[Corollary~11.19]{MV} and partial summation imply $V_\nu=\frac12\log2+o(1)$. In particular, $1/4\le V_\nu\le1/2$ for sufficiently large $x$, uniformly in the exponents in \eqref{common:eq:band-exponents}. The weight vanishes unless $n$ has a prime divisor from every set $\Pset_\nu$. Moreover,
\[
 0\le\W(n)=\frac1{\mathcal V(\Pset)}
       \sum_{\substack{a\in\mathcal A(\Pset)\\a\mid n}}1
       \le4^K d(n).
\]
Here $d(n)=\sum_{a\mid n}1$. All primes in $\Pset_*$ lie below $z$ for sufficiently large $x$, and a product of one prime from each set is $x^{o(1)}$. Choose once a nonnegative $\eta\in C_c^\infty([1,4])$ with $\sum_{j\in\mathbb Z}\eta(t/2^j)=1$ for $t>0$. In particular, $0\le\eta\le1$. For all sufficiently large $x$, \eqref{common:eq:family-range} implies $\log L/[10(KJ+1)]>\log2$. Hence all $KJ$ prime sets associated with \eqref{common:eq:family-exponents} are pairwise disjoint, uniformly in $J$. The gap within each family still satisfies \eqref{common:eq:band-exponents}. The integer $J_0$ introduced in the mean square proof counts shared prime labels. Write $\Pset^{(j)}=(\Pset_{j,1},\ldots,\Pset_{j,K})$ for the family obtained from \eqref{common:eq:application-bands} with exponents $\vartheta_{j,\nu}$. Let $V_{j,\nu}$, $\omega_{j,\nu}$ and $\W_j$ denote the quantities in \eqref{common:eq:weight} for this family, and put
\[
 \mathcal E=\bigcup_{j=1}^{J}\bigcup_{\nu=1}^K\Pset_{j,\nu}.
\]

\subsection{Vaughan's identity}
The prime sets are given by \eqref{common:eq:application-bands}. Extend $b$ by zero on nonpositive integers, and put $f(n)=b(n)\W(n)$. Then $0\le f(n)\le4^K d(n)$. Let $\Phi\in C_c^\infty((0,\infty))$ be supported in $[1/2,1]$. For every integer $B\ge0$, put $\mathcal N_B(\Phi)=\sum_{r=0}^B\|\Phi^{(r)}\|_\infty$. We record the dependence on $\Phi$ through these seminorms. Put
\[
 \Theta_{b,\W}(x,\Phi)=\sum_{n\ge2}\Lambda(n)b(n-1)\W(n-1)\Phi(n/x).
\]
Vaughan's identity \cite[Proposition~13.4, (13.39)]{IK} takes the form
\begin{equation}\label{common:eq:vaughan}
 \Theta_{b,\W}=T_1-T_2+T_3+O_{A,K}(\|\Phi\|_\infty xL^{-A}),
\end{equation}
where
\begin{align}
 T_1&=\sum_{\substack{d\le U_0\\P^-(d)>z}}\mu(d)
       \sum_{\substack{j\ge1\\P^-(j)>z}}(\log j)
       f(dj-1)\Phi(dj/x),\notag\\
 T_2&=\sum_{\substack{d,\ell\le U_0\\P^-(d\ell)>z}}\mu(d)\Lambda(\ell)
       \sum_{\substack{j\ge1\\P^-(j)>z}}
       f(d\ell j-1)\Phi(d\ell j/x),\notag\\
 T_3&=\sum_{\substack{d,t>U_0\\P^-(dt)>z}}\mu(d)c_{U_0}(t)
       f(dt-1)\Phi(dt/x).\label{common:eq:T3}
\end{align}
Here
\[
 U_0=x^{1/5},\qquad
 c_{U_0}(t)=\sum_{\substack{\ell\mid t\\\ell>U_0}}\Lambda(\ell)\le\log t.
\]
Indeed, the restriction to $z$-rough integers only removes prime powers $p^j$ with $p\le z$ and $j\ge2$. Their contribution is $O_{K,\epsilon}(\|\Phi\|_\infty zx^\epsilon L)=O_{A,K}(\|\Phi\|_\infty xL^{-A})$.

With $n=dj-1$ in $T_1$, we have $n\equiv-1\pmod d$; removing the restriction $P^-(j)>z$ by M\"obius inversion introduces moduli $dc$ in $T_1$ and $d\ell c$ in $T_2$, where $c\mid P(z)$. Section~\ref{sec:typeI} uses the distribution of $f$ in arithmetic progressions to evaluate $T_1-T_2$.

\section{The Type~I asymptotic}\label{sec:typeI}
We retain the notation of Section~\ref{sec:common}. We first establish a Bombieri--Vinogradov estimate for $f=b\W$ and reduce $T_1-T_2$ to a weighted mean of $b$. A Dirichlet convolution with $b$ and an evaluation of the Euler factors then determine its main term.

\subsection{A Bombieri--Vinogradov theorem for the weighted indicator}
For $4\nmid q$ and $(c,q)=1$, define
\[
 \Delta_f(y;q,c)=\sum_{\substack{n\le y\\n\equiv c\modq}}f(n)
       -\frac1{\varphi(q)}\sum_{\substack{n\le y\\(n,q)=1}}f(n).
\]

\begin{lemma}\label{common:lem:distribution}
For every fixed $A>0$,
\begin{equation*}
 \sup_{x/4\le y\le x}
 \sum_{\substack{q\le x^{0.41}\\4\nmid q}}
       \max_{(c,q)=1}|\Delta_f(y;q,c)|
       \ll_{A,K}xL^{-A}.
\end{equation*}
The prime sets defining $f$ remain fixed by $x$ throughout the supremum.
\end{lemma}

\begin{proof}
For complex parameters $\boldsymbol t=(t_1,\ldots,t_K)$ with $|t_\nu-1|=1/2$, consider the multiplicative function
\[
 b_{\boldsymbol t}(n)=b(n)\prod_{\nu=1}^K t_\nu^{\omega_\nu(n)}.
\]
For each fixed $n$, this is a polynomial in the parameters. On these circles $|t_\nu|\le3/2$, so $|b_{\boldsymbol t}(p^j)|\le3/2\le d(p^j)$ for $j\ge1$, and hence $|b_{\boldsymbol t}(n)|\le d(n)$. Its values at primes are
\[
 b_{\boldsymbol t}(p)=
 \begin{cases}
  b(p)t_\nu,&p\in\Pset_\nu,\\
  b(p),&p\notin\Pset_*.
 \end{cases}
\]
Outside $\Pset_*$, the values at odd primes are $1$ and $0$ in the classes $1$ and $3$ modulo $4$, respectively.

As in \cite[Definition~1.1]{FT}, use the increasing sequence formed by
\[
 2,\quad \exp(L^{\vartheta_\nu}),\quad
 \exp(2L^{\vartheta_\nu})\quad(1\le\nu\le K),\quad
 2^j\exp(2L^{\vartheta_K})\quad(j\ge1).
\]
For sufficiently large $x$, successive terms have ratio at least $2$, by \eqref{common:eq:band-exponents}. Thus they satisfy the spacing condition in that definition. If two primes lie in the same interval of this sequence and are congruent modulo $4$, the displayed prime values agree. It follows that $b_{\boldsymbol t}\in\mathcal F(4,2)$, uniformly in $\boldsymbol t$ on these circles.

For $q>1$ with $4\nmid q$, put $q_4=(q,4^\infty)=(q,4)\in\{1,2\}$. For $(c,q)=1$, the main term in \cite[(1.15)]{FT} is
\[
 \frac1{\varphi(q/q_4)}
 \sum_{\substack{n\le y\\n\equiv c\,({\rm mod}\,q_4)\\(n,q/q_4)=1}}
       b_{\boldsymbol t}(n)
 =\frac1{\varphi(q)}\sum_{\substack{n\le y\\(n,q)=1}}
       b_{\boldsymbol t}(n).
\]
For odd $q$, every reduced residue class has an odd representative; for even $q$, every reduced representative is odd. Hence the maxima over $(c,q)=1$ and $(c,4q)=1$ coincide. For $q=1$ the discrepancy is zero.

For every fixed $B>0$, $x^{0.41}<\sqrt y(\log y)^{-B}$ uniformly for $x/4\le y\le x$ once $x$ is sufficiently large. Theorem~1.8 of \cite{FT} is uniform over $\mathcal F(4,2)$. Applying it at each such $y$, with the prime sets and $b_{\boldsymbol t}$ still fixed by $x$, we obtain
\begin{equation}\label{common:eq:BV-family}
 \sup_{x/4\le y\le x}\ 
 \sup_{\substack{|t_\nu-1|=1/2\\1\le\nu\le K}}
 \sum_{\substack{q\le x^{0.41}\\4\nmid q}}
 \max_{(c,q)=1}|\Delta_{b_{\boldsymbol t}}(y;q,c)|
 \ll_A xL^{-A}.
\end{equation}
Here $\Delta_{b_{\boldsymbol t}}$ is defined in the same way as $\Delta_f$.

By \eqref{common:eq:weight},
\[
 f(n)=
 \biggl(\prod_{\nu=1}^K V_\nu^{-1}\frac{\partial}{\partial t_\nu}\biggr)
 b_{\boldsymbol t}(n)\bigg|_{t_1=\cdots=t_K=1}.
\]
Thus, by Cauchy's formula,
\begin{equation}\label{common:eq:BV-cauchy}
 \Delta_f(y;q,c)=
 \frac1{(2\pi i)^K\mathcal V(\Pset)}
 \oint\cdots\oint
 \frac{\Delta_{b_{\boldsymbol t}}(y;q,c)}
      {\prod_{\nu=1}^K(t_\nu-1)^2}
       \,\dd t_K\cdots\dd t_1,
\end{equation}
where the $\nu$th integral is over the positively oriented circle $|t_\nu-1|=1/2$. Take absolute values in \eqref{common:eq:BV-cauchy} and move the sum over $q$ and the maximum over $c$ inside the integrals. Since $V_\nu^{-1}\le4$ for sufficiently large $x$, we obtain
\[
\sup_{x/4\le y\le x}
   \sum_{\substack{q\le x^{0.41}\\4\nmid q}}
       \max_{(c,q)=1}|\Delta_f(y;q,c)|\le8^K
 \sup_{x/4\le y\le x}\ 
 \sup_{\substack{|t_\nu-1|=1/2\\1\le\nu\le K}}
   \sum_{\substack{q\le x^{0.41}\\4\nmid q}}
       \max_{(c,q)=1}|\Delta_{b_{\boldsymbol t}}(y;q,c)|
 \ll_{A,K}xL^{-A},
\]
by \eqref{common:eq:BV-family}. This proves the lemma.
\end{proof}

\subsection{Reduction of the Type~I sums}
Put
\begin{equation*}
 W(z)=\prod_{p\le z}(1-p^{-1}),\qquad
 w_z(n)=\prod_{\substack{p\le z\\p\nmid n}}
                      \Bigl(1-\frac1{p-1}\Bigr).
\end{equation*}
Here $0\le w_z(n)\le1$, and the factor at $2$ forces $w_z(n)=0$ for odd $n$.

\begin{lemma}\label{common:lem:typeI-model}
For every fixed $A>0$,
\begin{equation*}
 T_1-T_2=\frac1{W(z)}\sum_{n\ge1}b(n)\W(n)w_z(n)
       \Phi((n+1)/x)+O_{A,K}(\mathcal N_1(\Phi)xL^{-A}).
\end{equation*}
\end{lemma}
\begin{proof}
By homogeneity we may assume $\mathcal N_1(\Phi)\le1$. All implied constants in this proof may depend on $K$. We first evaluate $T_1$. By M\"obius inversion,
\[
 T_1=\sum_{\substack{d\le U_0\\(d,P(z))=1}}\mu(d)
     \sum_{c\mid P(z)}\mu(c)
     \sum_{\substack{n\ge1\\n\equiv-1\,({\rm mod}\,dc)}}f(n)\Psi_d(n),
\]
where $\Psi_d(n)=\Phi((n+1)/x)\log((n+1)/d)$. We truncate the $c$-sum at $R_0=\exp(L^{1/2})$. For each $d\le U_0$ with $(d,P(z))=1$, interchanging the two inner sums, we get
\[
\sum_{\substack{c\mid P(z)\\c>R_0}}\mu(c)
       \sum_{\substack{n\ge1\\n\equiv-1\,({\rm mod}\,dc)}}f(n)\Psi_d(n)=\sum_{\substack{n\ge1\\d\mid n+1}}f(n)\Psi_d(n)
       \sum_{\substack{c\mid(n+1,P(z))\\c>R_0}}\mu(c).
\]
Here $(d,c)=1$, and $f(n)\Psi_d(n)$ is independent of $c$. For each $n$, one has
\[
 \bigg|\sum_{\substack{c\mid(n+1,P(z))\\c>R_0}}\mu(c)\bigg|
 \le\sum_{\substack{c\mid(n+1,P(z))\\R_0<c\le zR_0}}1.
\]
By Shiu's theorem \cite[Theorem~1]{Shiu},
\[
 \sum_{\substack{n\asymp x\\n\equiv-1\,({\rm mod}\,q)}}d(n)
       \ll\frac{xL}{\varphi(q)}
       \qquad(q\le x^{0.41}).
\]
Since $f(n)\ll_K d(n)$, $|\Psi_d(n)|\ll L$, and $dc\le U_0zR_0<x^{0.41}$, the preceding two bounds show that the omitted contribution to $T_1$ is
\[
\begin{aligned}
 &\ll L\sum_{\substack{d\le U_0\\(d,P(z))=1}}
       \sum_{\substack{c\mid P(z)\\R_0<c\le zR_0}}
       \sum_{\substack{n\asymp x\\n\equiv-1\,({\rm mod}\,dc)}}d(n)\\
 &\ll xL^2
       \sum_{\substack{d\le U_0\\(d,P(z))=1}}\frac1{\varphi(d)}
       \sum_{\substack{c\mid P(z)\\R_0<c\le zR_0}}\frac1{\varphi(c)}\\
 &\ll xL^3\sum_{\substack{c\mid P(z)\\c>R_0}}\frac1{\varphi(c)},
\end{aligned}
\]
where $\sum_{d\le U_0}\varphi(d)^{-1}\ll L$. It is known from Rankin's trick that
\begin{equation}\label{common:eq:typeI-sieve-tail}
 \sum_{\substack{c\mid P(z)\\c>R_0}}\frac1{\varphi(c)}
 \le R_0^{-1/\log z}
       \prod_{p\le z}\Bigl(1+\frac{p^{1/\log z}}{p-1}\Bigr)
 \ll\exp\bigl(-L^{0.26}+O(\log L)\bigr)\ll_A L^{-A}.
\end{equation}
Applying \eqref{common:eq:typeI-sieve-tail} with exponent $A+4$, we obtain an error $O_A(xL^{-A})$ for the omitted contribution. Thus, we have
\begin{equation}\label{common:eq:typeI-T1-truncated}
 T_1=\sum_{\substack{d\le U_0\\(d,P(z))=1}}\mu(d)
     \sum_{\substack{c\mid P(z)\\c\le R_0}}\mu(c)
     \sum_{\substack{n\ge1\\n\equiv-1\,({\rm mod}\,dc)}}f(n)\Psi_d(n)
     +O_A(xL^{-A}).
\end{equation}

By partial summation,
\[
\sum_{\substack{n\ge1\\n\equiv-1\,({\rm mod}\,dc)}}f(n)\Psi_d(n)
 -\frac1{\varphi(dc)}
       \sum_{\substack{n\ge1\\(n,dc)=1}}f(n)\Psi_d(n)=-\int_{x/4}^x\Delta_f(y;dc,-1)\Psi_d'(y)\,\dd y.
\]
The derivative is supported in $[x/2-1,x-1]$ and satisfies $ \Psi_d'(y)\ll L/x$. Since $(d,P(z))=1$ and $c\mid P(z)$, we have $4\nmid dc$, and $dc$ determines $d$ and $c$. Here $U_0R_0<x^{0.41}$ for large $x$. Replacing the progression sum in \eqref{common:eq:typeI-T1-truncated} by its reduced-residue main term introduces an error which, by Lemma~\ref{common:lem:distribution} with exponent $A+1$, is
\[
 \ll\frac Lx\int_{x/4}^x
       \sum_{\substack{q\le x^{0.41}\\4\nmid q}}
       \max_{(a,q)=1}|\Delta_f(y;q,a)|\,\dd y
 \ll L\sup_{x/4\le y\le x}
       \sum_{\substack{q\le x^{0.41}\\4\nmid q}}
       \max_{(a,q)=1}|\Delta_f(y;q,a)|\ll_A xL^{-A}.
\]
Since $(d,c)=1$, we obtain
\begin{equation}\label{common:eq:typeI-T1-residue-mean}
 T_1=\sum_{\substack{d\le U_0\\(d,P(z))=1}}\frac{\mu(d)}{\varphi(d)}
     \sum_{\substack{c\mid P(z)\\c\le R_0}}\frac{\mu(c)}{\varphi(c)}
     \sum_{\substack{n\ge1\\(n,dc)=1}}f(n)\Psi_d(n)
     +O_A(xL^{-A}).
\end{equation}

For $(d,P(z))=1$ and $d\le U_0^2$, every prime factor of $d$ exceeds $z$, so
\[
 \sum_{p\mid d}\frac1p\le\frac{\log d}{z\log z}
       \ll\frac L{z\log z},\qquad
 \frac d{\varphi(d)}=1+O\!\Bigl(\frac L{z\log z}\Bigr).
\]
Also $\sum_{c\mid P(z)}\varphi(c)^{-1}\ll L$. Since $f(n)\ll_K d(n)$, we have
\[
 \sum_{n\le2x}f(n)\ll_K xL,\qquad
 \sum_{\substack{n\le2x\\p\mid n}}d(n)
       \le2\sum_{v\le2x/p}d(v)
       \ll\frac{xL}p.
\]
Replacing $1/\varphi(d)$ by $1/d$ and removing $(n,d)=1$ from \eqref{common:eq:typeI-T1-residue-mean} contribute at most
\[
 \ll\frac{xL^4}{z\log z}
       \sum_{\substack{d\le U_0\\(d,P(z))=1}}\frac1d
       +xL^3\sum_{\substack{d\le U_0\\(d,P(z))=1}}
                    \frac1d\sum_{p\mid d}\frac1p
 \ll\frac{xL^5}{z\log z}\ll_A xL^{-A}.
\]
The extension of the $c$-sum to all $c\mid P(z)$ contributes
\[
 \ll xL^3\sum_{\substack{c\mid P(z)\\c>R_0}}\frac1{\varphi(c)}
       \ll_A xL^{-A}
\]
by \eqref{common:eq:typeI-sieve-tail}. Interchanging the order of summation and evaluating the $c$-sum, we obtain
\begin{equation}\label{common:eq:typeI-T1-rough-moments}
\begin{aligned}
 T_1={}&\sum_{n\ge1}f(n)w_z(n)\Phi((n+1)/x)
       \sum_{\substack{d\le U_0\\(d,P(z))=1}}
                   \frac{\mu(d)}d\log\frac{n+1}{d}
       +O_A(xL^{-A})\\
 ={}&\sum_{n\ge1}f(n)w_z(n)\Phi((n+1)/x)
        \bigl(M_0\log(n+1)-M_1\bigr)+O_A(xL^{-A}),
\end{aligned}
\end{equation}
where
\[
 M_0=\sum_{\substack{d\le U_0\\(d,P(z))=1}}\frac{\mu(d)}d,
 \qquad
 M_1=\sum_{\substack{d\le U_0\\(d,P(z))=1}}\frac{\mu(d)\log d}d.
\]

We use the identity
\[
 \mu(n)\ind_{(n,P(z))=1}
       =\sum_{\substack{vm=n\\P^+(v)\le z}}\mu(m)
\]
to get
\[
M_0=\sum_{\substack{a\le U_0\\P^+(a)\le z}}\frac1a
\sum_{m\le U_0/a}\frac{\mu(m)}m,\qquad
M_1=\sum_{\substack{a\le U_0\\P^+(a)\le z}}\frac1a
\sum_{m\le U_0/a}\frac{\mu(m)\log(am)}m.
\]
By \cite[(6.18) and Exercise~16(a)]{MV} and partial summation, for $y\ge2$,
\[
 \sum_{m\le y}\frac{\mu(m)}m\ll_A(\log y)^{-A-7},\qquad
 \sum_{m\le y}\frac{\mu(m)\log m}m
       =-1+O_A((\log y)^{-A-5}).
\]
For $a\le\sqrt{U_0}$, these estimates yield the main terms $0$ and $-1$ in the inner sums of $M_0$ and $M_1$, with total error $O_A(L^{-A-4})$. For $a>\sqrt{U_0}$, we use Rankin's trick again and get
\[
 \sum_{\substack{a>\sqrt{U_0}\\P^+(a)\le z}}\frac1a
 \le U_0^{-1/(2\log z)}
       \prod_{p\le z}\bigl(1-p^{-1+1/\log z}\bigr)^{-1}
 \ll\exp\!\Bigl(-\frac{L^{0.76}}{10}+O(L^{0.24})\Bigr)
 \ll_A L^{-A-7}.
\]
Both inner sums are $O(L^2)$ for $a>\sqrt{U_0}$, and $\sum_{P^+(a)\le z}a^{-1}=W(z)^{-1}\ll L$. Thus
\begin{equation}\label{common:eq:typeI-rough-M0-M1}
 M_0\ll_A L^{-A-4},\qquad
 M_1=-W(z)^{-1}+O_A(L^{-A-4}).
\end{equation}
Inserting \eqref{common:eq:typeI-rough-M0-M1} into \eqref{common:eq:typeI-T1-rough-moments}, we conclude
\begin{equation}\label{common:eq:typeI-T1}
 T_1=\frac1{W(z)}\sum_{n\ge1}f(n)w_z(n)\Phi((n+1)/x)
       +O_A(xL^{-A}).
\end{equation}

For $T_2$, the same argument applies with $m=d\ell$. Indeed,
\[
mc\le U_0^2zR_0=x^{0.4}\exp(L^{1/2}+L^{0.24})<x^{0.41}
\qquad(c\le zR_0),
\]
and
\[
\sum_{\substack{d\ell=m\\d,\ell\le U_0}}
|\mu(d)|\Lambda(\ell)\le\sum_{\ell\mid m}\Lambda(\ell)=\log m\ll L.
\]
Thus
\begin{equation}\label{common:eq:typeI-T2}
\begin{aligned}
 T_2={}&\sum_{\substack{d,\ell\le U_0\\(d\ell,P(z))=1}}
       \mu(d)\Lambda(\ell)
       \sum_{\substack{c\mid P(z)\\c\le R_0}}\mu(c)
       \sum_{\substack{n\ge1\\n\equiv-1\,({\rm mod}\,d\ell c)}}
                    f(n)\Phi((n+1)/x)+O_A(xL^{-A})\\
 ={}&M_0\bigg(\sum_{\substack{\ell\le U_0\\(\ell,P(z))=1}}
                         \frac{\Lambda(\ell)}\ell\bigg)
       \sum_{n\ge1}f(n)w_z(n)\Phi((n+1)/x)+O_A(xL^{-A})\\
 \ll{}&xL^{-A},
\end{aligned}
\end{equation}
where the last bound follows from \eqref{common:eq:typeI-rough-M0-M1}. Combining \eqref{common:eq:typeI-T2} with \eqref{common:eq:typeI-T1} completes the proof.
\end{proof}

To evaluate the main term, we separate the factor at $2$ and write
\begin{equation}\label{common:eq:wz-factorization}
 w_z(n)=A_z\ind_{2\mid n}R_z(n),
\end{equation}
where
\[
 A_z=\prod_{3\le p\le z}\frac{p-2}{p-1},\qquad
 R_z(n)=\prod_{\substack{3\le p\le z\\p\mid n}}\frac{p-1}{p-2}.
\]
The function $R_z$ is multiplicative and satisfies $R_z(2n)=R_z(n)$. Moreover,
\[
 \frac{A_z}{W(z)}
 =2\prod_{3\le p\le z}\Bigl(1-\frac1{(p-1)^2}\Bigr)\asymp1.
\]

\subsection{Evaluation of the main term}
Since every auxiliary prime exceeds $\exp(L^{0.1})$, we have
\[
 \sum_{p\in\Pset_\nu}\Bigl(\frac1{p-1}-\frac1p\Bigr)
 \ll\exp(-L^{0.1}).
\]
Consequently,
\begin{equation}\label{mot:eq:mean-W}
 \prod_{\nu=1}^K
 \biggl(\frac1{V_\nu}\sum_{p\in\Pset_\nu}\frac1{p-1}\biggr)
 =1+O_K\bigl(\exp(-L^{0.1})\bigr),
\end{equation}
uniformly in the family index.

\begin{proposition}\label{common:prop:typeI}
Let $K\ge1$ be fixed. Uniformly for $J$ in \eqref{common:eq:family-range} and every $\Phi\in C_c^\infty((0,\infty))$ supported in $[1/2,1]$, we have
\[
 T_1-T_2=\frac{x}{\sqrt L}\biggl(\kappa\Spr
       \int_0^\infty\Phi(t)\,\dd t
       +O_K\bigl(\mathcal N_1(\Phi)L^{-1/2}\bigr)\biggr).
\]
\end{proposition}

\begin{proof}
By Lemma~\ref{common:lem:typeI-model} with $A=1$, it suffices to prove
\begin{equation}\label{mot:eq:meanmodel-asymp}
 \frac1{W(z)}\sum_{n\ge1}b(n)\W(n)w_z(n)\Phi((n+1)/x)
 =\frac{x}{\sqrt L}
 \biggl(\kappa\Spr\int_0^\infty\Phi(t)\,\dd t
       +O_K\bigl(\mathcal N_1(\Phi)L^{-0.76}\bigr)\biggr).
\end{equation}
Since $w_z(n)=0$ for odd $n$, it follows from the identities $b(2n)=b(n)$, $\W(2n)=\W(n)$ and \eqref{common:eq:wz-factorization} that
\begin{equation}\label{mot:eq:even-model}
 \sum_{n\le y}b(n)\W(n)w_z(n)
 =A_z\sum_{m\le y/2}b(m)R_z(m)\W(m).
\end{equation}
Put
\[
 G_z(s)=
 \prod_{\substack{p\le z\\p\equiv1\modf}}
       \Bigl(1+\frac{p^{-s}}{p-2}\Bigr)
 \prod_{\substack{p\le z\\p\equiv3\modf}}
       \Bigl(1+\frac{p^{-2s}}{p-2}\Bigr).
\]
For $\Re s>1$, put $\mathcal B(s)=\sum_{n\ge1}b(n)n^{-s}$. Comparing Euler factors, we find
\[
 \sum_{n\ge1}\frac{b(n)R_z(n)\W(n)}{n^s}
 =\mathcal B(s)\sum_{d\ge1}\frac{h(d)}{d^s},
\]
where
\begin{equation}\label{mot:eq:marked-euler-product}
 \sum_{d\ge1}\frac{h(d)}{d^s}
 =G_z(s)\prod_{\nu=1}^K
 \biggl(\frac1{V_\nu}\sum_{p\in\Pset_\nu}
       \frac{(p-1)p^{-s}}{p-2+p^{-s}}\biggr).
\end{equation}
Indeed, choosing a prime $p$ in a band replaces the factor $1+p^{-s}/(p-2)$ by $(p-1)p^{-s}/(p-2)$. The chosen primes belong to distinct bands, so these factors cancel separately in \eqref{mot:eq:marked-euler-product}. In particular, the right side is a Dirichlet polynomial with nonnegative coefficients. Since $p^{1/\log z}\le e$ for $p\le z$, we obtain
\begin{equation}\label{mot:eq:h-uniform-norms}
 \sum_{d\ge1}\frac{h(d)}{d^{1-1/\log z}}\ll_K1,
 \qquad
 \sum_{d\ge1}\frac{h(d)\log d}{d}\ll_K\log z.
\end{equation}
For the first estimate, the product of the unselected Euler factors is bounded by $\exp(O(\sum_p p^{-2}))$, and each selected prime sum, after division by $V_\nu$, is bounded. The second follows from $(\log d)d^{-1/\log z}\ll\log z$.

By \cite[Theorem~2.1]{BTSD}, with $f=b$, $r=1/2$, $A=2$ and $\sigma=3/4$, we have
\[
 \sum_{n\le v}b(n)
 =\frac{\kappa v}{\sqrt{\log v}}
  +O\Bigl(\frac{v}{(\log v)^{3/2}}\Bigr)
 \qquad(v\ge2).
\]
We apply this to $b*h=bR_z\W$, first for $d\le\sqrt y$. The replacement of $(\log(y/d))^{-1/2}$ by $(\log y)^{-1/2}$ costs $O_K(y\log z/(\log y)^{3/2})$ by \eqref{mot:eq:h-uniform-norms}. The terms with $d>\sqrt y$ are $O_K(y^{1-1/(2\log z)})$, since
\[
 \sum_{d>\sqrt y}\frac{h(d)}d
 \le y^{-1/(2\log z)}
       \sum_d\frac{h(d)}{d^{1-1/\log z}}.
\]
At $s=1$, it follows from \eqref{mot:eq:marked-euler-product} and \eqref{mot:eq:mean-W} that
\[
 \sum_d\frac{h(d)}d
 =G_z(1)\bigl(1+O_K(\exp(-L^{0.1}))\bigr).
\]
Consequently, uniformly for $y\asymp x$,
\begin{equation}\label{mot:eq:convolved-sum}
 \sum_{n\le y}b(n)R_z(n)\W(n)
 =\frac{\kappa y}{\sqrt{\log y}}
       \bigl(G_z(1)+O_K(L^{-0.76})\bigr).
\end{equation}
Here $\log z=L^{0.24}$ and $y^{-1/(2\log z)}\ll(\log y)^{-2}$.

The Euler factor of $A_zG_z(1)/(2W(z))$ is $1$ at $p=2$ and at $p\equiv1\modf$. For $p\le z$ with $p\equiv3\modf$, it is
\[
 \frac{p-2}{p-1}\Bigl(1+\frac1{p^2(p-2)}\Bigr)
       \Bigl(1-\frac1p\Bigr)^{-1}
 =1-\frac1{p(p-1)}.
\]
Thus
\begin{equation}\label{mot:eq:euler-identity}
 \frac{A_zG_z(1)}{2W(z)}
 =\prod_{\substack{p\le z\\p\equiv3\modf}}
       \Bigl(1-\frac1{p(p-1)}\Bigr)
 =\Spr+O(z^{-1}).
\end{equation}
Combining \eqref{mot:eq:even-model}, \eqref{mot:eq:convolved-sum} and \eqref{mot:eq:euler-identity}, we obtain
\[
 \frac1{W(z)}\sum_{n\le y}b(n)\W(n)w_z(n)
 =\frac{y}{\sqrt{\log y}}
       \bigl(\kappa\Spr+O_K(L^{-0.76})\bigr).
\]
Partial summation proves \eqref{mot:eq:meanmodel-asymp}, with the dependence on $\Phi$ bounded by $\mathcal N_1(\Phi)$.
\end{proof}

\section{A mean square estimate and Type~II sums}\label{core:sec:quadratic}
\subsection{The mean square estimate}
We retain the notation of Section~\ref{sec:common}, with $K$ fixed as $x\to\infty$. For the mean square estimate, let $\Pset_1,\ldots,\Pset_K$ be any sets of primes satisfying
\begin{equation}\label{core:eq:restricted-band-hyp}
\begin{gathered}
 \Pset_\nu\subseteq
 [\exp(L^{\vartheta_\nu}),\exp(2L^{\vartheta_\nu})],\qquad
 \Pset_\nu\cap\Pset_\mu=\varnothing\quad(\nu\ne\mu),\\
 0.1<\vartheta_1<\cdots<\vartheta_K<0.2,\qquad
 \vartheta_{\nu+1}-\vartheta_\nu\ge\frac1{20K}\quad(1\le\nu<K),\qquad
 \frac14\le V_\nu\le\frac12.
\end{gathered}
\end{equation}
The quantities in \eqref{common:eq:weight} and \eqref{common:eq:products} now refer to these sets. Put $p_0=\min\Pset_*$, so $p_0\ge\exp(L^{0.1})$. The application sets in \eqref{common:eq:application-bands} satisfy \eqref{core:eq:restricted-band-hyp} for sufficiently large $x$. The proof of Theorem~\ref{core:quadratic-theorem} uses only \eqref{core:eq:restricted-band-hyp}.

Fix $\delta,C>0$, and suppose $X=MN\asymp x$ and $M,N\ge x^\delta$. Let $(\alpha_m)$ and $(\beta_n)$ be complex sequences satisfying
\begin{equation}\label{core:eq:coefficient-bounds}
\begin{gathered}
 \supp\alpha\subseteq\{m\in[M,2M]:P^-(m)>z\},\qquad
 \supp\beta\subseteq\{n\in[N,2N]:P^-(n)>z\},\\
 |\alpha_m|,|\beta_n|\le L^C .
\end{gathered}
\end{equation}
Using the cutoff $\eta$ fixed in Section~\ref{sec:common}, define, for $Y>0$, the mean square
\begin{equation}\label{core:eq:typeII-Q}
 Q_Y=\sum_{h\ge1}\biggl|\mathcal V(\Pset)^{-1}
   \sum_{\substack{a\in\mathcal A(\Pset),\,m,n\\mn-1=ah}}
       \eta(a/Y)\alpha_m\beta_n\biggr|^2.
\end{equation}

\begin{samepage}
\begin{theorem}\label{core:quadratic-theorem}
For every $R>0$, there is an integer $K_0=K_0(\delta,C,R,\eta)$ with the following property. Let $K\ge K_0$ be fixed, and let the exponents and prime sets satisfy \eqref{core:eq:restricted-band-hyp}. If
\[
 \alpha_m=\mu(m)m^{iv}\ind_{m\in I}\ind_{P^-(m)>z},
 \quad I\subseteq[M,2M]\text{ an interval},\quad v\in\mathbb R,\ |v|\le L^C,
\]
then, for $x\ge x_0(\delta,C,R,K,\eta)$, where $x_0$ may be ineffective,
\[
 Q_Y\ll_{\delta,C,R,K,\eta}XYL^{-R}
 \qquad(Y>0).
\]
\end{theorem}
\end{samepage}

The proof of Theorem~\ref{core:quadratic-theorem} occupies Sections~\ref{core:sec:minor-reduction}--\ref{core:sec:typeII-major}. We first deduce the estimates for the Type~II sums in the proof of Theorem~\ref{mot:thm:main}.

\subsection{The Type~II sums}\label{common:sec:typeII}
In this subsection the prime sets are given by \eqref{common:eq:application-bands}.

\begin{proposition}\label{common:prop:typeII}
For every fixed $\delta,C,A>0$, there is an integer $K_0=K_0(\delta,C,A)$ with the following property. Let $K\ge K_0$ be fixed, $M,N\ge x^\delta$ and $MN\asymp x$. Suppose that
\[
 \alpha_m=\mu(m)m^{iv}\ind_{m\in I}\ind_{P^-(m)>z},
 \qquad I\subseteq[M,2M]\text{ an interval},
 \qquad v\in\mathbb R,\quad |v|\le L^C,
\]
and that $\beta$ is supported on the $z$-rough integers in $[N,2N]$ with $|\beta_n|\le L^C$. Then
\[
 \sum_{m,n}\alpha_m\beta_n b(mn-1)\W(mn-1)
 \ll_{\delta,C,A,K}xL^{-A}.
\]
\end{proposition}
\begin{proof}
By \eqref{common:eq:weight}, for $u\ge1$ we have
\[
 b(u)\W(u)=\frac1{\mathcal V(\Pset)}
 \sum_{\substack{a\in\mathcal A(\Pset)\\a\mid u}}b(u).
\]
Write $a=p_1\cdots p_K$, with $p_\nu\in\Pset_\nu$, and $u=ah$. Since every auxiliary prime is congruent to $1$ modulo $4$, \eqref{mot:eq:b-def} implies $b(ah)=b(h)$ without a coprimality condition.

Using $\sum_Y\eta(a/Y)=1$, with $Y$ running over powers of $2$, we divide the $a$-sum into dyadic intervals. Since $a$ is a product of $K$ auxiliary primes, only $O_K(L)$ values of $Y$ occur, and for these $1\le Y\le\exp(2KL^{0.2})$ for sufficiently large $x$. On the support of $\eta(a/Y)$, we have $h\ll x/Y$. Since $0\le b(h)\le1$,
\[
 \sum_{h\ll x/Y}b(h)^2\ll\frac{x}{Y}.
\]
It follows from Cauchy's inequality and \eqref{core:eq:typeII-Q} that
\[
 \biggl|\sum_{m,n}\alpha_m\beta_n b(mn-1)\W(mn-1)\biggr|
 \ll\sum_Y\Bigl(\frac{x}{Y}Q_Y\Bigr)^{1/2}.
\]
By Theorem~\ref{core:quadratic-theorem}, with $R>2A+4$,
\[
 \sum_Y\Bigl(\frac{x}{Y}Q_Y\Bigr)^{1/2}
 \ll_{\delta,C,R,K}xL^{1-R/2}
 \ll xL^{-A}.
\]
This proves the proposition.
\end{proof}

\begin{corollary}\label{common:cor:T3}
For every $A>0$, there is an integer $K_0=K_0(A)$ such that, for each fixed $K\ge K_0$ and every $\Phi\in C_c^\infty((0,\infty))$ supported in $[1/2,1]$, the sum $T_3$ defined by \eqref{common:eq:T3} satisfies
\[
 T_3\ll_{A,K}\mathcal N_{\lceil A\rceil+7}(\Phi)xL^{-A}.
\]
\end{corollary}
\begin{proof}
Partition $d,t>U_0$ into $(M,2M]$ and $(N,2N]$, where $M=2^iU_0$, $N=2^jU_0$, and $i,j\ge0$. Only $O(L^2)$ pairs of intervals meet the support of $\Phi(dt/x)$. For each such pair,
\[
 x/8\le MN\le x,\qquad M,N\ge x^{1/5}.
\]
Mellin inversion gives
\[
 \Phi(dt/x)=\frac1{2\pi}\int_{\mathbb R}
 \widehat\Phi(iv)\Bigl(\frac{dt}{x}\Bigr)^{iv}\dd v,
 \qquad
 \widehat\Phi(iv)=\int_0^\infty\Phi(y)y^{-1-iv}\dd y.
\]
Integration by parts in the variable $\log y$ shows that $\widehat\Phi(iv)\ll_B\mathcal N_B(\Phi)(1+|v|)^{-B}$ for every integer $B\ge2$. We have $b\W\ll_K d$ and
\[
\begin{aligned}
 \sum_{2\le dt\le4x}d(dt-1)
 &\le\sum_{1\le u\le4x}d(u+1)d(u)\\
 &\le\biggl(\sum_{1\le u\le4x}d(u+1)^2\biggr)^{1/2}
       \biggl(\sum_{1\le u\le4x}d(u)^2\biggr)^{1/2}
 \ll xL^3,
\end{aligned}
\]
The last estimate follows from Lemma~\ref{core:lem:complete-divisors} with $q=1$. Together with $c_{U_0}(t)\le\log t\ll L$, this shows that the integral over $|v|>L$, summed over these pairs of intervals, contributes $O_K(\mathcal N_B(\Phi)xL^{5-B})$. Taking $B=\lceil A\rceil+7$ bounds this by $O_{A,K}(\mathcal N_B(\Phi)xL^{-A})$.

For $|v|\le L$, take
\[
 \alpha_d=\mu(d)d^{iv}\ind_{M<d\le2M}\ind_{P^-(d)>z},
 \qquad
 \beta_t=c_{U_0}(t)t^{iv}\ind_{N<t\le2N}\ind_{P^-(t)>z}.
\]
Here $2N\le2x^{4/5}<x$ for large $x$, so $|\beta_t|\le L$. Proposition~\ref{common:prop:typeII} applies with $\delta=1/5$, $C=1$, and saving $L^{-A-3}$. Summing over the $O(L^2)$ pairs and using $\int_{\mathbb R}|\widehat\Phi(iv)|\dd v\ll\mathcal N_2(\Phi)$ proves the assertion.
\end{proof}

Combining \eqref{common:eq:vaughan}, Lemma~\ref{common:lem:typeI-model} and Corollary~\ref{common:cor:T3}, we obtain, for every fixed $A>0$ and all sufficiently large fixed $K$,
\begin{equation}\label{common:eq:weighted-model}
 \Theta_{b,\W}(x,\Phi)
 =\frac1{W(z)}\sum_{n\ge1}b(n)\W(n)w_z(n)\Phi((n+1)/x)
       +O_{A,K}(\mathcal N_{\lceil A\rceil+7}(\Phi)xL^{-A}).
\end{equation}

\section{Reduction of the minor arc estimate}\label{core:sec:minor-reduction}

We return to the prime sets satisfying \eqref{core:eq:restricted-band-hyp}. The coefficients are arbitrary complex sequences satisfying \eqref{core:eq:coefficient-bounds}. Here $\delta,C>0$ and $K$ are fixed as $x\to\infty$, $X=MN\asymp x$, $M,N\ge x^\delta$, and $Y>0$. We retain the notation of Section~\ref{core:sec:quadratic}, including $p_0=\min\Pset_*$ and $Q_Y$ from \eqref{core:eq:typeII-Q}. In this section and Section~\ref{core:sec:typeII-minor}, all estimates are uniform in the exponents and prime sets satisfying \eqref{core:eq:restricted-band-hyp}.

We first separate the major and minor arcs. For the latter, we expand a power of an operator as a sum over paths. Section~\ref{core:sec:typeII-minor} bounds this sum, and Section~\ref{core:sec:typeII-major} estimates the contribution from the major arcs. The two bounds prove Theorem~\ref{core:quadratic-theorem}.

\subsection{The major and minor arcs}

For the minor arc argument we may replace $\delta$ by $\min(\delta,1/4)$. We may also assume $\mathcal A(\Pset)\cap[Y,4Y]\ne\varnothing$, since otherwise $Q_Y=0$. Then $p_0/4\le Y\le\exp(2KL^{0.2})$.

In the expansion of \eqref{core:eq:typeII-Q}, the relations $mn-1=ah$ and $rs-1=ch$ imply
\begin{equation}\label{core:eq:common-h-determinant}
 cmn-ars=c-a.
\end{equation}
Conversely, for $(a,c)=1$ the identity $c(mn-1)=a(rs-1)$ implies $a\mid mn-1$ and $c\mid rs-1$, with equal quotients. These quotients are positive for sufficiently large $x$, since $mn,rs\gg x^{2\delta}$. We detect \eqref{core:eq:common-h-determinant} by Fourier orthogonality, as in the dispersion method \cite[Section~19.2]{FI}.

Fix once and for all an even real $\psi\in C_c^\infty(\mathbb R)$ that is one on $[-4,4]$ and supported on $[-5,5]$. Write $\mathbb T=\mathbb R/\mathbb Z$ and let $\|\cdot\|$ denote distance to the nearest integer. For $A_0>0$ put
\begin{equation}\label{core:eq:major-arcs-def}
  \mathfrak M=\bigcup_{\substack{q\le L^{A_0}\\0\le u<q,\ (u,q)=1}}
 \biggl\{\theta\in\mathbb R/\mathbb Z:
   \Bigl\|\theta-\frac uq\Bigr\|\le\frac{2L^{A_0}}Y\biggr\}.
\end{equation}
Define the kernel and the corresponding contribution by
\begin{align}
 H_{\mathfrak M}(t,a,c)
   &=\psi(t/Y)\int_{\mathfrak M}e\bigl(\theta(t-c+a)\bigr)\,\dd\theta,\notag
\\
 Q_Y(\mathfrak M)
   &=\mathcal V(\Pset)^{-2}
     \sum_{\substack{a,c\in\mathcal A(\Pset)\\m,n,r,s}}
       \eta(a/Y)\eta(c/Y)\alpha_m\overline{\alpha_r}
       \beta_n\overline{\beta_s}
       H_{\mathfrak M}(cmn-ars,a,c).                 \label{core:eq:major-Q}
\end{align}

The sum defining $Q_Y(\mathfrak M)$ includes all pairs $a,c$. When $mn-1=ah$ and $rs-1=ch$, $t=cmn-ars=c-a$ satisfies $|t|\le3Y$, so $\psi(t/Y)=1$. Put $\mathfrak m=\mathbb T\setminus\mathfrak M$ and $Q_Y(\mathfrak m)=Q_Y-Q_Y(\mathfrak M)$. This difference is the sum in \eqref{core:eq:major-Q} with $H_{\mathfrak M}(t,a,c)$ replaced by
\begin{equation}\label{core:eq:minor-bracket}
 \mathcal B_Y(a,c;m,n,r,s)
 =\psi(t/Y)\biggl(\sum_{h\ge1}\ind_{mn-1=ah}\ind_{rs-1=ch}
       -\int_{\mathfrak M} e\bigl(\theta(t-c+a)\bigr)\dd\theta\biggr).
\end{equation}
When $(a,c)=1$, \eqref{core:eq:common-h-determinant} and Fourier orthogonality identify this bracket with $\psi(t/Y)\int_{\mathfrak m}e(\theta(t-c+a))\dd\theta$. For $(a,c)>1$ we retain the sum over the positive integer $h$ in \eqref{core:eq:minor-bracket}.

\begin{proposition}\label{core:lem:restricted-determinant}
For fixed $0<\delta<1/2$ and $C,R>0$, there exist $A_0>0$ and $K_0\ge1$, depending only on $\delta,C,R,\eta$, such that for each fixed $K\ge K_0$,
\begin{equation*}
 |Q_Y-Q_Y(\mathfrak M)|\ll_{\delta,C,R,K,\eta}XYL^{-R}\qquad(Y>0).
\end{equation*}
\end{proposition}

The rest of this section reduces Proposition~\ref{core:lem:restricted-determinant} to \eqref{core:eq:normalized-high-moment}, apart from exceptional terms. Section~\ref{core:sec:typeII-minor} proves this bound and estimates the exceptional terms.

\subsection{Inserting sums over auxiliary primes}

Put $J_0=\lfloor L^{0.01}\rfloor$ and $m_0=J_0+1$. For each $\nu$, we adjoin $J_0$ shared prime labels from $\Pset_\nu$ to the two endpoints, weighting each choice of $p$ by $(V_\nu p)^{-1}$ as in \cite[\S2.2]{Helfgott}. Thus
\begin{equation}\label{core:eq:harmonic-switching}
  1=\sum_{\substack{p_{\nu,h}\in\Pset_\nu\\
                  1\leq\nu\leq K,\ 1\leq h\leq J_0}}
  \prod_{\nu=1}^K\prod_{h=1}^{J_0}\frac1{V_\nu p_{\nu,h}}.
\end{equation}
For each choice of the shared primes, put $D=\prod_{\nu=1}^K\prod_{h=1}^{J_0}p_{\nu,h}$. We divide the products $D$ into $O_K(L)$ dyadic intervals $[d_0,2d_0)$. Multiplying the sum for $Q_Y(\mathfrak m)$ by \eqref{core:eq:harmonic-switching}, we obtain
\begin{align}
 Q_Y(\mathfrak m)
 ={}&
 \sum_{d_0\ {\rm dyadic}}
 \sum_{\substack{p_{\nu,h}\in\Pset_\nu\\d_0\leq D<2d_0}}
 \biggl(\prod_{\nu=1}^K\prod_{h=1}^{J_0}\frac1{V_\nu p_{\nu,h}}\biggr)
 \frac1{\mathcal V(\Pset)^2}
 \notag\\
 &\quad\times
 \sum_{\substack{a_\nu,c_\nu\in\Pset_\nu\\1\leq\nu\leq K}}
 \eta(a/Y)\eta(c/Y)
 \sum_{m,n,r,s}
 \alpha_m\beta_n\overline{\alpha_r}\,\overline{\beta_s}
 \mathcal B_Y(a,c;m,n,r,s),
\label{core:eq:switched-tuple-sum}
\end{align}
where
\[
  a=\prod_{\nu=1}^{K}a_\nu,\qquad c=\prod_{\nu=1}^{K}c_\nu.
\]
Fix $d_0$ and put
\begin{equation}\label{core:eq:endpoint-coordinates}
 P=\begin{pmatrix}Dcm\\s\end{pmatrix},\qquad
 Q=\begin{pmatrix}Dar\\n\end{pmatrix}.
\end{equation}
Then
\[
 \det(P,Q)=D(cmn-ars).
\]
Put
\begin{equation}\label{core:eq:scales}
 H=d_0Y,\qquad U=HM,\qquad V=N .
\end{equation}
On the support of the sum we have $P,Q\in[U,16U]\times[V,2V]$ and $|\det(P,Q)|\le10H$. Suppose $\mu^2(Dac)=1$ and the coefficients and bracket are nonzero. Then $t=cmn-ars\ne0$, since $cmn=ars$ would imply $a\mid mn$. Indeed, since $\mu^2(Dac)=1$, we have $(a,c)=1$, and all auxiliary primes are less than $z$, hence coprime to $m,n,r,s$. The endpoint gcds are $(m,s)$ and $(r,n)$. Each divides $t$ and has all its prime factors greater than $z$. Since $0<|t|\le5Y<z$ for fixed $K$ and large $x$, both points are primitive.

For each $\nu$, attach an ordered tuple of $m_0=J_0+1$ primes from $\Pset_\nu$ to each endpoint. Along an edge the two tuples share $J_0$ entries, so, after reordering, they take the form
\begin{equation}\label{core:eq:edge-labels}
 \mathbf L_\nu=(p_{\nu,1},\ldots,p_{\nu,J_0},c_\nu),
 \qquad
 \mathbf L'_\nu=(p_{\nu,1},\ldots,p_{\nu,J_0},a_\nu).
\end{equation}
The $p_{\nu,h}$ are the shared labels, $c_\nu$ the leaving label, and $a_\nu$ the entering label. With $\mathbf L=(\mathbf L_1,\ldots,\mathbf L_K)$ and $\Pi(\mathbf L)=\prod_\nu\prod_{p\in\mathbf L_\nu}p$, we have
\begin{equation*}
 \Pi(\mathbf L)=Dc,\qquad \Pi(\mathbf L')=Da.
\end{equation*}

\begin{definition}
A \emph{state} consists of a primitive point $P\in[U,16U]\times[V,2V]$ and, for each $\nu$, a tuple $\mathbf L_\nu\in\Pset_\nu^{m_0}$ whose entries are distinct prime divisors of $P_1$. A directed edge $(P,\mathbf L)\to(Q,\mathbf L')$ is \emph{admissible} if, for every $\nu$, the two tuples share exactly $J_0$ entries, and the product $D$ of all shared entries lies in $[d_0,2d_0)$. Writing $c_\nu$ and $a_\nu$ for the remaining entries as in \eqref{core:eq:edge-labels}, we then have $\mu^2(Dac)=1$.

\end{definition}

An admissible edge determines the shared product $D$ and the four arithmetic variables through
\[
m=\frac{P_1}{\Pi(\mathbf L)},\qquad
s=P_2,\qquad
r=\frac{Q_1}{\Pi(\mathbf L')},\qquad
n=Q_2.
\]

Let $\mathscr S_{d_0}$ be the set of all states, including isolated states. Assign to each state the mass $w_{\rm st}=\prod_\nu V_\nu^{-m_0}$ and put
\[
 \mathcal H=\ell^2(\mathscr S_{d_0},w_{\rm st}).
\]

\subsection{An operator on labelled points}\label{core:sec:diophantine-separation}
For \(D\in[d_0,2d_0)\), put $t=\det(P,Q)/D$ and define the kernel
\begin{equation}\label{core:eq:coarse-kernel}
 \mathcal K(P,Q;D,a,c)
 =
 \frac{d_0}{D}\eta(c/Y)\eta(a/Y)\psi(t/Y)
 \biggl(
  \ind_{t=c-a}
  -
  \int_{\Mmajor}e\bigl(\theta(t-c+a)\bigr)\dd\theta
 \biggr).
\end{equation}
When $(a,c)=1$, the first term imposes equality of the two quotients, as in \eqref{core:eq:common-h-determinant}. The factors associated with auxiliary primes involve only the selected labels. An auxiliary prime divisor of an endpoint coordinate that is not selected receives no additional factor.

Each unordered state has $(m_0!)^K$ ordered representatives. For two unordered tuples drawn from the same set of primes and sharing $J_0$ entries, there are $(m_0!)^2$ ways to order the two tuples and $J_0!$ ways to order their shared entries in \eqref{core:eq:switched-tuple-sum}. We therefore put $\kappa_{J_0}=\frac{J_0!}{(m_0!)^2}=\frac1{(J_0+1)(J_0+1)!}$. For an admissible edge $(x,y)$, define
\[
 \kappa(x,y)=\frac{\kappa_{J_0}^K}{\mathcal V(\Pset)}\,
                 \mathcal K(P,Q;D,a,c),
\]
and set $\kappa(x,y)=0$ on all other pairs. The associated operator and its form are
\[
 (\mathsf K f)(x)=\sum_y\kappa(x,y)f(y),\qquad
 \langle F,\mathsf K G\rangle
     =w_{\rm st}\sum_{x,y}\overline{F(x)}\,\kappa(x,y)G(y).
\]
Thus the kernel relative to the state measure is $\kappa(x,y)/w_{\rm st}$. For unordered states, sum the coefficients over the ordered representatives of the target state \cite[\S2.3]{BH}. There are $(m_0!)^{K}$ such representatives, so a pair of unordered states contributes $(m_0!)^{2K} w_{\rm st}\kappa(x,y)$ to the form on their indicators.

For a primitive endpoint $P=(u,v)^t$, put
\[
 \tau(P)=\frac{\xi_P}{u}\pmod1,\qquad
 v\xi_P\equiv-1\pmod u,\quad 0\le\xi_P<u.
\]
For a tuple $\mathbf L$, let $\mathcal D$ consist of the products obtained by omitting one entry from each tuple. Thus $|\mathcal D|\le m_0^K$.

For $1\le s\le K$ set
\[
 \mathcal N_s=\biggl\{D\prod_{\nu<s}p_\nu:
       D\in\mathcal D,\ p_\nu\in\Pset_\nu\biggr\},\qquad
 \mathcal F_s=\biggl\{\prod_{\nu>s}q_\nu:q_\nu\in\Pset_\nu\biggr\}.
\]
Repeated representations of the same integer in $\mathcal N_s$ are counted once. Assign to $F=\prod_{\nu>s}q_\nu$ the weight $w_s(F)=\prod_{\nu>s}(V_\nu q_\nu)^{-1}$, so that $\sum_{F\in\mathcal F_s}w_s(F)=1$. Empty products are $1$. For $F\in\mathcal F_s$ let $\mathcal E_s(F)$ be the set of $r$ for which two distinct $n,n'\in\mathcal N_s$ satisfy
\[
 \|(n-n')Fr\|\le 100e^{-L^{\vartheta_s}}.
\]
We call $r$ \emph{separated} for $\mathbf L$ if
\begin{enumerate}
\item $\|lDr\|>Y^{-0.7}$ for all $D\in\mathcal D$ and integers $1\le l\le Y^{0.2}$;
\item for each $1\le s\le K$ one has
\begin{equation}\label{core:eq:spacing-weighted-exception}
 W_s(r):=\sum_{F\in\mathcal F_s}w_s(F)\ind_{r\in\mathcal E_s(F)}
      \le e^{-L^{\vartheta_s}/4}.
\end{equation}
\end{enumerate}
Let $\ind_{\rm sep}(r,\mathbf L)$ denote the indicator of these two conditions. The first condition will allow us to choose a common basis for the determinant estimate. The second condition ensures that, outside a set of $F$ of weight at most $e^{-L^{\vartheta_s}/4}$, the values $nFr\bmod1$ are then separated as $n$ varies. At an integral point put
\begin{equation}\label{core:eq:separated-state-indicator}
 \ind_{\rm sep}(P,\mathbf L)=\ind_{\rm sep}(\tau(P),\mathbf L).
\end{equation}

Let $\mathsf G$ denote multiplication by $\ind_{\rm sep}(P,\mathbf L)$, and put
\begin{equation}\label{core:eq:A-factorization}
 \mathsf A=\mathsf G\mathsf K\mathsf G.
\end{equation}
Reversing an edge interchanges $a,c$ and sends $t$ to $-t$. Since $\psi$ is real and even and $\mathfrak M=-\mathfrak M$, the kernel in \eqref{core:eq:coarse-kernel} is real and symmetric. The state mass is constant, so $\mathsf K$ is self-adjoint. As $\mathsf G$ is multiplication by a real indicator, $\mathsf A$ is also self-adjoint.

For a state \((P,\mathbf L)\), define
\begin{equation}\label{core:eq:endpoint-vectors}
  F(P,\mathbf L)
     =\overline{\alpha_{P_1/\Pi(\mathbf L)}}\,\beta_{P_2},
\end{equation}
with coefficients extended by zero. For $P,Q$ as in \eqref{core:eq:endpoint-coordinates},
\[
 \overline{F(P,\mathbf L)}\,F(Q,\mathbf L')
 =\alpha_m\beta_n\overline{\alpha_r}\,\overline{\beta_s}.
\]
If $F(P,\mathbf L)\ne0$, then $P_1/\Pi(\mathbf L)$ is $z$-rough. The entries of $\mathbf L$ therefore account for all auxiliary primes dividing $P_1$, each to exponent one. Every tuple above $P$ is then a reordering of $\mathbf L$. Thus $F$ takes a constant value $f(P)$ on each fibre, with $f(P)=0$ on an empty fibre. All fibres, including those where $f(P)=0$, enter the powers of $\mathsf A$.

\begin{proposition}\label{core:prop:operator-identity}
We have
\begin{equation*}
 Q_Y-Q_Y(\Mmajor)
 =
 \sum_{d_0\ {\rm dyadic}}\frac1{d_0}\ip{F}{\mathsf A F}
 +E_1+E_2,
\end{equation*}
where $E_1$ and $E_2$ are the sums in \eqref{core:eq:switched-tuple-sum} over $\mu^2(Dac)=1$ with at least one endpoint not separated, and over $\mu^2(Dac)=0$, respectively.
\end{proposition}
\begin{proof}
For $\mu^2(Dac)=1$ every contributing pair is primitive, and \eqref{core:eq:common-h-determinant} identifies the two terms of \eqref{core:eq:coarse-kernel} with those of \eqref{core:eq:minor-bracket}. The two factors $\mathsf G$ select the pairs whose endpoints both satisfy the separation conditions.

Fix the sets of labels in the two endpoint tuples. Their $(m_0!)^{2K}$ orderings supply $(m_0!)^{2K}\kappa_{J_0}^K=(J_0!)^K$, exactly the number of ordered shared tuples in \eqref{core:eq:switched-tuple-sum}. For each such tuple the remaining scalar is
\[
 \frac{w_{\rm st}}{d_0}\,\frac{1}{\mathcal V(\Pset)}\,\frac{d_0}{D}
   =\frac1D\prod_{\nu=1}^K V_\nu^{-(J_0+2)}.
\]
This is precisely the weight of this choice in \eqref{core:eq:switched-tuple-sum}.
\end{proof}

\subsection{Reduction to a moment estimate}

Take $T=2\lceil L^{1/2}/2\rceil$. All estimates for powers of the operator below hold for even $T\asymp L^{1/2}$ with fixed implied constants. In each dyadic range for $D$, \eqref{core:eq:scales} and $\log d_0\le2J_0\sum_\nu L^{\vartheta_\nu}$ imply $H=\exp(O_K(L^{0.21}))$ and $x^\delta\le U,V\le x^{O(1)}$. For each primitive $P$, let $u_P$ be the indicator of its full tuple fibre in $\mathcal H$, interpreted as zero when the fibre is empty. For every prescribed $E>0$ we shall prove
\begin{equation}\label{core:eq:normalized-high-moment}
 \frac1{UV}
 \sum_{\substack{P\in[U,16U]\times[V,2V]\\(P_1,P_2)=1}}
       \ip{u_P}{\mathsf A^T u_P}
 \ll L^{-ET}
\end{equation}
with the parameters of Proposition~\ref{core:prop:path-bound}. The usual trace sums walks returning to the same state \cite[\S3.2]{Helfgott}. Here the endpoints have the same position $P$, and their tuples are summed independently.

\begin{lemma}\label{core:lem:spectral-interpolation}
Let $\mathsf A$ act on a finite-dimensional Hilbert space, let $T\ge2$ be even, and put $\mathsf B=(\mathsf A\mathsf A^*)^{T/2}$. For a finite family $(u_i)$ write $b_{ij}=\langle u_i,\mathsf B u_j\rangle$. Suppose each row of $(b_{ij})$ has at most $D_0$ nonzero entries. If $F=\sum_i f_i u_i\ne0$, then for any $G$,
\begin{equation*}
 |\langle F,\mathsf A G\rangle|
 \le \|G\|\,\|F\|^{1-2/T}
 \biggl(D_0\|f\|_\infty^2\sum_i b_{ii}\biggr)^{1/T}.
\end{equation*}
\end{lemma}
\begin{proof}
The spectral measure of $\mathsf A\mathsf A^*$ at $F/\|F\|$ is a probability measure on $[0,\infty)$. Jensen's inequality for $t\mapsto t^{T/2}$ gives the H\"older--McCarthy bound \cite[Theorem~1.4(i)]{FMPSeo}
\[
 \|\mathsf A^*F\|
 \le \|F\|^{1-2/T}\langle F,\mathsf B F\rangle^{1/T}.
\]
The matrix $(b_{ij})$ is the Gram matrix of $(\mathsf B^{1/2}u_i)$, so $2|b_{ij}|\le b_{ii}+b_{jj}$. Its support is symmetric. Hence
\[
 \langle F,\mathsf B F\rangle
 \le \frac{\|f\|_\infty^2}{2}\sum_{b_{ij}\ne0}(b_{ii}+b_{jj})
 \le D_0\|f\|_\infty^2\sum_i b_{ii}.
\]
Together with $|\langle F,\mathsf A G\rangle|\le\|\mathsf A^*F\|\,\|G\|$, this proves the assertion.
\end{proof}

To apply Lemma~\ref{core:lem:spectral-interpolation}, write $\sigma(P)=P_2/P_1$. For distinct primitive positions, Farey spacing gives
\begin{equation}\label{core:eq:position-slope}
 |\sigma(P)-\sigma(Q)|=\frac{|\det(P,Q)|}{P_1Q_1}
                  \ge\frac1{(16U)^2}.
\end{equation}
See \cite[proof of Corollary~1]{MV73}. On an edge the same identity and $|\det(P,Q)|\le10H$ imply $|\sigma(P)-\sigma(Q)|\le10H/U^2$. Thus at most $O(1+sH)$ primitive positions can be reached from $P$ in $s$ steps. This use of slope to restrict a determinant is also made in \cite[\S12, following (12.6) and (12.9)]{FI98}.

\begin{proposition}\label{core:prop:conditional-comparison}
Suppose \eqref{core:eq:normalized-high-moment} holds in every dyadic range for $D$ with exponent $E>0$. Then
\[
 Q_Y-Q_Y(\Mmajor)=E_1+E_2+O\!\bigl(XYL^{-E+4C+1}\bigr).
\]
\end{proposition}

\begin{proof}
By \eqref{core:eq:endpoint-vectors} we have $F=\sum_P f(P)u_P$ with $\|f\|_\infty\le L^{2C}$. For a fixed ordered tuple $\mathbf L$, at most $16U/\Pi(\mathbf L)$ positive multiples of $\Pi(\mathbf L)$ lie in the first-coordinate range, and there are at most $2V$ choices for $P_2$. Dropping primitivity and summing with the state mass, we obtain
\[
 \|F\|^2
 \le 32 UV L^{4C}
 \sum_{\mathbf L}\prod_{\nu=1}^K\prod_{p\in\mathbf L_\nu}
                         \frac1{V_\nu p}
 \le 32 UV L^{4C}.
\]
The last sum is at most one since the weights $(V_\nu p)^{-1}$ sum to one on each set of primes. Apply Lemma~\ref{core:lem:spectral-interpolation} with the fibre indicators $u_P$, taking both input vectors to be $F$. Then \eqref{core:eq:normalized-high-moment} bounds the sum \(\sum_i b_{ii}\) by $UVL^{-ET}$, and the preceding slope bound implies $D_0\ll1+TH$. Hence
\[
 |\ip{F}{\mathsf A F}|
 \ll UV L^{-E+4C}(1+TH)^{1/T}\ll UV L^{-E+4C},
\]
as $\log H=o(T)$. The assertion follows from Proposition~\ref{core:prop:operator-identity}, since $UV/d_0=XY$ in each of the $O_K(L)$ ranges.
\end{proof}

\subsection{Expanding powers of the operator}

With the even parameter $T$ as in the preceding subsection, take for each primitive initial position $P_0=(u_0,v_0)^t$ the canonical frame

\[
 g=g_{P_0}:=\begin{pmatrix}
 u_0&\xi_{P_0}\\v_0&(1+v_0\xi_{P_0})/u_0
 \end{pmatrix}\in\operatorname{SL}_2(\mathbb Z),
\]
where $0\le\xi_{P_0}<u_0$ and $v_0\xi_{P_0}\equiv-1\bmod u_0$. Then $ge_1=P_0$, where $e_1=(1,0)^t$.

For a walk of length $T$ returning to $P_0$, write its positions uniquely as
\[
  P_j=g\mathbf z_j,\qquad
  \mathbf z_0=\mathbf z_T=e_1,
\]
where the \(\mathbf z_j\) are primitive integer vectors. Every contributing walk satisfies
\begin{equation*}
  (P_j)_1\asymp U,\qquad (P_j)_2\asymp V,\qquad
  |\det(P_j,P_{j+1})|\ll H.
\end{equation*}

For a primitive \(\mathbf z\), choose \(\mathbf w_{\mathbf z}\in\mathbb Z^2\) with \(\det(\mathbf z,\mathbf w_{\mathbf z})=1\). Whenever $(g\mathbf z)_1\ne0$, put
\begin{equation*}
 \alpha_g(\mathbf z):=
 \frac{(g\mathbf w_{\mathbf z})_1}{(g\mathbf z)_1}\pmod1.
\end{equation*}
The choice of complement affects the quotient by an integer. In an integral frame, $g\mathbf w_{\mathbf z}$ completes $P=g\mathbf z$ to a basis of determinant one, so $\alpha_g(\mathbf z)=\tau(P)$. For $p\in\Pset_*$ and $g\in\operatorname{SL}_2(\mathbb Z)$, write $\ell_p=\ell_p(g)$ for the kernel of the first row of $g\bmod p$. This line belongs to $\mathbb P^1(\mathbb F_p)$, the set of $p+1$ lines in $\mathbb F_p^2$.

If \([\mathbf z]_p\) denotes the line spanned by a primitive vector \(\mathbf z\) modulo \(p\), then
\begin{equation}\label{core:eq:line-divisibility}
  p\mid(g\mathbf z)_1
  \quad\Longleftrightarrow\quad
  [\mathbf z]_p=\ell_p.
\end{equation}

To average over $P_0$, replace the parameters in its canonical frame by real variables, putting
\begin{equation}\label{core:eq:continuous-frame}
 g(u,v,r)=
 \begin{pmatrix}
 Uu&Uur\\ Vv&Vvr+(Uu)^{-1}
 \end{pmatrix},
 \qquad (u,v,r)\in\mathcal D_\infty:=[1,16]\times[1,2]\times[0,1].
\end{equation}
It has determinant one.  For an integral initial point, its parameters are $(u_0/U,v_0/V,\xi_{P_0}/u_0)$.  For the continuous average, choose $\ell_p$ uniformly and independently in $\mathbb P^1(\mathbb F_p)$ for all $p\in\Pset_*$. Write $\mathbb E_{\boldsymbol\ell}$ for this product expectation. For a test function $\Phi$ of the frame and the lines, our averaging functionals are
\begin{align*}
 \mathcal A_{\rm ar}(\Phi)
 &=\frac1{UV}\sum_{\substack{(u_0,v_0)=1\\
                  u_0/U\in[1,16],\ v_0/V\in[1,2]}}
       \Phi\bigl(g_{P_0},(\ell_p(g_{P_0}))_p\bigr),\\
 \mathcal A_{\rm mod}(\Phi)
 &=\frac1{\zeta(2)}\int_{\mathcal D_\infty}
       \mathbb E_{\boldsymbol\ell}\Phi(g(u,v,r),\boldsymbol\ell)
       \,\dd u\,\dd v\,\dd r.
\end{align*}
The integral with independent choices of lines has total mass $15/\zeta(2)$; the sum over primitive initial points has total mass $O(1)$.

A single finite family of paths suffices for both averages. Indeed, summing the signed slope differences in \eqref{core:eq:position-slope}, we obtain
\[
 \det(P_0,P_j)
 = (P_0)_1(P_j)_1
    \sum_{i<j}\frac{\det(P_i,P_{i+1})}
                        {(P_i)_1(P_{i+1})_1}.
\]
The support of $\psi$ and $D_i<2d_0$ imply $|\det(P_i,P_{i+1})|\le10H$.  Since $\det g=1$ and $P_0=ge_1$,
\[
 |z_{j,2}|=|\det(P_0,P_j)|\le2560TH,
 \qquad
 z_{j,1}=\frac{(P_j)_1}{(P_0)_1}-r z_{j,2},
 \qquad |z_{j,1}|\le16+|z_{j,2}|.
\]
These statements also hold for a continuous frame in \eqref{core:eq:continuous-frame}. We may therefore fix $C_*=10^4$ and use the finite set
\begin{equation}\label{core:eq:uniform-position-set}
 \mathcal Z_0=
 \{\mathbf z\in\mathbb Z^2:\gcd(z_1,z_2)=1,
       \ |z_2|\le C_*TH,\ |z_1|\le16+C_*TH\}.
\end{equation}
Let $\Gamma$ consist of all labelled paths
\begin{equation}\label{core:eq:uniform-path-family}
 \gamma=(\mathbf z_0,\ldots,\mathbf z_T;
                 \mathbf L_0,\ldots,\mathbf L_T),
 \qquad \mathbf z_j\in\mathcal Z_0,
 \quad \mathbf z_0=\mathbf z_T=e_1,
\end{equation}
with the following restrictions. For each $\nu$, each list contains $m_0$ distinct entries from $\Pset_\nu$, and consecutive lists share exactly $J_0$ such entries. Their shared, entering and leaving products $D_j,a_j,c_j$ satisfy
\[
 \mu(D_ja_jc_j)^2=1,\qquad d_0\le D_j<2d_0,
 \qquad D_j\mid\det(\mathbf z_j,\mathbf z_{j+1}),
 \qquad |\det(\mathbf z_j,\mathbf z_{j+1})|\le10H.
\]
Only the endpoint positions are identified. Membership in $\Gamma$ depends on the path in integral coordinates and its labels, so it is independent of the frame.

For $\gamma\in\Gamma$ put
\[
 \begin{split}
 K_\gamma^{\rm box}(g)
  &=\prod_{j=0}^T
       \mathbf1_{U\le(g\mathbf z_j)_1\le16U}
       \mathbf1_{V\le(g\mathbf z_j)_2\le2V},\\
 K_\gamma^{\rm sep}(g)
  &=\prod_{j=0}^T
       \mathbf1_{\rm sep}(\alpha_g(\mathbf z_j),\mathbf L_j),
 \qquad K_\gamma=K_\gamma^{\rm box}K_\gamma^{\rm sep}.
 \end{split}
\]
We set $K_\gamma=0$ outside the coordinate bounds before evaluating the quotients. On the set defined by these bounds, their denominators are nonzero. Primitivity of $\mathbf z_j$ is part of the definition of $\mathcal Z_0$.

For a prime $p$, let
\[
 A_p=\{0\le j\le T:p\in\mathbf L_j\},\qquad
 r_p=\#\{j\in A_p:j=0\text{ or }j-1\notin A_p\}.
\]
The number $r_p$ counts the maximal intervals of consecutive indices $j$ for which $p\in\mathbf L_j$; we call these intervals the runs of $p$. Write $\mathcal P_{\rm act}(\gamma)=\{p:A_p\ne\varnothing\}$. The divisibility conditions for this prime are represented by
\begin{equation}\label{core:eq:local-line-factor}
 L_{\gamma,p}(\ell)=\prod_{j\in A_p}\ind_{\ell=[\mathbf z_j]_p}.
\end{equation}
It equals one for $p\notin\mathcal P_{\rm act}(\gamma)$. For each $\gamma\in\Gamma$, define
\[
 c_\gamma=\kappa_{J_0}^{KT}
   \biggl(\prod_{j=0}^{T-1}
        \mathcal K(\mathbf z_j,\mathbf z_{j+1};D_j,a_j,c_j)\biggr)
   \prod_{\nu=1}^K\prod_{p\in\Pset_\nu}V_\nu^{-r_p}.
\]
The kernel uses $t_j=\det(\mathbf z_j,\mathbf z_{j+1})/D_j$ and hence is independent of $g$. At each step, one entering label is chosen from each set of primes, so $\sum_{p\in\Pset_\nu}r_p=m_0+T$. By \eqref{core:eq:uniform-position-set} each position admits $\exp(O_K(L^{0.21}))$ choices. The same bound holds for its prime tuples since there are $Km_0$ entries, each at most $\exp(2L^{0.2})$. With $T\asymp L^{1/2}$, it follows that $|\Gamma|\le\exp(O_K(L^{0.71}))$. Moreover, $|\mathcal K|\ll1$, $\kappa_{J_0}\le1$ and $V_\nu\asymp1$, so $|c_\gamma|\le\exp(O_K(T+J_0))$. In particular, for a fixed $C_\Gamma$,
\begin{equation}\label{core:eq:path-coefficient-size}
 |\Gamma|,\ \max_{\gamma\in\Gamma}|c_\gamma|
       \le\exp(C_\Gamma L^{0.73}).
\end{equation}
For a family $\Gamma$ and coefficients satisfying this bound, define
\begin{equation}\label{core:eq:joint-test-function}
\begin{aligned}
\mathcal F(g,\boldsymbol\ell)
   & =\sum_{\gamma\in\Gamma}c_\gamma K_\gamma(g)
                  \prod_{p\in\Pset_*}L_{\gamma,p}(\ell_p),\\
\mathcal M_T(g)
    &=\mathbb E_{\boldsymbol\ell}\mathcal F(g,\boldsymbol\ell)
\end{aligned} 
\end{equation}
An explicit expression for \(\mathcal M_T(g)\) appears in \eqref{core:eq:T0-definition}.

\begin{lemma}

For the full path family and the coefficients $c_\gamma$ defined above, every primitive $P_0$ in the coordinate box satisfies
\begin{equation}\label{core:eq:closed-walk-pointwise}
 \bigl\langle u_{P_0},\mathsf A^T u_{P_0}\bigr\rangle
 =\mathcal F\bigl(g_{P_0},(\ell_p(g_{P_0}))_p\bigr).
\end{equation}
Hence
\begin{equation}\label{core:eq:closed-walk-average}
 \frac1{UV}
 \sum_{\substack{P_0\ \mathrm{primitive}\\
                 (P_0)_1/U\in[1,16],\ (P_0)_2/V\in[1,2]}}
 \bigl\langle u_{P_0},\mathsf A^T u_{P_0}\bigr\rangle
 =\mathcal A_{\rm ar}(\mathcal F).
\end{equation}
\end{lemma}

\begin{proof}
Let $a_{\rm c}(x,y)$ for the coefficient of $f(y)$ in $(\mathsf A f)(x)$. Expanding $\mathsf A^T$, we obtain
\begin{equation*}
 \big\langle u_{P_0},\mathsf A^T u_{P_0}\big\rangle
 =w_{\rm st} \sum_{\substack{x_0,\ldots,x_T\in\mathscr S_{d_0}\\
x_0,x_T\text{ above }P_0}}
 \prod_{j=0}^{T-1}a_{\rm c}(x_j,x_{j+1}).
\end{equation*}
The single factor $w_{\rm st}$ comes from the inner product. With $g=g_{P_0}$ fixed, write $x_j=(P_j,\mathbf L_j)$ and $\mathbf z_j=g^{-1}P_j$. Unimodularity preserves primitivity and determinants, and the slope bound puts each $\mathbf z_j$ in \eqref{core:eq:uniform-position-set}. Conversely, a contributing path in $\Gamma$ determines the unique walk $x_j=(g\mathbf z_j,\mathbf L_j)$.

The indicators of the coordinate bounds and separation conditions have product $K_\gamma(g)$. When $p\in\mathbf L_j$, the divisibility indicator is $\ind_{\ell_p(g)=[\mathbf z_j]_p}$ by \eqref{core:eq:line-divisibility}. Thus the conditions at all the primes are expressed by the product of \eqref{core:eq:local-line-factor} over the primes occurring on the path.

Finally $\det(g\mathbf z_j,g\mathbf z_{j+1}) =\det(\mathbf z_j,\mathbf z_{j+1})$, and the scalar left after the kernels have been extracted is
\[
 w_{\rm st}\biggl(\kappa_{J_0}^K\prod_\nu V_\nu^{-1}\biggr)^T=\kappa_{J_0}^{KT}\prod_{\nu=1}^K\prod_{p\in\Pset_\nu}V_\nu^{-r_p},
\]
since, for each $\nu$, $m_0$ runs start at index $0$ and one new run starts at each edge. This proves \eqref{core:eq:closed-walk-pointwise} with the stated $c_\gamma$. Averaging over $P_0$ proves \eqref{core:eq:closed-walk-average}.
\end{proof}

The following proposition supplies the moment estimate in \eqref{core:eq:normalized-high-moment}. We prove it in Section~\ref{core:sec:typeII-minor}.

\begin{proposition}\label{core:prop:path-bound}
For every $E>0$, there exist $A_0>0$ and $K_0\ge1$, depending only on $E$, such that for each fixed $K\ge K_0$ and even $T\asymp L^{1/2}$, $|\mathcal M_T(g)|\ll L^{-ET}$, uniformly in the real frame $g$ and the dyadic range for $D$. In particular,
\begin{equation*}
 \frac1{UV}\sum_{(P_1,P_2)=1}
 \bigl\langle u_P,\mathsf A^{T}u_P\bigr\rangle
 \ll L^{-ET},
\end{equation*}
where the sum is over the primitive positions in the coordinate box.
\end{proposition}

\section{The moment estimate and exceptional terms}\label{core:sec:typeII-minor}

Throughout this section the prime sets satisfy \eqref{core:eq:restricted-band-hyp}. We prove Proposition~\ref{core:prop:path-bound} and estimate the exceptional terms in Proposition~\ref{core:prop:operator-identity}. We first replace the sum over primitive initial points by an integral over real frames and an average over independent projective lines. For a fixed real frame, we construct an auxiliary operator that records repeated prime labels. Lattice and exponential sum estimates bound its norm and imply the required moment estimate. We treat the exceptional terms last and apply Proposition~\ref{core:prop:conditional-comparison} to complete the proof of Proposition~\ref{core:lem:restricted-determinant}.

\subsection{Averaging over the initial point}

For squarefree $S$, write $G_S=\operatorname{SL}_2(\mathbb Z/S\mathbb Z)$.

\begin{lemma}\label{core:lem:finite-line-average}
Let $S$ be squarefree, and let $f_p$ be a complex-valued function on $\mathbb P^1(\mathbb F_p)$ for each $p\mid S$. If $\ell_p(g_0)$ is the kernel of the first row of $g_0\bmod p$, then
\begin{equation*}
 \frac{1}{|G_S|}\sum_{g_0\in G_S}
 \prod_{p\mid S} f_p\bigl(\ell_p(g_0)\bigr)
 =
 \prod_{p\mid S}\biggl(
   \frac{1}{p+1}\sum_{\ell_p\in\mathbb P^1(\mathbb F_p)} f_p(\ell_p)
 \biggr).
\end{equation*}
\end{lemma}

\begin{proof}
For $g_p=g_0\bmod p$ we have $\ell_p(g_0)=g_p^{-1}[0:1]$. By transitivity on $\mathbb P^1(\mathbb F_p)$ \cite[Lemma~3.15]{Pascadi}, the distribution at each $p$ is uniform. Since $S$ is squarefree, independence over $p\mid S$ follows from the Chinese remainder theorem.
\end{proof}

\begin{lemma}
\label{core:lem:joint-root-boundary-geometry} Let $\Omega_\gamma=\{(u,v,r)\in\mathcal D_\infty: K_\gamma(g(u,v,r))=1\}$. For each primitive $P_0=(u_0,v_0)^t$ in the initial box,
\begin{equation}\label{core:eq:discrete-continuous-match}
 K_\gamma(g_{P_0})
 =\ind_{\Omega_\gamma}(u_0/U,v_0/V,\xi_{P_0}/u_0).
\end{equation}
For the grid of mesh $\varepsilon=1/n$, $n\ge2$, there are unions of cubes $\Omega^-\subseteq\Omega_\gamma\subseteq\Omega^+$ and $1\le\mathcal N_\gamma\le\exp(O_K(L^{0.21}))$ such that
\begin{equation}\label{core:eq:R4-geometry-corrected}
 |\Omega^+\setminus\Omega^-|\ll\mathcal N_\gamma\varepsilon,\qquad
 \#\{\text{cubes in }\Omega^+\setminus\Omega^-\}
       \ll\mathcal N_\gamma\varepsilon^{-2}.
\end{equation}
\end{lemma}
\begin{proof}
Substituting $(u,v,r)=(u_0/U,v_0/V,\xi_{P_0}/u_0)$ in \eqref{core:eq:continuous-frame} recovers $g_{P_0}$, hence \eqref{core:eq:discrete-continuous-match}. For $\mathbf z=(z_1,z_2)^t$ put $\tau_{\mathbf z}=z_1+rz_2$. Then
\begin{equation*}
 (g\mathbf z)_1/U=u\tau_{\mathbf z},\qquad
 (g\mathbf z)_2/V=v\tau_{\mathbf z}+z_2/(UVu).
\end{equation*}
Thus the box boundaries at each vertex lie in the zero sets of the polynomials $u\tau_{\mathbf z}-a$ for $a=1,16$, and $uv\tau_{\mathbf z}-v_*u+z_2/(UV)$ for $v_*=1,2$. All have degree at most three.

On the set defined by the coordinate bounds, $r$ lies in the interval $I_{\mathbf z}=\{r\in[0,1]:1/16\le\tau_{\mathbf z}\le16\}$. There we use the real-valued lift $\phi=\tau_{\mathbf w_{\mathbf z}}/\tau_{\mathbf z}$ with $\phi'=\tau_{\mathbf z}^{-2}\le256$. Each value modulo one has at most 257 preimages since $\phi$ is strictly increasing. An empty interval makes the region empty, and a singleton confines it to one plane. On a nondegenerate interval $I_{\mathbf z}$ the separation conditions change only at the preimages of the endpoints of the circle inequalities in Section~\ref{core:sec:diophantine-separation}. This also holds for the weighted condition, a finite sum of indicators. The first condition has at most $m_0^{K}Y^{0.2}$ frequencies $lD$. In the second condition the frequencies are $(n-n')F$, with $n,n'\in\mathcal N_s$ distinct and $F\in\mathcal F_s$. There are at most $\sum_s|\mathcal N_s|^2|\mathcal F_s| \le m_0^{2K}\exp(O_K(L^{0.2}))$ of them. Every frequency $h$ in either condition satisfies $0<|h|\le\exp(O_K(L^{0.21}))$ and contributes at most $2|h|$ endpoints modulo one. Including the ends of the $I_{\mathbf z}$ and the faces of the box for the initial point, the boundary of $\Omega_\gamma$ lies in $\mathcal N_\gamma\le\exp(O_K(L^{0.21}))$ planes and algebraic hypersurfaces of degree at most three.

Partition $\mathcal D_\infty$ into $15n^3$ disjoint grid cubes, assigning each common face to one side. Take for $\Omega^-$ the cubes whose closures lie in the interior of $\Omega_\gamma$, and for $\Omega^+$ those whose closures meet $\overline{\Omega_\gamma}$. Each cube of the difference lies within $2\varepsilon$ of one of the covering zero sets. In a fixed ball containing $\mathcal D_\infty$, Wongkew's tube estimate \cite[Main Theorem]{Wongkew} bounds the volume of the $2\varepsilon$-neighbourhood of each such set by $O(\varepsilon)$, uniformly for nonzero polynomials of degree at most three. This proves the volume bound in \eqref{core:eq:R4-geometry-corrected}. Dividing by $\varepsilon^3$ bounds the number of cubes.
\end{proof}

We count the canonical matrices associated with initial points in each grid cube, uniformly in the sides of the cube and in the modulus $S$. For a related count at fixed level see Boca \cite[\S3]{BocaAngles}.

For $P=(u,v)^t$ abbreviate the canonical frame $g_P$ to $g_{u,v}=\bigl(\begin{smallmatrix}u&c\\v&d\end{smallmatrix}\bigr)$, where $c=\xi_P$ and $d=(1+vc)/u$.

\begin{samepage}
\begin{lemma}
Fix $\delta,C>0$.\label{core:lem:uniform-primitive-frame} Suppose $U,V\ge x^\delta$ and $1\le S\le\exp\bigl(C(\log x)^{0.98}\bigr)$ is squarefree. For sufficiently large $x$, uniformly in $g_0\in G_S$ and intervals $I_1\subset[1,16]$, $I_2\subset[1,2]$, and $I_3\subset[0,1]$, we have, with $c_\delta=\delta/4$,
\begin{align*}
& \#\bigl\{
  (u,v)\in UI_1\times V I_2:
  (u,v)=1,\ 
  c/u\in I_3,\ 
  g_{u,v}\equiv g_0\,({\rm mod}\, S)
 \bigr\}\\
= &
 \frac{UV|I_1||I_2||I_3|}
 {\zeta(2)|\operatorname{SL}_2(\mathbb Z/S\mathbb Z)|}
 +O_{\delta,C}\!\bigl(UVx^{-c_\delta}\bigr).
\end{align*}
\end{lemma}
\end{samepage}

\begin{proof}
We first take $U\le V$, allowing both $I_1,I_2\subseteq[1,16]$. Write $g_0=\bigl(\begin{smallmatrix}u_0&c_0\\v_0&d_0\end{smallmatrix}\bigr)$ and choose integer representatives for its entries. For fixed $u\equiv u_0\pmod S$, the residue conditions on the canonical frame are equivalent to
\begin{equation}\label{core:eq:frame-congruences}
 c\equiv c_0 \,({\rm mod}\, S),\qquad v\equiv v_0 \,({\rm mod}\, S),\qquad
 cv\equiv ud_0-1\,({\rm mod}\,{uS}),
\end{equation}
with $0\le c<u$. The last congruence implies both $(u,v)=1$ and $d=(1+cv)/u\equiv d_0\pmod S$.

Write $u=d_S q$, where $d_S=\prod_{p\mid S}p^{v_p(u)}$, and set $B=d_S S$. Then $uS=Bq$ with $(B,q)=1$. The reduction of \eqref{core:eq:frame-congruences} modulo $B$ admits exactly $d_S$ residue pairs $(c_*,v_*)$. Indeed, for $p\mid S$ with $p^e\Vert u$, $e\ge1$, the determinant condition implies $c_0v_0\equiv-1\pmod p$. Each of the $p^e$ lifts of $c_0\bmod p$ to $p^{e+1}$ is a unit and determines a unique lift of $v_0$. At each $p\mid S$ with $p\nmid u$ the pair is fixed modulo $p$. The count follows by the Chinese remainder theorem.

Modulo $q$ the remaining condition is $cv\equiv-1\pmod q$. For one of the $d_S$ pairs, write $c=c_*+Bi$, $v=v_*+Bj$. If $q>1$ and $B\bar B\equiv1\pmod q$, set
\[
 s=i+\bar B c_*,\qquad t=j+\bar B v_*.
\]
The remaining congruence is $st\equiv-\bar B^2\pmod q$. The two new variables range over independent translates of intervals of lengths $u|I_3|/B$ and $V|I_2|/B$. Taking $f\equiv1$ and subdivision parameter $2$ in \cite[Proposition~4.4]{BZdirections}, we obtain, for intervals $\mathcal I,\mathcal J$ and $(h,q)=1$,
\begin{align*}
 \#\{(s,t)\in(\mathcal I\times\mathcal J)\cap\mathbb Z^2:
          st\equiv h\modq\}
 =\frac{\varphi(q)}{q^2}|\mathcal I|\,|\mathcal J|+O_\varepsilon\!\biggl(q^{1/2+\varepsilon}
      \Bigl(1+\frac{|\mathcal I|}{q}\Bigr)
      \Bigl(1+\frac{|\mathcal J|}{q}\Bigr)\biggr).
\end{align*}
The condition $(h,q)=1$ forces both coordinates to be units modulo $q$, so fixing an endpoint changes the count by $O(1+|\mathcal I|/q+|\mathcal J|/q)$. For $q=1$ the formula follows by counting the integers in each interval. Since $u/B=q/S$, summing over the $d_S$ residue pairs shows that the number of canonical matrices for this $u$ is
\begin{equation*}
 N(u)=\frac{V|I_2||I_3|}{S^2}\frac{\varphi(q)}q
       +O_\varepsilon\!\biggl(d_S q^{1/2+\varepsilon}
                     \Bigl(1+\frac{V}{uS}\Bigr)\biggr).
\end{equation*}

For the total error it suffices to sum over all $u\le16U$. Writing $d_S=d_S(u)$, we have
\[
 \sum_{u\le16U} d_S(u)^{1/2}
 \le 16U\prod_{p\mid S}(1-p^{-1/2})^{-1}\le 16U S^2.
\]
Indeed, a given $d_S$ occurs at most $16U/d_S$ times, and $(1-p^{-1/2})^{-1}\le p^2$. On $u\asymp U$ one has $d_S q^{1/2+\varepsilon}\ll U^{1/2+\varepsilon} d_S^{1/2}$. Since $U\le V$ and $1+V/(uS)\ll V/U$ on this range,
\begin{equation}\label{core:eq:frame-total-error}
 \sum_{u/U\in I_1} d_S q^{1/2+\varepsilon}
       \Bigl(1+\frac{V}{uS}\Bigr)
 \ll_\varepsilon VU^{1/2+\varepsilon}S^2 .
\end{equation}
Take $\varepsilon=1/8$. For large $x$ in terms of $\delta,C$ we have $S^2\le U^{1/8}$, so the total error is $\ll VU^{3/4}\le UVx^{-\delta/4}$.

For the main term use $\varphi(q)/q=\sum_{r\mid u,\,(r,S)=1}\mu(r)/r$. The conditions $r\mid u$ and $u\equiv u_0\pmod S$ specify one class modulo $rS$. Thus
\begin{equation}\label{core:eq:frame-u-average}
 \sum_{\substack{u/U\in I_1\\ u\equiv u_0\,({\rm mod}\,S)}}
    \frac{\varphi(q)}q
 =\frac{U|I_1|}{S\zeta(2)}\prod_{p\mid S}(1-p^{-2})^{-1}
     +O(\log(2U)),
\end{equation}
since the count in each class has error $O(1)$ and the tail $(U/S)\sum_{r>16U}r^{-2}$ is bounded. The group has order $|G_S|=S^3\prod_{p\mid S}(1-p^{-2})$. Multiplying \eqref{core:eq:frame-u-average} by $V|I_2||I_3|/S^2$, we obtain the asserted main term. Its error $O(V\log(2U)/S^2)$ is absorbed by \eqref{core:eq:frame-total-error}.

If $V<U$, exchange $u,v$. The canonical frame is then
\[
 g_{v,u}=\begin{pmatrix}v&v-d\\u&u-c\end{pmatrix},
 \qquad \frac{v-d}{v}=1-\frac{c}{u}-\frac1{uv}.
\]
Indeed $0\le v-d<v$ and $u(v-d)\equiv-1\pmod v$. This involution sends $g_0$ to a single prescribed class modulo $S$. Put $J_3=1-I_3$ and $\Delta=(UV)^{-1}$. The transformed third parameter differs from $1-c/u$ by at most $\Delta$. The desired count lies between the swapped counts with $J_3$ contracted and enlarged at both endpoints by $\Delta$, respectively, intersected with $[0,1]$ and allowing an empty inner interval. Their lengths differ from $|I_3|$ by $O(\Delta)$, so the change in the main term is $O(1)$. Applying the first case with scales $(V,U)$ bounds the error by $O_{\delta,C}(UVx^{-\delta/4})$.
\end{proof}

\begin{lemma}\label{core:lem:bounded-modulus}
Let $K$ be an indicator satisfying the identity and grid bounds of Lemma~\ref{core:lem:joint-root-boundary-geometry}, with support $\Omega$ and constant $\mathcal N$. Let $S$ be squarefree with $\log S\le C L^{0.98}$, and let $0\le\Phi\le1$ be a function of the lines at the primes dividing $S$. For $U,V\ge x^\delta$ and $\varepsilon=1/n$, $n\ge2$, we have
\begin{equation*}
 \Bigl|\mathcal A_{\rm ar}(K\Phi)
       -\frac{|\Omega|}{\zeta(2)}\mathbb E_{\boldsymbol\ell}\Phi\Bigr|
 \ll_{\delta,C}\mathcal N\varepsilon+
        x^{-c_\delta}\varepsilon^{-3} S^3 .
\end{equation*}
\end{lemma}

\begin{proof}
On each grid cube $B$, sum Lemma~\ref{core:lem:uniform-primitive-frame} over $G_S$ with weight $\Phi$. Lemma~\ref{core:lem:finite-line-average}, applied to the indicators of individual tuples of lines, identifies the main term as $|B|\mathbb E_{\boldsymbol\ell}\Phi/\zeta(2)$. The error is $O_{\delta,C}(x^{-c_\delta}|G_S|)$. Since $0\le\Phi\le1$, summing over the inner and outer approximations to $\Omega$ bounds $\mathcal A_{\rm ar}(K\Phi)$ on either side. Now $|G_S|\le S^3$, each union has at most $15\varepsilon^{-3}$ cubes and their volumes differ by $O(\mathcal N\varepsilon)$. This proves the lemma.
\end{proof}

When $S,\mathcal N=x^{o(1)}$, the error in Lemma~\ref{core:lem:bounded-modulus} with $n=\lceil x^{c_\delta/4}\rceil$ is $x^{-c_\delta/4+o(1)}$. We use this choice both for a single edge and for a path of length $T$. For a path, the modulus in Lemma~\ref{core:lem:bounded-modulus} is the product of its distinct prime labels.

\begin{proposition}\label{core:prop:joint-root-replacement}
Let $T\asymp L^{1/2}$ be even, with fixed implied constants, and $U,V\ge x^\delta$. For any subfamily of \eqref{core:eq:uniform-path-family} and coefficients independent of the frame and lines satisfying \eqref{core:eq:path-coefficient-size}, the function \eqref{core:eq:joint-test-function} satisfies
\begin{equation*}
 \mathcal A_{\rm ar}(\mathcal F)
   =\mathcal A_{\rm mod}(\mathcal F)+O(x^{-\delta/20}).
\end{equation*}
The estimate is uniform in $U,V$ and these choices, with its constant and the lower threshold for $x$ depending only on $C_\Gamma$ and the fixed parameters.
\end{proposition}
\begin{proof}
For a given path put $q_\gamma=\prod_{p\in\mathcal P_{\rm act}(\gamma)}p$. There are $m_0$ initial labels from each set of primes and at most $T$ further choices from that set, so $\log q_\gamma\ll_K(m_0+T)L^{0.2}\ll_K L^{0.70}$. The product of divisibility indicators in \eqref{core:eq:joint-test-function} takes values in $[0,1]$ and depends only on the lines at $p\mid q_\gamma$. By Lemma~\ref{core:lem:joint-root-boundary-geometry}, $\mathcal N_\gamma=x^{o(1)}$. Thus the error in Lemma~\ref{core:lem:bounded-modulus} with $S=q_\gamma$ and $n=\lceil x^{c_\delta/4}\rceil$ is $x^{-\delta/16+o(1)}$ per path, since $c_\delta=\delta/4$. By \eqref{core:eq:path-coefficient-size}, $\sum_\gamma|c_\gamma|\le\exp(2C_\Gamma L^{0.73})=x^{o(1)}$. Summing proves the proposition.
\end{proof}

For each prime $p$ occurring on the path, let $I_p(\gamma)$ be the indicator that the lines $[\mathbf z_j]_p$ with $j\in A_p$ all coincide. Their common line has probability $1/(p+1)$. By independence,
\begin{equation}\label{core:eq:T0-definition}
 \mathcal M_T(g)=
 \sum_{\gamma:K_\gamma(g)=1}\kappa_{J_0}^{KT}
 \prod_{j<T}\mathcal K_j(\gamma)\,
 \prod_{\nu=1}^K\prod_{\substack{p\in\Pset_\nu\\ A_p\ne\varnothing}}
 \frac{I_p(\gamma)}{V_\nu^{r_p}(p+1)},
\end{equation}
where $\mathcal K_j(\gamma)$ is the kernel on the $j$th edge. 
Consequently,
\begin{equation*}
\mathcal A_{\rm mod}(\mathcal F)
= \frac1{\zeta(2)} \int_{\mathcal D_\infty}
       \mathbb E_{\boldsymbol\ell} \mathcal M_T(g(u,v,r),\boldsymbol\ell)
       \,\dd u\,\dd v\,\dd r.
\end{equation*}

\subsection{Construction of an auxiliary operator}

Fix a frame $g$ in \eqref{core:eq:continuous-frame}. In \eqref{core:eq:T0-definition}, a prime $p\in\Pset_\nu$ with $r_p$ runs contributes $V_\nu^{-r_p}/(p+1)$ if the prescribed lines agree in all its runs. When $p$ next occurs, we must impose the same line condition as in its preceding run. We therefore record the prime and its line between consecutive runs. In Pilatte's additive paths, a prime dividing two positions likewise divides the sum of the intervening increments \cite[Definitions~8.2 and~8.4(2)]{Pilatte}.

Let $\mathscr X_g$ consist of triples $\xi=(\mathbf z,\mathbf L,\mathcal C)$, where $\mathbf z\in\mathcal Z_0$, $g\mathbf z\in[U,16U]\times[V,2V]$, and $\alpha_g(\mathbf z)$ is separated for $\mathbf L$. For each $\nu$, the tuple $\mathbf L_\nu$ consists of $m_0$ distinct primes from $\Pset_\nu$.
Here $\mathcal C=(\mathcal C_1,\ldots,\mathcal C_K)$, where each $\mathcal C_\nu$ is a set of pairs $(p,\ell)$ with $p\in\Pset_\nu$ and $\ell\in\mathbb P^1(\mathbb F_p)$. Every prime recorded in $\mathcal C$ occurs in just one pair and is absent from the current tuple. We require $|\mathcal C|:=\sum_\nu |\mathcal C_\nu|\le B_*:=KT$. Each step adds at most one pair to each of these sets, so every sequence of $T$ transitions starting with $\mathcal C=\varnothing$ satisfies this bound.

For a prime $p\in\Pset_\nu$, put
\[
 \pi_\nu(p)=\frac{1}{V_\nu(p+1)},\qquad
 s_\nu=\sum_{p\in\Pset_\nu}\pi_\nu(p)\le 1.
\]
Equip $\mathscr X_g$ with the measure
\begin{equation*}
 \varpi(\mathbf z,\mathbf L,\mathcal C)=
 \prod_{\nu=1}^K
 \biggl(\prod_{p\in\mathbf L_\nu}\pi_\nu(p)\biggr)
 \biggl(\prod_{(p,\ell)\in\mathcal C_\nu}\frac1{p+1}\biggr)
\end{equation*}
with counting measure on positions. The function $\mathbf b=\ind_{\mathbf z=e_1}\ind_{\mathcal C=\varnothing}$ satisfies
\begin{equation*}
 \|\mathbf b\|^2\le\prod_{\nu=1}^K s_\nu^{m_0}\le1.
\end{equation*}
Indeed, the product on the right is obtained by dropping the separation and distinctness restrictions.

In a transition from $\xi$ to $\xi'$, the two tuples share $J_0$ primes from each set. Write $c_\nu$ for the leaving prime, $a_\nu$ for the entering prime, $D$ for the product of the shared primes, and $c=\prod_\nu c_\nu$, $a=\prod_\nu a_\nu$. We require
\begin{equation*}
 \mu(Dac)^2=1,\qquad d_0\le D<2d_0,\qquad
 D\mid\det(\mathbf z,\mathbf z'),\qquad
 |\det(\mathbf z,\mathbf z')|\le 10H.
\end{equation*}
We evaluate $\mathcal K=\mathcal K(\mathbf z,\mathbf z';D,a,c)$ from \eqref{core:eq:coarse-kernel} at $t=\det(\mathbf z,\mathbf z')/D$.

A pair $\boldsymbol\epsilon,\boldsymbol\delta\in\{0,1\}^K$ specifies how $\mathcal C$ changes. If $\delta_\nu=1$, choose the pair $(a_\nu,[\mathbf z']_{a_\nu})$ from $\mathcal C_\nu$ and remove it. If $\delta_\nu=0$, choose $a_\nu$ outside the primes occurring in $\mathcal C_\nu$. If $\epsilon_\nu=1$, adjoin $(c_\nu,[\mathbf z]_{c_\nu})$ to the resulting set. If $\epsilon_\nu=0$, do not record $c_\nu$. These choices determine $\mathcal C'_\nu$. Define the transition operator by
\begin{equation}\label{core:eq:forward-memory-operator}
 (\mathsf T_{\boldsymbol\epsilon,\boldsymbol\delta}f)(\xi)
 =\sum_{\xi'}\kappa_{J_0}^K\,\mathcal K
 \prod_{\nu=1}^K
 \pi_\nu(a_\nu)^{1-\delta_\nu} V_\nu^{-\delta_\nu} f(\xi'),
 \qquad
 \mathsf T=\sum_{\boldsymbol\epsilon,\boldsymbol\delta}
           \mathsf T_{\boldsymbol\epsilon,\boldsymbol\delta}.
\end{equation}
Here the sum is over the admissible transitions with $\xi'\in\mathscr X_g$. The kernel relative to $\varpi$ is the coefficient of $f(\xi')$ in this sum divided by $\varpi(\xi')$.

Reversing an edge exchanges addition to and removal from $\mathcal C$. To check the adjoint, fix an index and omit the factors from the shared tuple and from pairs unchanged in $\mathcal C$. The source mass times the coefficient associated with the entering and leaving primes is
\[
 \pi_\nu(c_\nu)(a_\nu+1)^{-\delta_\nu}
 \pi_\nu(a_\nu)^{1-\delta_\nu}V_\nu^{-\delta_\nu}
 =\pi_\nu(c_\nu)\pi_\nu(a_\nu),
\]
because $1/(V_\nu(p+1))=\pi_\nu(p)$. In reverse the same expression uses $c_\nu$ in place of $a_\nu$ and $\epsilon_\nu$ in place of $\delta_\nu$. The kernel is real and symmetric. Consequently
\begin{equation}\label{core:eq:fine-pattern-adjoint}
 \mathsf T_{\boldsymbol\epsilon,\boldsymbol\delta}^*
     =\mathsf T_{\boldsymbol\delta,\boldsymbol\epsilon},
 \qquad \mathsf T^*=\mathsf T .
\end{equation}

Let $\flat$ restrict the expansion of $\langle\mathbf b,\mathsf T^T\mathbf b\rangle$ to sequences for which the initial tuple entries and all values $a_\nu$ chosen with $\delta_\nu=0$ are distinct. A sequence specifies all the choices in \eqref{core:eq:forward-memory-operator}, including the pairs added to and removed from $\mathcal C$.

\begin{proposition}\label{core:prop:exact-fine-transfer}
For every frame $g$ in \eqref{core:eq:continuous-frame},
\begin{equation*}
 \mathcal M_T(g)=\langle\mathbf b,\mathsf T^T\mathbf b\rangle^\flat .
\end{equation*}
\end{proposition}

\begin{proof}
Fix a path of positions and tuples in \eqref{core:eq:T0-definition}. Within a run of a prime $p$, the line $[\mathbf z_j]_p$ is constant. Indeed, at each internal edge $p\mid D$, so $\det(\mathbf z_j,\mathbf z_{j+1})\equiv0\pmod p$, and both vectors are primitive.

On its first run $p\in\Pset_\nu$ contributes $\pi_\nu(p)$, through the initial state mass if the run starts at time $0$ and through the coefficient for $\delta_\nu=0$ otherwise. Between consecutive runs, $\flat$ forces the pair to be added to $\mathcal C_\nu$ at departure and selected from $\mathcal C_\nu$ at the next entry. This selection requires the same line and contributes $V_\nu^{-1}$. At a departure after the last run we must have $\epsilon_\nu=0$, since the terminal set $\mathcal C$ is empty.

Thus the path contributes zero if its prescribed lines disagree. Otherwise the choices of $\boldsymbol\epsilon,\boldsymbol\delta$ are unique, and each prime with $r_p\ge1$ runs contributes
\[
 \pi_\nu(p)V_\nu^{-(r_p-1)}
   =\frac{V_\nu^{-r_p}}{p+1},
\]
while each absent prime contributes $1$. The product of these factors with $\kappa_{J_0}^{KT}\prod_{j<T}\mathcal K_j$ is precisely the summand in \eqref{core:eq:T0-definition}.
\end{proof}

\subsection{Lattice and exponential sum estimates}

We next choose an integral basis using the first separation condition. In these coordinates, the required bound follows from the additive large sieve. Write $\mathfrak M_{A_1}$ for the arcs in \eqref{core:eq:major-arcs-def} with exponent $A_1>0$ in place of $A_0$, and $\mathfrak m_{A_1}=\mathbb T\setminus\mathfrak M_{A_1}$.

\begin{lemma}\label{core:lem:balanced-box}
Let $C_0>0$ and let $\alpha\in\mathbb R$ satisfy $\|l\alpha\|>Y^{-0.7}$ for $1\le |l|\le Y^{0.2}$. There are $A,B>0$ and a basis $\{b_1,b_2\}$ with $(b_1\ b_2)\in\operatorname{SL}_2(\mathbb Z)$ such that every $(h,l)^t\in\mathbb Z^2$ with $|h+\alpha l|\le C_0$ and $|l|\le C_0Y$ has the form $mb_1+nb_2$ with
\[
 |m|\ll A,\qquad |n|\ll B,
 \qquad AB\ll_{C_0}Y,
 \qquad Y^{0.1}\ll A,B\ll Y^{0.9}.
\]
\end{lemma}

\begin{proof}
Use
\[
 \Phi_\alpha(h,l)=(h+\alpha l,l/Y),
 \qquad \Lambda_\alpha=\Phi_\alpha(\mathbb Z^2),
 \qquad \operatorname{covol}(\Lambda_\alpha)=Y^{-1}.
\]
If the shortest nonzero vector had length $<Y^{-0.8}$, its second coordinate would force $1\le |l|<Y^{0.2}$ and its first coordinate would imply $\|l\alpha\|<Y^{-0.8}$, a contradiction.  Dirichlet approximation supplies $1\le l\le Y^{1/2}$ with $\|l\alpha\|\le Y^{-1/2}$, so, writing the first minimum as $\lambda_1$,
\[
 Y^{-0.8}\le\lambda_1\ll Y^{-1/2}.
\]
Use the reduced basis in Heath-Brown \cite[pp.~7--8]{HB12}, taking $E=Y^2$, $M=Y$ and $s_0=-\alpha$ in his lattice. Its vectors $v_1,v_2$ have lengths $\lambda_1\le\lambda_2$, satisfy $\lambda_1\lambda_2\asymp Y^{-1}$, and have the following coordinate bound. If $mv_1+nv_2$ belongs to the fixed box $[-C_0,C_0]^2$, then $|m|\ll_{C_0}\lambda_1^{-1}$ and $|n|\ll_{C_0}\lambda_2^{-1}$. Pull the basis back by $\Phi_\alpha$ and choose its orientation to have determinant one. The choices $A=\lambda_1^{-1}$ and $B=\lambda_2^{-1}\asymp Y\lambda_1$ satisfy the asserted ranges and $AB\asymp Y$.
\end{proof}

\begin{lemma}\label{core:lem:minor-determinant-matrix}
Fix $A_1\ge G>0$, $0<\sigma<1/2$, $\upsilon>0$ and $\psi\in C_c^\infty(\mathbb R)$. Suppose $L^\upsilon\ll\log Y$, $Y^\sigma\le A,B\le Y^{1-\sigma}$ and $AB\ll Y$. Then for every $\theta\in\mathfrak m_{A_1}$,
\[
 \Bigl\|\bigl(Y^{-1}\psi(\det(u,v)/Y)e(\theta\det(u,v))\bigr)
 _{u\in\mathcal B_1,v\in\mathcal B_2}\Bigr\|_{\ell^2(\mathcal B_2)\to\ell^2(\mathcal B_1)}
 \ll_{G,A_1,\sigma,\upsilon,\psi}L^{-G}
\]
for rectangles $\mathcal B_i\subset\mathbb Z^2$ of side lengths $O(A)\times O(B)$.
\end{lemma}

\begin{proof}
Fourier inversion gives
\[
 Y^{-1}\psi(\det(u,v)/Y)e(\theta\det(u,v))
 =Y^{-1}\int_{\mathbb R}\widehat\psi(s)
 e\!\biggl(\Bigl(\theta+\frac{s}{Y}\Bigr)\det(u,v)\biggr)\dd s.
\]
Since both rectangles have $O(Y)$ points, the part $|s|>L^{A_1}/2$ contributes $O_R(L^{-R})$ to the operator norm in the lemma, for every $R>0$. For the remaining part put $\theta'=\theta+s/Y$ and $Q=\lfloor Y/L^{A_1}\rfloor$.  Dirichlet's theorem provides a reduced $b_0/q$ with $q\le Q$ and $|\theta'-b_0/q|\le(qQ)^{-1}$. If $q\le L^{A_1}$, then $\|\theta-b_0/q\|\le(qQ)^{-1}+L^{A_1}/(2Y)<2L^{A_1}/Y$ for large $Y$, contrary to $\theta\in\mathfrak m_{A_1}$. Thus $L^{A_1}<q\le Y/L^{A_1}$.

Put $T_{I,J}(\theta')=(e(\theta' mn))_{m\in I,n\in J}$. Writing $\mathcal B_1=I_1\times I_2$ and $\mathcal B_2=J_1\times J_2$, we obtain
\[
 \bigl(e(\theta'\det(u,v))\bigr)_{u,v}
 \cong T_{I_1,J_2}(\theta')\otimes
 \overline{T_{I_2,J_1}(\theta')}.
\]
For integer intervals $I,J$ of cardinalities $M,N$, partition $I$ into $O(M/q+1)$ consecutive blocks of at most $\lfloor q/2\rfloor$ integers.  Since $|\theta'-b_0/q|\le q^{-2}$, distinct $m,m'$ in one block satisfy
\[
 \|(m-m')\theta'\|
 \ge\Bigl\|\frac{(m-m')b_0}{q}\Bigr\|
              -\frac{|m-m'|}{q^2}\ge\frac1{2q}.
\]
The additive large sieve of Montgomery--Vaughan \cite[Theorem~1]{MV73}, applied to these $1/(2q)$-separated frequencies, implies
\[
 \sum_{m\text{ in one block}}
       \biggl|\sum_{n\in J}f_n e(\theta' mn)\biggr|^2
 \le (N+2q)\sum_{n\in J}|f_n|^2.
\]
Summing over the blocks, we obtain
\[
 \|T_{I,J}(\theta')\|^2
 \ll (M/q+1)(N+q).
\]
Consequently the tensor-product norm satisfies
\[
 \frac{\|T_{I_1,J_2}\|\|T_{I_2,J_1}\|}{Y}
 \ll \frac{(A+q)(B+q)}{qY}
 =\frac{AB}{qY}+\frac{A+B}{Y}+\frac qY
 \ll L^{-A_1}+Y^{-\sigma}\ll L^{-G}.
\]
Here $AB\ll Y$, $A,B\le Y^{1-\sigma}$ and $L^{A_1}<q\le Y/L^{A_1}$. The assertion follows on integrating against $|\widehat\psi(s)|$.
\end{proof}

To apply Lemma~\ref{core:lem:minor-determinant-matrix} to a transition with $\boldsymbol\epsilon=\boldsymbol\delta=0$, we use the same integral basis at both ends of each block.

\begin{lemma}\label{core:lem:common-fresh-basis}
Fix the ordered shared labels with product $D\in[d_0,2d_0)$, the sets $\mathcal C$, and a projective class modulo $D$. Partition the slopes $(g\mathbf z)_2/(g\mathbf z)_1$ into intervals of length $H/U^2$. For two intervals whose indices differ by at most $11$, if the source interval contains a state with these data, both position sets admit common integral coordinates $u,v$ in rectangles as in Lemma~\ref{core:lem:balanced-box}, and
\[
 \frac{\det(\mathbf z,\mathbf z')}{D}=\det(u,v).
\]
The norm of the operator with $\boldsymbol\epsilon=\boldsymbol\delta=0$ is at most an absolute constant times the largest block norm.
\end{lemma}

\begin{proof}
On $\ell^2(\mathscr X_g,\varpi)$, let $\mathsf S$ average the tuple orders while fixing $\mathbf z$ and $\mathcal C$. Since $\varpi$ is invariant under these permutations, $\|\mathsf S\|\le1$. Define $\mathsf T_{\rm fixed,0}$ by changing only the last slot of each tuple, preserving the first $J_0$ slots in order, and using the coefficient in \eqref{core:eq:forward-memory-operator} with $\boldsymbol\epsilon=\boldsymbol\delta=\mathbf0$ and the factor $\kappa_{J_0}^K$ omitted. For each index $\nu$, exactly $J_0!$ of the $(m_0!)^2$ pairs of endpoint orders have the shared entries in the first $J_0$ slots in the same order. Since $J_0!/(m_0!)^2=\kappa_{J_0}$, we obtain
\[
 \mathsf T_{\mathbf0,\mathbf0}
 =\mathsf S\mathsf T_{\rm fixed,0}\mathsf S,
 \qquad
 \|\mathsf T_{\mathbf0,\mathbf0}\|
 \le \|\mathsf T_{\rm fixed,0}\|.
\]
The operator $\mathsf T_{\rm fixed,0}$ preserves the ordered shared tuple and $\mathcal C$. Since $D$ is squarefree and both vectors are primitive, the condition $D\mid\det(\mathbf z,\mathbf z')$ implies that the source and target have the same projective class modulo every prime $p\mid D$. These data define orthogonal blocks.

Since $|\det(\mathbf z,\mathbf z')|\le10H$, the two slopes differ by at most $10H/U^2$. Hence the indices of their intervals differ by at most $11$.  For a fixed offset the source intervals and target intervals are separately disjoint, so the operator norm is the supremum of the block norms.  Summing the bounded number of offsets multiplies the bound by an absolute constant.

Choose the position $\mathbf z_\circ$ of any source state in the interval and an integral complement $w_\circ$ with $\det(\mathbf z_\circ,w_\circ)=1$.  All positions in the block belong to the same index-$D$ lattice
\[
 \Lambda=\{\mathbf z:D\mid\det(\mathbf z_\circ,\mathbf z)\}
       =\mathbb Z\mathbf z_\circ+\mathbb Z D w_\circ.
\]
Write $\mathbf z=h\mathbf z_\circ+Dl w_\circ$ and put $r_\circ=(g w_\circ)_1/(g\mathbf z_\circ)_1$. The common or neighboring slope intervals ensure $|\det(\mathbf z_\circ,\mathbf z)|\ll H$, while the first-coordinate box imposes
\[
 |l|\ll H/D\ll Y,\qquad
 |h+Dl r_\circ|
       =\frac{(g\mathbf z)_1}{(g\mathbf z_\circ)_1}\le16.
\]
The shared $D$ occurs among the products tested at the chosen source, so \eqref{core:eq:separated-state-indicator} applies to $D r_\circ$.  Apply Lemma~\ref{core:lem:balanced-box} to $(h,l)$ once, using the same basis for both endpoints.  The map from these new coordinates to $\mathbf z$ has determinant $D$, proving the determinant identity. In these coordinates the restriction $|\det(\mathbf z,\mathbf z')|\le10H$ holds throughout the support of $\psi(t/Y)$ since $D<2d_0$ and $|t|\le5Y$. The box, primitivity, membership in $\mathcal Z_0$, separation, and the requirement that no tuple prime occurs in the fixed sets $\mathcal C$ are all conditions on the individual endpoints. The only remaining restriction coupling the two label tuples is $a_\nu\ne c_\nu$ for every index. The chosen basis therefore also applies after this restriction is removed.
\end{proof}

\subsection{The moment estimate}

For $\boldsymbol\epsilon=\boldsymbol\delta=0$, we use an exponential sum estimate. For the remaining choices, Schur's test \cite{Schur} bounds the operator norm on $\ell^2(\varpi)$ by the geometric mean of the maximal absolute row sum and the maximal weighted column sum. In the column sum, the absolute coefficient of the transition from $\xi$ to $\xi'$ is weighted by $\varpi(\xi)/\varpi(\xi')$.

\begin{lemma}\label{core:lem:transfer-schur}
For each $(\boldsymbol\epsilon,\boldsymbol\delta)$, the row and weighted column sums of the absolute transition coefficients in \eqref{core:eq:forward-memory-operator} are $\ll_K(1+B_*)^K(1+L^{3A_0})$. In either orientation, fixing $a_\nu=p$ with $\delta_\nu=0$ contributes an additional factor $\pi_\nu(p)$ in the row bound.
\end{lemma}
\begin{proof}
Fix a source $(\mathbf z,\mathbf L,\mathcal C)$. For each $D,a,c$ write $\mathbf z'=h\mathbf z+Dt w_{\mathbf z}$ with $\det(\mathbf z,w_{\mathbf z})=1$. For a given $t$, the box for the first coordinate permits $O(1)$ integers $h$. Since the measure of $\mathfrak M$ is $O(L^{3A_0}/Y)$,
\begin{equation}\label{core:eq:positive-edge-majorant}
 |\mathcal K|\ll |\eta(c/Y)\eta(a/Y)\psi(t/Y)|
                 (\ind_{t=c-a}+L^{3A_0}/Y).
\end{equation}
Summing this upper bound over $h,t$, we obtain $O(1+L^{3A_0})$. For each tuple, the $m_0$ choices of the prime that leaves it and the $m_0!$ orders of the target tuple cancel $\kappa_{J_0}=1/(m_0m_0!)$. For $\delta_\nu=0$, summation over $a_\nu$ contributes $s_\nu\le1$. For $\delta_\nu=1$, there are at most $|\mathcal C_\nu|$ choices, each of weight $1/V_\nu$. Thus the sum of absolute values in the row is bounded by
\[
 \ll_K(1+L^{3A_0})\prod_\nu
                 (s_\nu+|\mathcal C_\nu|/V_\nu).
\]
Fixing $a_\nu=p$ when $\delta_\nu=0$ replaces $s_\nu$ by $\pi_\nu(p)$. Reversing the transition as in \eqref{core:eq:fine-pattern-adjoint} proves the column bound.
\end{proof}

We allow a separate multiplier on each $a_\nu$ chosen with $\delta_\nu=0$. This uniformity allows us to impose the restriction $\flat$ in the proof of Proposition~\ref{core:prop:path-bound}.

\begin{lemma}\label{core:lem:fresh-reused-edges}
Suppose $2\le T\ll L^{1/2}$ and $L^{0.1}\ll\log Y\ll_K L^{0.2}$. Given $G>0$, take $A_0=G+3$ and $K_0$ with $0.005K_0>G+1$. For each fixed $K\ge K_0$, uniformly in $g,\boldsymbol\epsilon,\boldsymbol\delta$,
\[
 \|\mathsf T_{\boldsymbol\epsilon,\boldsymbol\delta}\|_{\ell^2(\varpi)\to\ell^2(\varpi)}
 \ll_{K,G} L^{-G}.
\]
The same bound holds in either orientation with factors $\chi_\nu(a_\nu)$, $|\chi_\nu|\le1$, for $\delta_\nu=0$.
\end{lemma}

\begin{proof}
First suppose $\boldsymbol\epsilon=\boldsymbol\delta=0$. By the norm reduction in Lemma~\ref{core:lem:common-fresh-basis}, it suffices to estimate the operator with the first $J_0$ slots of each tuple fixed. Keep the distinctness condition within each endpoint tuple, but extend its kernel to pairs with $a_\nu=c_\nu$ for some $\nu$. For fixed source, this equality determines $a_\nu$ and contributes at most $\Delta=\max_{\nu,p}\pi_\nu(p)\ll_K p_0^{-1}$. The leaving slots and target order are fixed, so \eqref{core:eq:positive-edge-majorant} bounds the added row sum by $O_K(\Delta(1+L^{3A_0}))$, after summing over $\nu$. The same column bound follows by reversal. By Schur's test the extension changes the norm by $O_A(L^{-A})$ for every fixed $A>0$, uniformly in the multipliers.

Within a block of Lemma~\ref{core:lem:common-fresh-basis}, write $m(c)=\prod_\nu\pi_\nu(c_\nu)$ and similarly for $a$. The source and target have masses $M_{\rm sh}m(c)$ and $M_{\rm sh}m(a)$, where $M_{\rm sh}$ is constant on the block. Conjugating to counting measure changes the coefficient $m(a)$ to $\sqrt{m(c)m(a)}$. In common integral coordinates put
\[
 B_\theta(u,v)=\psi(\det(u,v)/Y)e(\theta\det(u,v)),\qquad
 (\mathcal U_\theta F)(u)
     =\sum_c\sqrt{m(c)}\eta(c/Y)e(\theta c)F(u,c).
\]
Use superscripts $\mathrm s,\mathrm t$ for the latter maps at the two endpoints. Fourier orthogonality for the integer $t=\det(\mathbf z,\mathbf z')/D=\det(u,v)$ identifies the enlarged block with
\begin{equation}\label{core:eq:fresh-frequency-factorization}
 \frac{d_0}{D}P_{\rm s}\int_{\mathbb T\setminus\mathfrak M_{A_0}}
    (\mathcal U_\theta^{\rm s})^* B_\theta
       \mathcal U_\theta^{\rm t}\,\dd\theta\,P_{\rm t}.
\end{equation}
Here $P_{\rm s},P_{\rm t}$ are the diagonal contractions imposing the endpoint conditions of Lemma~\ref{core:lem:common-fresh-basis}. Since the shared entries are fixed, the factor $\prod_\nu\chi_\nu(a_\nu)$ depends only on the target tuple and may be included in $P_{\rm t}$. Distinct tuples $c$ determine distinct integer products and hence distinct integer frequencies $c$. Parseval therefore yields
\begin{equation}\label{core:eq:fresh-label-parseval}
 \int_{\mathbb T}\|\mathcal U_\theta F\|_2^2\,\dd\theta
 =\sum_{u,c}m(c)|\eta(c/Y)|^2|F(u,c)|^2
 \ll_K Y^{-1}\|F\|_2^2 ,
\end{equation}
since $m(c)\le(\mathcal V(\Pset)c)^{-1}$ and $c\asymp Y$ on the support of $\eta$. By Cauchy--Schwarz in $\theta$ and \eqref{core:eq:fresh-label-parseval}, the block norm in \eqref{core:eq:fresh-frequency-factorization} is $\ll_K Y^{-1}\sup_{\theta\in\mathfrak m_{A_0}}\|B_\theta\|$. Lemma~\ref{core:lem:minor-determinant-matrix} with $A_1=A_0$, $\sigma=1/20$ and $\upsilon=0.1$ bounds the supremum by $O(YL^{-G})$. The side lengths in Lemma~\ref{core:lem:balanced-box} satisfy its hypotheses for large $Y$. Combining the block estimates and applying the norm reduction in Lemma~\ref{core:lem:common-fresh-basis} proves the estimate for $\boldsymbol\epsilon=\boldsymbol\delta=\mathbf0$.

\smallskip\noindent

We next estimate the term $t=c-a$ for the remaining choices of $\boldsymbol\epsilon,\boldsymbol\delta$. Fix a source $(\mathbf z,\mathbf L,\mathcal C)$ and write $P_*=\Pi(\mathbf L)=Dc$, which is independent of the leaving choices. When $\boldsymbol\delta=0$, no entering prime is chosen from $\mathcal C$, and for each choice of labels $t=c-a$ permits $O(1)$ targets. Summing the corresponding coefficients as in Lemma~\ref{core:lem:transfer-schur} bounds this row sum by $O_K(1)$.

Suppose $\boldsymbol\delta\ne0$, and let $s=\max\{\nu:\delta_\nu=1\}$. Fix the pair $(p_*,\ell_*)$ selected from $\mathcal C_s$, so that $a_s=p_*$. There are at most $B_*$ possibilities. Write
\[
 a_{<s}=\prod_{\nu<s}a_\nu,\qquad F=\prod_{\nu>s}a_\nu,\qquad
 n=D a_{<s}\in\mathcal N_s.
\]
For every $\nu>s$ we have $\delta_\nu=0$. Using the equality $t=c-a$, we have
\begin{equation*}
 d_*:=\det(\mathbf z,\mathbf z')=P_*-p_*nF,\qquad
 \mathbf z'=h\mathbf z+d_*w_{\mathbf z}.
\end{equation*}
Since $p_*\nmid P_*$, the target in the basis $(\mathbf z,w_{\mathbf z})$ modulo $p_*$ has coordinates $(h,P_*)$ with nonzero second entry. The intersection of $\ell_*$ with the affine line of points $(h,P_*)$ is empty or fixes one residue $h\equiv h_*\pmod{p_*}$. This residue depends only on the source and the selected pair. Let $r_{\mathbf z}=(g w_{\mathbf z})_1/(g\mathbf z)_1$. The first-coordinate box imposes $|h+d_*r_{\mathbf z}|\le16$, whence
\begin{equation}\label{core:eq:reused-phase-condition}
 \biggl\|nF r_{\mathbf z}
       -\frac{P_*r_{\mathbf z}+h_*}{p_*}\biggr\|\le16/p_* .
\end{equation}
If two distinct $n,n'$ satisfy this condition, then $\|(n-n')F r_{\mathbf z}\|\le32/p_*\le100e^{-L^{\vartheta_s}}$.

For $\mathcal F_{\rm exc}=\{F:r_{\mathbf z}\in\mathcal E_s(F)\}$, the separation condition \eqref{core:eq:spacing-weighted-exception} and $\pi_\nu(p)\le(V_\nu p)^{-1}$ imply
\begin{equation*}
 \sum_{F\in\mathcal F_{\rm exc}}\prod_{\nu>s}\pi_\nu(a_\nu)
      \le e^{-L^{\vartheta_s}/4}.
\end{equation*}
For each other $F$, \eqref{core:eq:reused-phase-condition} determines at most one $n$. Now $(a_{<s},P_*)=1$ and $D\mid P_*$, so $D=(n,P_*)$ and $a_{<s}=n/D$ determine all the leaving choices and all lower incoming primes. Thus the target orders contribute $(m_0!\kappa_{J_0})^K=m_0^{-K}$, and the fixed determinant allows $O(1)$ choices of $h$. For $\nu<s$, each choice with $\delta_\nu=0$ contributes $\pi_\nu(a_\nu)$; a prime chosen with $\delta_\nu=1$ occurs at most once in $\mathcal C_\nu$ and has weight $V_\nu^{-1}$. The sum of the weights for $\nu>s$ is at most one. Summing also over the at most $B_*$ choices of $(p_*,\ell_*)$ contributes $O_K(B_*m_0^{-K})$.

For exceptional $F$, allow all leaving choices and target orders, whose number cancels $\kappa_{J_0}^K$. The choices for which $\delta_\nu=1$, including the pair with index $s$, contribute $O_K(B_*^K)$ in total. For each remaining index $\nu<s$ the sum of the weights is at most one. Thus the row sum for $t=c-a$ when $\boldsymbol\delta\ne0$ is $O_K(\varepsilon_K)$, where
\begin{equation*}
 \varepsilon_K=\frac{B_*}{m_0^K}+B_*^K e^{-L^{0.1}/4}.
\end{equation*}
Reversing the transitions exchanges the absolute row and column sums and interchanges $\boldsymbol\epsilon$ and $\boldsymbol\delta$. The column bound is $O_K(1)$ when $\boldsymbol\epsilon=0$ and $O_K(\varepsilon_K)$ otherwise. At least one bound is $O_K(\varepsilon_K)$, and the other is $O_K(1)$. Since $B_*\ll_K L^{1/2}$, $m_0>L^{0.01}$ and $K_0/200>G+1$, we have $\sqrt{\varepsilon_K}\ll_{K,G} L^{1/4-K/200}+L^{-G}\ll_{K,G}L^{-G}$. Schur's test proves the same bound for the operator norm.

\smallskip\noindent

It remains to estimate the major arc term for these choices of $\boldsymbol\epsilon,\boldsymbol\delta$. If $\boldsymbol\delta\ne0$, fix the source, leaving choices, entering primes and all pairs selected from $\mathcal C$. Choose one such pair $(p,\ell)$. Writing $\mathbf z'=h\mathbf z+Dt w_{\mathbf z}$ as above places $(h,t)$ in the centrally symmetric parallelogram
\[
 |h+Dt r_{\mathbf z}|\le16,\qquad |t|\le5Y ,
\]
of area $O(Y)$. Since $p\nmid D$, the inverse image of $\ell$ under $(h,t)\mapsto h\mathbf z+Dt w_{\mathbf z}\pmod p$ is an index-$p$ lattice. In a basis of this lattice the parallelogram has area $O(Y/p)$, and admissible pairs have primitive coordinates in that basis, since a common divisor would divide both entries of $\mathbf z'$. Heath-Brown \cite[Lemma~1(vii)]{HB02} therefore bounds the number of targets by $O(1+Y/p)$. Multiplying by the major-arc measure and summing the label choices as in Lemma~\ref{core:lem:transfer-schur} bounds the row sum by
\[
 O_K\!\bigl(L^{3A_0}(1+B_*)^K(p_0^{-1}+Y^{-1})\bigr)
 \ll_A L^{-A}
\]
for any fixed $A$. If $\boldsymbol\delta=0$, some $\epsilon_\nu=1$ and the same column bound follows by the reverse argument. For the opposite orientation, use the bound for the sum of absolute values from Lemma~\ref{core:lem:transfer-schur} for the contribution from the major arcs. Taking $A$ sufficiently large, we obtain the required norm bound from Schur's test. All estimates for these choices of $\boldsymbol\epsilon,\boldsymbol\delta$ use absolute values, so they hold with the stated multipliers.
\end{proof}

It remains to impose the restriction $\flat$ in the identity of Proposition~\ref{core:prop:exact-fine-transfer} and to replace the expectation over lines by the sum over primitive initial points.

\begin{proof}[Proof of Proposition~\ref{core:prop:path-bound}]
Take $G=E+4$ and choose $A_0,K_0$ as in Lemma~\ref{core:lem:fresh-reused-edges}. Let $t_p$, $p\in\Pset_*$, be independent circle variables with the degree-one Fej\'er density $1+\cos(2\pi t)=|1+e(t)|^2/2$, and put $Z_p=2e(t_p)$. Fourier orthogonality gives
\begin{equation}\label{core:eq:distinct-introduction-moments}
 \mathbb E Z_p^j
 =2^j\int_{\mathbb T} e(jt)
       \bigl(1+\tfrac12 e(t)+\tfrac12 e(-t)\bigr)\,\dd t
 =\begin{cases}1,&j=0,1,\\0,&j\ge2\end{cases}
 \qquad (j\in\mathbb Z_{\ge0}).
\end{equation}
Let $\mathsf T(Z)$ be the sum of the operators in \eqref{core:eq:forward-memory-operator}, with an additional factor $Z_{a_\nu}$ whenever $\delta_\nu=0$. Put
\[
 \mathbf b_Z(\xi)=\mathbf b(\xi)
       \prod_{\nu=1}^K\prod_{p\in\mathbf L_\nu}\overline{Z_p}.
\]
Conjugation in the first argument of the inner product introduces one factor $Z_p$ for each initial entry $p$. For a sequence of transitions, let $j_p$ count the initial entries equal to $p$ and the choices $a_\nu=p$ with $\delta_\nu=0$. Its multiplier is then $\prod_p Z_p^{j_p}$. By independence and \eqref{core:eq:distinct-introduction-moments}, its expectation equals $\ind_{j_p\le1\ \mathrm{for\ all}\ p}$. Proposition~\ref{core:prop:exact-fine-transfer} now implies
\[
 \mathcal M_T(g)=\mathbb E\langle\mathbf b_Z,\mathsf T(Z)^T\mathbf b\rangle .
\]

We have $\|\mathbf b_Z\|\le 2^{Km_0}$ and $\|\mathbf b\|\le1$. For each $\boldsymbol\epsilon,\boldsymbol\delta$, apply Lemma~\ref{core:lem:fresh-reused-edges} with the unit-modulus multiplier $Z_p/2$ and restore a factor $2$ for each index with $\delta_\nu=0$. Summing the $4^K$ choices, we obtain $\|\mathsf T(Z)\|\le 8^K C_{K,G}L^{-G}$ uniformly in $Z$. It follows that
\[
 |\mathcal M_T(g)|
 \le 2^{Km_0}(8^K C_{K,G}L^{-G})^T\ll L^{-ET},
\]
since $K$ is fixed, $m_0=o(T)$ and $G=E+4$. By \eqref{core:eq:closed-walk-average} and Proposition~\ref{core:prop:joint-root-replacement}, the average over primitive initial points is
\[
 \frac1{\zeta(2)}\int_{\mathcal D_\infty}\!
 \mathcal M_T(g(u,v,r))\,\dd u\,\dd v\,\dd r+O(x^{-\delta/20}).
\]
The domain has volume $15$, and $T\log L=o(\log x)$, which proves the asserted mean value estimate.
\end{proof}

\subsection{Exceptional terms}

The term $E_1$ in Proposition~\ref{core:prop:operator-identity} is the sum over edges with squarefree $Dac$ at which separation fails at one or both endpoints.

\begin{lemma}\label{core:lem:single-edge-comparison}
Let $\Gamma_1$ be a family of edges $\gamma=(e_1,\mathbf z';\mathbf L_0,\mathbf L_1)$ with primitive $\mathbf z'$. Assume the label, dyadic and determinant conditions in \eqref{core:eq:uniform-path-family}, with the coordinate bounds holding on a nonempty set of real frames. Suppose that
\[
 |\Gamma_1|,\ \max_{\gamma\in\Gamma_1}|d_\gamma|
       \le\exp(C_1L^{0.21}),
\]
where $C_1$ is fixed and $d_\gamma$ depends only on $\gamma$. At each endpoint choose one of $1$, $\ind_{\rm sep}$ or $1-\ind_{\rm sep}$. Let $K_{\gamma,1}$ be their product on the coordinate boxes and zero outside. Put
\[
 \mathcal F_1(g,\boldsymbol\ell)
  =\sum_{\gamma\in\Gamma_1}d_\gamma K_{\gamma,1}(g)
       \prod_{p\in\mathcal P_{\rm act}(\gamma)}
           \prod_{j\in A_p(\gamma)}
                       \ind_{\ell_p=[\mathbf z_j]_p},
 \qquad \mathbf z_0=e_1,\quad\mathbf z_1=\mathbf z'.
\]
Here \(A_p(\gamma)\subseteq\{0,1\}\) records the two endpoint lists. Then
\begin{equation*}
 \mathcal A_{\rm ar}(\mathcal F_1)=\mathcal A_{\rm mod}(\mathcal F_1)+O_{A,C_1}(L^{-AT})
 \qquad(A>0).
\end{equation*}
\end{lemma}

\begin{proof}
The slope bound on an edge implies $|z'_2|\ll H$ and $|z'_1|\ll H+1$. Replacing a separation indicator by its complement introduces no further boundary points, so by Lemma~\ref{core:lem:joint-root-boundary-geometry} we have $\mathcal N\le\exp(O_K(L^{0.21}))$. All line conditions are modulo primes dividing $S=Dac$, with $\log S\ll_K(J_0+2)L^{0.2}$. Thus the error in Lemma~\ref{core:lem:bounded-modulus} with the choice of $n$ above is $x^{-c_\delta/4+o(1)}$ per edge. Since $\sum_{\gamma\in\Gamma_1}|d_\gamma|=x^{o(1)}$, the total error is $\ll_A L^{-AT}$.
\end{proof}

At the initial point, $\alpha_g(e_1)=r$. It therefore suffices to bound the measure of the values for which the separation conditions fail.

\begin{lemma}\label{core:lem:good-phase-measure}
There is $c_1>0$ such that, uniformly in tuples $\mathbf L$ with $m_0$ distinct entries from each $\Pset_\nu$, the set of $r\in\mathbb T$ failing separation for $\mathbf L$ has measure $O(e^{-c_1L^{0.1}})$.
\end{lemma}
\begin{proof}
For a nonzero integer $h$, multiplication by $h$ preserves circle measure, so $\operatorname{meas}\{r:\|hr\|\le t\}=2t$ for $0\le t\le 1/2$. Summing over $l,D$ bounds the exceptional measure for the first condition by $2m_0^{K}Y^{-1/2}$.

For the second condition, the integer $(n-n')F$ is nonzero whenever $n\ne n'$. Fix $s$ and $F\in\mathcal F_s$, and put $t=100e^{-L^{\vartheta_s}}$. Summing the same identity over distinct pairs $n,n'\in\mathcal N_s$, we obtain $|\mathcal E_s(F)|\le 2t\binom{|\mathcal N_s|}{2}$, as in the first moment estimate of Rudnick \cite[(2), (10)]{Rudnick18}. The lower bound for the gaps in \eqref{core:eq:restricted-band-hyp} implies
\[
 \log|\mathcal N_s|\le K\log m_0+\sum_{\nu<s}\log|\Pset_\nu|
 \le K\log m_0+2K L^{\vartheta_s-1/(20K)}
 =o_K(L^{\vartheta_s}),
\]
uniformly in the allowed exponents, since $L^{\vartheta_s}\ge L^{0.1}$. Averaging the preceding measure bound over $F$, we obtain $\int_{\mathbb T}W_s(r)\,\dd r\ll\exp(-L^{\vartheta_s}+o_K(L^{\vartheta_s}))$. By Markov's inequality, \eqref{core:eq:spacing-weighted-exception} fails on a set of measure $\ll e^{-L^{\vartheta_s}/2}$. Now sum over $s$ and use $\log Y\ge L^{\vartheta_1}-O(1)$ in the first bound.
\end{proof}

For the terms in which prime labels are repeated we count factorizations of $ah+1$ using the following divisor estimate.

\begin{lemma}\label{core:lem:complete-divisors}
Uniformly for $Z\ge2$ and $1\le q\le Z^{1/4}$,
\[
 \sum_{\substack{n\le Z\\n\equiv1\,({\rm mod}\,q)}}d(n)^2
 \ll \frac{Z(\log Z)^3}{q}.
\]
\end{lemma}
\begin{proof}
See \cite[Theorem~2]{Shiu}.
\end{proof}

\begin{proposition}\label{core:prop:exceptional-arithmetic}
Let $K$ and $A_0>0$ be fixed. Under \eqref{core:eq:coefficient-bounds}, the terms $E_1,E_2$ in Proposition~\ref{core:prop:operator-identity} satisfy
\[
 |E_1|+|E_2|\ll_R XYL^{-R}\qquad(R>0).
\]
\end{proposition}
\begin{proof}
We first estimate $E_1$. Let $\mathcal A_Y=\{a\in\Aset(\Pset):\eta(a/Y)\ne0\}$. For $a,c\in\mathcal A_Y$, the absolute value of the kernel in \eqref{core:eq:coarse-kernel} is $\ll_{\eta,\psi}(d_0/D)\lambda_{a,c}(t)$, where
\[
 \lambda_{a,c}(t)=\ind_{t=c-a}
       +L^{3A_0}Y^{-1}\ind_{|t|\le5Y},
 \qquad \sum_t\lambda_{a,c}(t)\ll1+L^{3A_0}.
\]
We bound $|E_1|$ by the sum of the absolute contributions from failure at the source and at the target. For the source, bound the four coefficients by $L^{4C}$ and retain the primitive endpoints, their coordinate boxes and the divisibility conditions. Place the ordered shared primes first in the two corresponding tuples, followed by $c_\nu$ at the source and $a_\nu$ at the target. Each term of \eqref{core:eq:switched-tuple-sum} occurs once. In each dyadic range the source contribution is therefore $\ll L^{4C}(UV/d_0)\mathcal A_{\rm ar}(\mathcal F_1)$, with $1-\ind_{\rm sep}$ at the source, no separation condition at the target and
\[
 d_\gamma=\frac{d_0}{D}\lambda_{a,c}(t)
                  \prod_{\nu=1}^K V_\nu^{-(J_0+2)}.
\]
Here $\mathbf z'=(h,Dt)^t$ with $|h|\ll H$ and $|t|\le5Y$. Each edge has $O_K(J_0)$ selected primes. Thus $|\Gamma_1|$ and $\max|d_\gamma|$ are at most $\exp(O_K(L^{0.21}))$, as required in Lemma~\ref{core:lem:single-edge-comparison}. After averaging over the independent lines, for fixed labels with $Dac$ squarefree, the line conditions have probability
\[
 \prod_{p\mid Dca}\frac1{p+1}\le\frac1{Dca}.
\]
For shared primes this uses $D\mid\det(e_1,\mathbf z')$, so their source and target lines agree.

For fixed $u,r,t$ the first-coordinate condition $1\le u(h+Dt r)\le16$ allows $O(1)$ integers $h$. The separation conditions at the source are evaluated at $r$. By Lemma~\ref{core:lem:good-phase-measure},
\[
 \int_0^1\bigl(1-\ind_{\rm sep}(r,\mathbf L)\bigr)\,\dd r
             \ll e^{-c_1L^{0.1}}.
\]
It remains to sum $\lambda_{a,c}(t)$ over $t$ and the line probabilities over the labels. Dropping the dyadic and coprimality conditions, we have
\[
 \sum_{D,c,a}^{\rm ordered}
    \biggl(\prod_{\nu=1}^{K}V_\nu^{-(J_0+2)}\biggr)\frac1{Dca}
 \le\prod_{\nu=1}^K
       V_\nu^{-(J_0+2)}\biggl(\sum_{p\in\Pset_\nu}\frac1p\biggr)^{J_0+2}
 =1.
\]
Since $UV/d_0=XY$ and $d_0/D\le1$, the contribution from failure at the source in one dyadic range satisfies
\[
 |E_{1,{\rm src}}|
 \ll XY L^{4C}(1+L^{3A_0})e^{-c_1L^{0.1}}
                       +O_A(XYL^{-AT}).
\]
The same bound holds at the target on reversing the edge, since $\lambda_{a,c}(t)=\lambda_{c,a}(-t)$. Summing over the $O_K(L)$ dyadic ranges, we obtain
\[
 |E_1|\ll_R XYL^{-R}.
\]

We next estimate $E_2$. For $a\in\mathcal A_Y$ we have $a\asymp Y=X^{o(1)}$. The arguments $ah+1$ and $ch+1$ are positive for every $h\ge1$. Taking $Z\asymp X$ sufficiently large in Lemma~\ref{core:lem:complete-divisors}, first with $q=a$ and then with $q=1$, we obtain
\begin{equation}\label{core:eq:progression-divisor-square}
 \sum_{h\ll X/Y}d(ah+1)^2\ll \frac XYL^3,
 \qquad
 \sum_{u\ll X}d(u)^2\ll XL^3.
\end{equation}
For fixed $a,c\in\mathcal A_Y$, Cauchy's inequality and \eqref{core:eq:progression-divisor-square} imply
\[
 L^{4C}\sum_{h\ll X/Y}d(ah+1)d(ch+1)
 \ll L^{4C+3}X/Y
\]
for the terms with $mn-1=ah$ and $rs-1=ch$. For the major-arc term put $u=mn$ and $v=rs$, so $u,v\ll X$ and $|cu-av|\le5Y$. Since $a,c\asymp Y$, each $u$ is paired with $O(1)$ values of $v$, and conversely. Applying $2d(u)d(v)\le d(u)^2+d(v)^2$ and \eqref{core:eq:progression-divisor-square} bounds the sum of $d(u)d(v)$ over these pairs by $O(XL^3)$. As $\operatorname{meas}(\Mmajor)\ll L^{3A_0}/Y$, the combined absolute weight for any $\mathcal S\subseteq\mathcal A_Y^2$ is
\[
\ll L^{3A_0+4C+3}\frac XY\,
\mathcal V(\Pset)^{-2}|\mathcal S|.
\]
Since $|\mathcal A_Y|\le4Y$, we have $\mathcal V(\Pset)^{-2}|\mathcal A_Y|^2\ll_K Y^2$. For fixed $a,c$, under the product measure $\prod_{\nu,h}(V_\nu p_{\nu,h})^{-1}$, the probability of a repetition among the $J_0$ shared labels chosen from one set or a common prime factor of $D$ and $ac$ is $O_K(J_0^2p_0^{-1})=O_K(J_0^2e^{-L^{0.1}})$. Indeed, each pairwise coincidence has probability at most $4/p_0$, as does equality with any fixed prime factor of $ac$. Counting pairs for each common prime, we obtain
\[
\#\{(a,c)\in\mathcal A_Y^2:(a,c)>1\}
\le \sum_{p\in\Pset_*}(4Y/p)^2\ll_K Y^2p_0^{-1}
\ll_K Y^2e^{-L^{0.1}}.
\]
Consequently
\[
 |E_2|\ll_K XYL^{3A_0+4C+3}J_0^2e^{-L^{0.1}}
 \ll_R XYL^{-R}.
\]
\end{proof}

\begin{proof}[Proof of Proposition~\ref{core:lem:restricted-determinant}]
Given $R>0$, take $E=R+4C+6$, $G=E+4$, $A_0=G+3$ and $K_0=\lfloor200(G+1)\rfloor+1$ in Lemma~\ref{core:lem:fresh-reused-edges}. For fixed $K\ge K_0$, Proposition~\ref{core:prop:path-bound} supplies \eqref{core:eq:normalized-high-moment} in every dyadic range for $D$. For the same $A_0$, Proposition~\ref{core:prop:exceptional-arithmetic} bounds $|E_1|+|E_2|$ by $O_R(XYL^{-R})$. Since $E>R+4C+1$, the required bound follows from Proposition~\ref{core:prop:conditional-comparison}.
\end{proof}

\section{The major arc estimate}\label{core:sec:typeII-major}

We now estimate $Q_Y(\mathfrak M)$ for the M\"obius coefficients $\alpha_m$ in Theorem~\ref{core:quadratic-theorem}. We combine M\"obius cancellation with a mean square estimate on the intervals where a prime Dirichlet polynomial takes a large value. The estimate holds for every fixed $A_0>0$ and $K\ge1$, so it applies to the arcs chosen in Proposition~\ref{core:lem:restricted-determinant}. Adding that comparison completes the proof of Theorem~\ref{core:quadratic-theorem}.

\begin{proposition}\label{core:lem:typeII-major}
Let the coefficients have the form specified in Theorem~\ref{core:quadratic-theorem}, and let $\mathfrak M$ be defined by \eqref{core:eq:major-arcs-def} for $A_0>0$. For every fixed $K\ge1$ and every $R>0$, we have
\begin{equation*}
  |Q_Y(\mathfrak M)|\ll xYL^{-R}
  \qquad(Y>0),
\end{equation*}
where the implied constant and an ineffective lower threshold for $x$ depend only on $\delta,C,R,A_0,K,\eta$.
\end{proposition}

\subsection{Estimates for Dirichlet polynomials}

\begin{lemma}\label{core:lem:sifted-mobius}
Let $\delta,A,B>0$. Suppose $x^\delta\le M\le x$ and $q\le L^B$, and let $\chi$ be any Dirichlet character modulo $q$. Then, uniformly for $t\in\mathbb R,\ |t|\le x^2$,
\begin{equation*}
 \sum_{\substack{m\le M\\P^-(m)>z}}\mu(m)\chi(m)m^{it}
 \ll M L^{-A},
\end{equation*}
where the implied constant depends on $\delta,A,B$.
\end{lemma}

\begin{proof}
	Let $\chi^*$ be the primitive character modulo $q^*\mid q$ inducing $\chi$.
	Since $q\le L^B<z$ for large $x$, every integer with $P^-(n)>z$ is
	coprime to $q$, and hence $\chi(n)=\chi^*(n)$. Thus, for $\Re s>1$,
	\[
	\begin{aligned}
		F_z(s,\chi)&:=\sum_{\substack{n\ge1\\P^-(n)>z}}
		\frac{\mu(n)\chi(n)}{n^s}
		=\prod_{p>z}\bigl(1-\chi^*(p)p^{-s}\bigr)\\
		&=\frac{1}{L(s,\chi^*)}
		\prod_{p\le z}\bigl(1-\chi^*(p)p^{-s}\bigr)^{-1}.
	\end{aligned}
	\]
	The exclusion of Landau--Siegel zeros
	(see \cite{BanksSiegel,GreenSiegel}), together with \cite[Theorem~24.20]{MVIII}, shows that there exists some absolute $c_0>0$ such that
	\[
	L(s,\chi^*)^{-1} \ll \log q^* + (\log(|u|+3))^{3/4} 
	\qquad \text{for } \sigma\ge1-
	\frac{c_0}{\log q^* + (\log(|u|+3))^{3/4}},
	\]
	where $s=\sigma+iu$.
	For $q^*=1$, the same bound follows from the estimate for
	$\zeta(s)^{-1}$ in \cite[Theorem~24.18]{MVIII}.
	We take $\sigma_0=1-\tfrac14c_0L^{-3/4}$.
	For $q^*\le L^B$,
	$|L(s,\chi^*)^{-1}|\ll_B L^{3/4}$ uniformly for
	$\Re s\ge\sigma_0$ and $|\Im s|\le2x^2$.
	Since $(1-\sigma_0)\log z=o(1)$, Mertens' theorem
	\cite[Theorem~2.7]{MV} shows that, uniformly for $\Re s\ge\sigma_0$,
	\[
	\log\left|\prod_{p\le z}\bigl(1-\chi^*(p)p^{-s}\bigr)^{-1}\right|
	\le\sum_{p\le z}p^{-\sigma_0}+O(1)
	=\log\log z+O(1).
	\]
	Consequently,
	\begin{equation}\label{eq:Fz-bound-strip}
		|F_z(s,\chi)|\ll_B L^{3/4}\log z\ll_B L
		\qquad(\Re s\ge\sigma_0,\ |\Im s|\le2x^2).
	\end{equation}
	
	Let $\sigma_1=1+1/\log M$ and $T=L^{A+3}$. We obtain from the truncated Perron formula \cite[Theorem~7.2]{Koukou} that
	\begin{equation*}
		\sum_{\substack{m\le M\\P^-(m)>z}}\mu(m)\chi(m)m^{it}
		=\frac{1}{2\pi i}\int_{\sigma_1-iT}^{\sigma_1+iT}
		F_z(s-it,\chi)\frac{M^s}{s}\,\dd s
		+O\!\left(\frac{M\log M}{T}+1\right).
	\end{equation*}
	The error is $O_{\delta,A}(ML^{-A})$.
	We shift the contour to $\Re s=\sigma_0$. The function
	$F_z(s-it,\chi)$ is holomorphic throughout the rectangle, where
	\eqref{eq:Fz-bound-strip} applies since
	$|\Im(s-it)|\le T+|t|\le2x^2$ for large $x$.
	The two horizontal integrals and the new vertical integral contribute
	in total
	\[
	\ll_B\frac{ML}{T}+M^{\sigma_0}L\log T
	\ll_{\delta,A,B}ML^{-A},
	\]
	since $M\geq x^\delta$. This completes the proof of Lemma \ref{core:lem:sifted-mobius}.
\end{proof}

\begin{lemma}\label{core:lem:prime-sampling}
Let $\delta,A>0$, and suppose $N\ge x^\delta$ and $\tfrac12L^{0.1}\le\log P\le2L^{0.2}$. For $|v_p|\le1$ and arbitrary complex coefficients $\beta_n$, we consider
\[
  D(t)=\frac1P\sum_{P\le p<2P}v_p p^{it},\qquad
  N_\beta(t)=\frac1N\sum_{N\le n\le2N}\beta_n n^{it}.
\]
We define
\[
 \mathcal I=\Bigl\{[j,j+1]:j\in\mathbb Z,\quad
   \exists\,t\in[j,j+1]\cap[-x^2,x^2]\ \text{with }|D(t)|>L^{-A}\Bigr\}.
\]
Then we have
\begin{equation*}
  |\mathcal I|\le\exp\bigl(O(L^{0.91})\bigr),\qquad
  \sum_{I\in\mathcal I}\sup_{t\in I}|N_\beta(t)|^2
  \ll\frac{L}{N}\sum_{N\le n\le2N}|\beta_n|^2,
\end{equation*}
where the two implied constants depend on $A$ and on $\delta,A$, respectively.
\end{lemma}

\begin{proof}
For each $I=[j,j+1]\in\mathcal I$, we choose $t_I\in I\cap[-x^2,x^2]$ with $|D(t_I)|>L^{-A}$. For $h=0,1$, the points $t_I$ with $j\equiv h\pmod2$ are $1$-separated. We apply \cite[Lemma~6.5]{KMT} to each parity class with $q=1$ and $\alpha=(A+1)\log L/\log P$, and obtain
\[
 |\mathcal I|\le \exp\bigl(O_A(L^{0.9}\log L)\bigr)
 \le \exp\bigl(O_A(L^{0.91})\bigr).
\]

We now prove the mean square bound. We may assume $\mathcal I\ne\varnothing$. For each $I\in\mathcal I$, we choose $u_I\in I$ at which $|N_\beta|$ attains its maximum. We order the intervals by their left endpoints and write $u_1,\ldots,u_{|\mathcal I|}$ for the corresponding points, which satisfy $|u_r|\le x^2+1$. We write $S_N(u)=\sum_{N\le n\le2N}n^{iu}$. For $|u|\le N$, we have by \cite[Theorem~16.13]{MVII}
\[
 S_N(u)=\int_N^{2N}v^{iu}\dd v+O(1) \ll\frac{N}{1+|u|}.
\]
For $N<|u|\le3x^2$, the higher derivative estimate \cite[Theorem~16.12]{MVII}, with a fixed derivative order greater than $2/\delta$, shows that $S_N(u)\ll_\delta Nx^{-c}$ for some $c=c(\delta)>0$, since $N\ge x^\delta$. Thus
\[
 S_N(u) \ll_\delta N\Bigl(\frac1{1+|u|}+x^{-c}\Bigr)
 \qquad(|u|\le3x^2).
\]
By the duality argument in \cite[proof of Theorem~9.6]{IK} and the ordering of the intervals, we obtain
\begin{equation*}
\begin{aligned}
 \sum_{r=1}^{|\mathcal I|}|N_\beta(u_r)|^2
 &\le\frac1{N^2}\max_{1\le r\le|\mathcal I|}
    \sum_{\ell=1}^{|\mathcal I|}|S_N(u_r-u_\ell)|
    \sum_{N\le n\le2N}|\beta_n|^2\\
 &\ll_\delta\frac{\log(2|\mathcal I|)+|\mathcal I|x^{-c}}{N}
    \sum_{N\le n\le2N}|\beta_n|^2
 \ll_{\delta,A}\frac{L}{N}\sum_{N\le n\le2N}|\beta_n|^2,
\end{aligned}
\end{equation*}
where the last step follows from $|\mathcal I|x^{-c}=o(1)$ and $\log(2|\mathcal I|)\ll_A L$.
\end{proof}

\begin{lemma}\label{core:lem:fourier-mellin}
For a finite complex sequence $(w_\ell)$ and $y_\ell\in[1,16]$, we consider
\[
 S(\xi)=\sum_\ell w_\ell e(\xi y_\ell),\qquad
 F(t)=\sum_\ell w_\ell y_\ell^{it}.
\]
For every $T\ge1$,
\[
 \int_{-T}^T|S(\xi)|^2\,\dd\xi
 \ll\int_{-200T}^{200T}|F(t)|^2\,\dd t
       +\frac1T\biggl(\sum_\ell|w_\ell|\biggr)^2.
\]
\end{lemma}

\begin{proof}
Let $g\in C_c^\infty((1/2,20))$ be fixed and equal to $1$ on $[1,16]$, and write
\[
 \widehat g_\xi(t)=\frac1{2\pi}\int_0^\infty
 g(y)e(\xi y)y^{-it}\,\frac{\dd y}{y}.
\]
For $|\xi|\le T$ and $|t|>200T$, we have $|2\pi\xi-t/y|\ge |t|/100$ on the support of $g$. Applying integration by parts twice, we find $\widehat g_\xi(t)\ll|t|^{-2}$. Mellin inversion \cite[(4.106)]{IK} therefore yields
\begin{equation}\label{core:eq:fourier-mellin-truncation}
 S(\xi)=\int_{|t|\le200T}\widehat g_\xi(t)F(t)\,\dd t
       +O\!\biggl(\frac1T\sum_\ell|w_\ell|\biggr)
       \qquad(|\xi|\le T).
\end{equation}
By Fubini's theorem,
\[
 \int_{|t|\le200T}\widehat g_\xi(t)F(t)\,\dd t
 =\frac1{2\pi}\int_0^\infty\frac{g(y)}y
   \biggl(\int_{|t|\le200T}F(t)y^{-it}\,\dd t\biggr)e(\xi y)\,\dd y.
\]
By Plancherel's formula \cite[(4.98)]{IK} with respect to $\xi$, we get
\[
\begin{aligned}
\int_{\mathbb R}\biggl|\int_{|t|\le200T}
	\widehat g_\xi(t)F(t)\,\dd t\biggr|^2\dd\xi &=\frac1{4\pi^2}\int_0^\infty\frac{|g(y)|^2}{y^2}
	\biggl|\int_{|t|\le200T}F(t)y^{-it}\,\dd t\biggr|^2\dd y\\
	&\ll\int_{\mathbb R}\biggl|\int_{|t|\le200T}
	F(t)e^{-itv}\,\dd t\biggr|^2\dd v
	\ll\int_{|t|\le200T}|F(t)|^2\,\dd t,
\end{aligned}
\]
where the first inequality follows from the substitution $v=\log y$ on the compact support of $g$, and the second from Plancherel's identity applied to $F(t)\ind_{|t|\le200T}$. Combining this with \eqref{core:eq:fourier-mellin-truncation} completes the proof.
\end{proof}

\subsection{Proof of Proposition~\ref{core:lem:typeII-major}}
\begin{proof}
	If $\mathcal A(\Pset)\cap[Y,4Y]=\varnothing$, then $Q_Y(\mathfrak M)=0$. By Fourier inversion, we have
	\begin{equation*}
		\begin{aligned}
			Q_Y(\mathfrak M)
			=\mathcal V(\Pset)^{-2}\int_{\mathfrak M}\int_{\mathbb R}
			\widehat\psi(\xi)
			\bigg|\sum_{\substack{c\in\mathcal A(\Pset)\\m,n}}
			\eta(c/Y)\alpha_m\beta_n
			e\!\biggl(\theta(cmn-c)+\frac{\xi cmn}{Y}\biggr)\bigg|^2
			\,\dd\xi\,\dd\theta.
		\end{aligned}
	\end{equation*}
	We parametrize the arcs in \eqref{core:eq:major-arcs-def} by $\theta=u/q+\lambda/Y$. Since $(ucmn,q)=1$ on the support of the sum, we use the character expansion \cite[(3.11)]{IK}
	\begin{equation*}
		e_q(ucmn)=\frac1{\varphi(q)}
		\sum_{\chi\bmod q}\tau(\overline\chi)\chi(ucmn),
	\end{equation*}
	where $\tau(\chi)=\sum_{a\bmod q}\chi(a)e_q(a)$ is the Gauss sum, and $|\tau(\chi)|\le\sqrt q$ for every character modulo $q$ \cite[Lemma~3.1]{IK}. For fixed $u,q,\chi,\lambda$, we write
	\[
	S_\lambda(\xi)
	=\frac{1}{YMN\mathcal V(\Pset)}
	\sum_{\substack{c\in\mathcal A(\Pset)\\m,n}}
	\eta(c/Y)\alpha_m\beta_n\chi(cmn)e_q(-uc)e(-\lambda c/Y)e(\xi cmn/Y).
	\]
	We apply Cauchy's inequality to the character sum. After summation over the major arcs, we obtain
	\begin{equation}\label{core:eq:major-Q-character-sum}
		\frac{|Q_Y(\mathfrak M)|}{xY}
		\ll x\sum_{q\le L^{A_0}}\frac{q}{\varphi(q)}
		\sum_{\substack{u\bmod q\\(u,q)=1}}\sum_{\chi\bmod q}\int_{-2L^{A_0}}^{2L^{A_0}}\int_{\mathbb R}
		|\widehat\psi(\xi)|\,|S_\lambda(\xi+\lambda)|^2
		\,\dd\xi\,\dd\lambda.
	\end{equation}
	
	We next reduce the inner integral in \eqref{core:eq:major-Q-character-sum} to a mean square of the Dirichlet polynomial $A_\lambda(t)=B_\lambda(t)M_\alpha(t)N_\beta(t)$, where
	\begin{align*}
		B_\lambda(t)&=\frac{1}{Y\mathcal V(\Pset)}
		\sum_{c\in\mathcal A(\Pset)}
		\eta(c/Y)e_q(-uc)e(-\lambda c/Y)\chi(c)c^{it},\\
		M_\alpha(t)&=\frac1M\sum_m\alpha_m\chi(m)m^{it},\qquad
		N_\beta(t)=\frac1N\sum_n\beta_n\chi(n)n^{it}.
	\end{align*}
	Since $\sum_{c\in\mathcal A(\Pset)}c^{-1}=\mathcal V(\Pset)$ and $c\le4Y$ on the support of $\eta(c/Y)$, the sum of the absolute values of the coefficients of $B_\lambda$ is $O(1)$, and that of $A_\lambda$ is $\ll L^C$.

	We split the $\xi$-integral at $|\xi|=L^{A_0+1}$. The coefficients of $M_\alpha N_\beta$, supported on $\ell=mn\asymp x$, have modulus $\ll L^C d(\ell)/x$. By the definition of $S_\lambda$, Cauchy's inequality in $c$ and direct integration, we have
	\[
	 \int_j^{j+1}|S_\lambda(\xi+\lambda)|^2\,\dd\xi
	 \ll\frac1{(MN)^2}\sum_{\ell\asymp x}
	       \biggl|\sum_{mn=\ell}\alpha_m\beta_n\chi(mn)\biggr|^2
	 \ll\frac{L^{2C+3}}x,
	\]
	for any $j\in \mathbb Z$, where the last bound follows from \eqref{core:eq:progression-divisor-square}. Thus, we use the rapid decay of $\widehat\psi$ to obtain
	\begin{equation}\label{core:eq:typeII-xi-tail}
	\begin{aligned}
	 \int_{|\xi|>L^{A_0+1}}|\widehat\psi(\xi)|
	       |S_\lambda(\xi+\lambda)|^2\,\dd\xi
	 &\ll\frac{L^{2C+3}}x
	       \sum_{\substack{j\in\mathbb Z\\|j|>L^{A_0+1}-1}}
	       \sup_{\xi\in[j,j+1]}|\widehat\psi(\xi)|\\
	     & \ll x^{-1}L^{-R-4A_0-1}
	\end{aligned}
	\end{equation}
for any $R>0$.

	We apply Lemma~\ref{core:lem:fourier-mellin} with $T=2MNL^{A_0+1}$ to $S_\lambda(\xi'/(MN))$, since $1\le cmn/(YMN)\le16$ on the support of the sum. The function $F$ in that lemma is $(YMN)^{-it}A_\lambda(t)$. Since $\widehat\psi$ is bounded, we obtain
	\begin{equation}\label{core:eq:typeII-xi-main}
	\begin{aligned}
	 \int_{|\xi|\le L^{A_0+1}}|\widehat\psi(\xi)|
	       |S_\lambda(\xi+\lambda)|^2\,\dd\xi
	 &\ll\frac1{MN}\int_{-T}^T
	       |S_\lambda(\xi'/(MN))|^2\,\dd\xi'\\
	 &\ll\frac1{MN}\int_{-200T}^{200T}|A_\lambda(t)|^2\,\dd t
	       +\frac{L^{2C}}{MN\,T}.
	\end{aligned}
	\end{equation}
	For large $x$, the $t$-interval in \eqref{core:eq:typeII-xi-main} is contained in $\mathcal T=[-xL^{A_0+2},xL^{A_0+2}]$. Combining  \eqref{core:eq:typeII-xi-tail} with \eqref{core:eq:typeII-xi-main}, we get
	\begin{equation*}
		\int_{\mathbb R}|\widehat\psi(\xi)|\,|S_\lambda(\xi+\lambda)|^2\,\dd\xi
		\ll \frac1x\int_{\mathcal T}|A_\lambda(t)|^2\,\dd t
		      +x^{-1}L^{-R-4A_0-1}.
	\end{equation*}
	We insert this bound into \eqref{core:eq:major-Q-character-sum} and obtain
	\begin{equation}\label{core:eq:Q-chi-cauchy}
		\frac{|Q_Y(\mathfrak M)|}{xY}
		\ll L^{4A_0}\sup_{u,q,\chi,\lambda}
		\int_{\mathcal T}|A_\lambda(t)|^2\,\dd t+O(L^{-R-1}),
	\end{equation}
	where the supremum is over the ranges of $u,q,\chi,\lambda$ in that sum.
	
	We reduce the mean square in \eqref{core:eq:Q-chi-cauchy} to a sum over primes as in Lemma~\ref{core:lem:prime-sampling},
	\[
	 D(t)=\frac1P\sum_{P\le p<2P}v_p p^{it},\qquad |v_p|\le1.
	\]
	We write $c=pc'$, with $p\in\Pset_1$. For fixed $c'$, only $Y/c'\le p<4Y/c'$ contributes. We split the $p$-sum at $2Y/c'$. After multiplication by $c'/Y$, the first part is a polynomial from Lemma~\ref{core:lem:prime-sampling} and the second is twice another such polynomial. The required coefficient bound follows from $0\le\eta\le1$, and \eqref{core:eq:restricted-band-hyp} ensures that each nonzero polynomial has the required scale. Since $\mathcal V(\Pset)^{-1}\sum_{c'}c'^{-1}=V_1^{-1}\ll1$, we have by Cauchy's inequality
	\begin{equation}\label{core:eq:typeII-small-prime-factorization}
		\int_{\mathcal T}|A_\lambda(t)|^2\,\dd t
		\ll\sup_D\int_{\mathcal T}
		|D(t)M_\alpha(t)N_\beta(t)|^2\,\dd t,
	\end{equation}
	where $D$ ranges over the finitely many polynomials just obtained.
	
	For a fixed $A>0$, to be chosen at the end, we split the integral according as $|D(t)|\le L^{-A}$ or $|D(t)|>L^{-A}$. For the first part, by the mean value theorem for Dirichlet polynomials \cite[Theorem~6.1]{Mont71} and the coefficient bound for $M_\alpha N_\beta$ above, we have
	\begin{equation*}
		\int_{\mathcal T}|M_\alpha(t)N_\beta(t)|^2\,\dd t
		\ll L^{2C+A_0+5}.
	\end{equation*}

	For the second part, we use Lemmas~\ref{core:lem:sifted-mobius} and \ref{core:lem:prime-sampling}. For large $x$, we have $2M\le x$ and $|t+v|\le x^2$ on $\mathcal T$. Lemma~\ref{core:lem:sifted-mobius} yields
	\begin{equation*}
		\sup_{t\in\mathcal T}|M_\alpha(t)| \ll_A L^{-A}.
	\end{equation*}
	Let $\mathcal I$ be the interval family in Lemma~\ref{core:lem:prime-sampling} for this $D$ and $A$. Every $t\in\mathcal T$ with $|D(t)|>L^{-A}$ lies in an interval $I_0\in\mathcal I$, since $\mathcal T\subset[-x^2,x^2]$ for large $x$. Applying the second bound in the lemma to $\beta_n\chi(n)$, we obtain
	\[
	\sum_{I_0\in\mathcal I}\sup_{t\in I_0}|N_\beta(t)|^2
	\ll\frac{L}{N}\sum_{N\le n\le2N}|\beta_n\chi(n)|^2
	\ll L^{2C+1}.
	\]
	The estimates for the two parts, together with $|D|\ll1$, imply
	\begin{equation}\label{core:eq:typeII-major-energy}
		\begin{aligned}
			\int_{\mathcal T}
			|D(t)M_\alpha(t)N_\beta(t)|^2\,\dd t
			&\ll L^{-2A}\int_{\mathcal T}|M_\alpha(t)N_\beta(t)|^2\,\dd t+L^{-2A}\sum_{I_0\in\mathcal I}\sup_{t\in I_0}|N_\beta(t)|^2\\
			&\ll L^{2C+A_0+5-2A}.
		\end{aligned}
	\end{equation}
Combining \eqref{core:eq:typeII-small-prime-factorization} and \eqref{core:eq:typeII-major-energy}, we obtain
\begin{equation*}
	\sup_{u,q,\chi,\lambda}
	\int_{\mathcal T}|A_\lambda(t)|^2\,\dd t
	\ll L^{2C+A_0+5-2A}.
\end{equation*}
Inserting this into \eqref{core:eq:Q-chi-cauchy} and choosing $2A>R+2C+5A_0+5$, we obtain the required bound.

\end{proof}

\subsection{Proof of Theorem~\ref{core:quadratic-theorem}}
\begin{proof}
We may replace $\delta$ by $\min(\delta,1/4)$. Choose $A_0$ and $K_0$ from Proposition~\ref{core:lem:restricted-determinant} with the desired $R$. Proposition~\ref{core:lem:typeII-major} applies to the same major arcs and the same prime sets at the given scale $x$. Since $Q_Y$ is nonnegative,
\[
 Q_Y\le |Q_Y-Q_Y(\mathfrak M)|+|Q_Y(\mathfrak M)|
       \ll XYL^{-R}.
\]
Both propositions are uniform in the exponents and prime sets under \eqref{core:eq:restricted-band-hyp}, with the stated dependence on the fixed parameters. This proves Theorem~\ref{core:quadratic-theorem}.
\end{proof}

\section{Removal of the auxiliary weights}\label{mot:sec:application}
We combine the Type~I asymptotic with the estimate for $T_3$ to obtain a weighted asymptotic formula, and then remove the auxiliary weights. Throughout this section $K$ is fixed, $J$ lies in \eqref{common:eq:family-range}, and the prime sets are defined by \eqref{common:eq:family-exponents} and \eqref{common:eq:application-bands}. All estimates below are uniform in $J$.

\subsection{A weighted asymptotic formula}
\begin{proposition}\label{mot:prop:weighted}
For sufficiently large fixed $K$ and every $\Phi\in C_c^\infty((0,\infty))$ supported in $[1/2,1]$, we have
\[
 \Theta_{b,\W}(x,\Phi)
 =\frac{x}{\sqrt L}\biggl(\kappa\Spr
       \int_0^\infty\Phi(t)\,\dd t
       +O_K\bigl(\mathcal N_8(\Phi)L^{-1/2}\bigr)\biggr).
\]
The formula holds for all such $\Phi$ once $x$ is sufficiently large in terms of $K$ alone.
\end{proposition}

\begin{proof}
Combine \eqref{common:eq:vaughan}, Proposition~\ref{common:prop:typeI} and Corollary~\ref{common:cor:T3} with $A=1$. The combined error is $O_K\bigl(\mathcal N_8(\Phi)xL^{-1}\bigr)$.
\end{proof}

\subsection{Averaging the weights}

Following \cite[Section~6]{PD} and \cite[Proposition~2.3 and Section~5.3]{LWW}, we average the weights $\W_j$ over the $J$ families of Section~\ref{sec:common}. To bound the squared difference between this average and $1$ over primes $x/2<p\le x$ with $b(p-1)=1$, we use a Selberg upper-bound sieve of dimension $3/2$ \cite[Section~7.1]{FI}. Cauchy's inequality will then remove the averaged weight.

Put $P_{\mathcal E}=\prod_{p\in\mathcal E}p$. For a function $\mathscr F(n)$ depending only on the indicators $(\ind_{q\mid n})_{q\in\mathcal E}$, define its average over residue classes by
\[
 \mathbb E_{b,\mathcal E}\bigl(\mathscr F\bigr)
 =\frac{1}{\varphi(P_{\mathcal E})}
 \sum_{\substack{a\bmod P_{\mathcal E}\\ (a+1,P_{\mathcal E})=1}}\mathscr F(a).
\]
For each $q\in\mathcal E$, the Chinese remainder theorem shows that the $q-1$ classes modulo $q$ other than $-1$ occur equally often in this sum. Thus $\mathbb E_{b,\mathcal E}\bigl(\ind_{q\mid a}\bigr)=1/(q-1)$. Since the bands are pairwise disjoint, the corresponding counts $\omega_{j,\nu}$ are independent under $\mathbb E_{b,\mathcal E}$. For each band, write $\mu_{j,\nu}=\sum_{p\in\Pset_{j,\nu}}1/(p-1)$. Then
\[
 \mathbb E_{b,\mathcal E}\omega_{j,\nu}=\mu_{j,\nu},\qquad
 \mathbb E_{b,\mathcal E}\omega_{j,\nu}^2
 =\mu_{j,\nu}^2+
       \sum_{p\in\Pset_{j,\nu}}\frac1{p-1}
                         \Bigl(1-\frac1{p-1}\Bigr)
 \le\mu_{j,\nu}^2+\mu_{j,\nu}.
\]
Since $\mu_{j,\nu}\asymp1$, uniformly in $j,\nu$, we deduce from \eqref{mot:eq:mean-W} that
\[
 \mathbb E_{b,\mathcal E}\W_j
 =\prod_{\nu=1}^K\frac{\mu_{j,\nu}}{V_{j,\nu}}
 =1+O_K\bigl(\exp(-L^{0.1})\bigr),\qquad
 \mathbb E_{b,\mathcal E}\W_j^2\ll_K1.
\]
Define
\begin{equation}\label{mot:eq:ZJ-def}
 Z^{(b)}_{J}(n)
 =\biggl(\frac1{J}\sum_{j=1}^{J}\W_j(n)-1\biggr)^2.
\end{equation}
For $i\ne j$, independence implies $\mathbb E_{b,\mathcal E}((\W_i-1)(\W_j-1)) =O_K(\exp(-2L^{0.1}))$. Consequently,
\begin{equation}\label{mot:eq:variance}
 \mathbb E_{b,\mathcal E}Z^{(b)}_{J}
 =\frac1{J^2}\sum_{j=1}^{J}
       \mathbb E_{b,\mathcal E}(\W_j-1)^2
       +O_K\bigl(\exp(-2L^{0.1})\bigr)
 \ll_K\frac1{J}+\exp(-2L^{0.1}),
\end{equation}
with an implied constant independent of $J$.

\subsection{An upper bound for the weighted sum over primes}
\begin{lemma}\label{mot:lem:sparse-upper}
	We have
	\begin{equation*}
		\sum_{x/2<p\le x} b(p-1)Z^{(b)}_{J}(p-1)
		\le C_K\frac{x}{JL^{3/2}}+O_K(xL^{-10}).
	\end{equation*}
\end{lemma}

Let $z_0=x^{1/100}$ and $D_0=z_0^5=x^{1/20}$, so that $\mathcal E\subset[3,z_0]$ for large $x$. If $n+1\in(x/2,x]$ is prime, then $p\nmid n+1$ for every prime $p\le z_0$. The condition $b(n)=1$ also excludes $p\parallel n$ when $p\equiv3\modf$. The corresponding excluded residue classes are
\[
 \Omega_p=
 \begin{cases}
  \{a\bmod p^2:p\mid a+1\}, 
    &\text{if}\; p=2\text{ or }p\equiv1\modf,\\
  \{a\bmod p^2:p\mid a+1\}\cup\{a\bmod p^2:p\parallel a\},
    & \text{if} \; p\equiv3\modf.
 \end{cases}
\]
The two sets in the second case are disjoint and contain $p$ and $p-1$ classes, respectively. Define the multiplicative function $g_b$ on the divisors of $P(z_0)$ by
\begin{equation}\label{mot:eq:terminal-local-densities}
	g_b(p) = \frac{|\Omega_p|}{p^2} = \begin{cases}
		1/p, & p=2\text{ or }p\equiv1\modf,\\
		2/p - 1/p^2, & p\equiv 3\modf.
	\end{cases}
\end{equation}
For $r\mid P(z_0)$, we define
\[
 \Omega_r=\{a\bmod r^2:a\bmod p^2\in\Omega_p
                  \text{ for every }p\mid r\}.
\]
By the Chinese remainder theorem, $|\Omega_r|=r^2g_b(r)$. We write $n\in\Omega_r$ for $n\bmod r^2\in\Omega_r$.

The next lemma provides the remainder estimate needed for Selberg's sieve applied to $Z^{(b)}_{J}$.

\begin{lemma}\label{mot:lem:weighted-local-counting}
	For $r\mid P(z_0)$, we have
	\begin{equation*}
		\sum_{\substack{x/2<n+1\le x\\n\in\Omega_r}} Z^{(b)}_{J}(n)
		= \frac{x}{2}g_b(r)\mathbb E_{b,\mathcal E}\bigl(Z^{(b)}_{J}\bigr) + R(r),
	\end{equation*}
	where, for every fixed $B>0$,
	\begin{equation*}
		\sum_{\substack{d_1,d_2\le D_0\\ d_1,d_2\mid P(z_0)}}
		\bigl|R([d_1,d_2])\bigr|
		=O_{B,K}(xL^{-B}).
	\end{equation*}
\end{lemma}

\begin{proof}
Expanding the square in \eqref{mot:eq:ZJ-def} and then the products of prime sums in \eqref{common:eq:weight}, we obtain the exact identity
\begin{equation}\label{mot:eq:Z-signed-expansion}
 Z^{(b)}_{J}(n)
 =\sum_{\substack{Q\mid P_{\mathcal E}\\\omega(Q)\le2K}}
       c_Q\ind_{Q\mid n}.
\end{equation}
In a squared weight the two selected primes in one band may coincide; their divisibility conditions combine by taking the least common multiple. Thus every $Q$ is squarefree and
\[
 Q\le\exp(4KL^{0.2}),\qquad
 \sum_Q|c_Q|\le\exp(O_K(L^{0.2})).
\]
Indeed, the expansion of the square contains a constant, an average of single weights and an average of products of two weights. Before collecting terms, these scalar coefficients have total absolute sum $1+2+1=4$, and each weight contains $K$ prime selections. Since $V_{j,\nu}\asymp1$ and $\sum_{p\in\Pset_{j,\nu}}(p-1)^{-1}\asymp1$, the coefficients in this expansion satisfy
\begin{equation}\label{mot:eq:coefficient-reciprocal-bounds}
 \sum_Q\frac{|c_Q|}{\varphi(Q)}\ll_K1,
 \qquad
 \sum_{\substack{Q\\q\mid Q}}\frac{|c_Q|}{\varphi(Q)}
       \ll_K\frac1q\quad(q\in\mathcal E).
\end{equation}
For the first estimate, apply the second moment calculation above to $1$, $\W_j$ and $\W_i\W_j$ before collecting coefficients. For the second, prescribe a selected prime to be $q$. In a repeated band $\Pset_{j,\nu}$, the total contribution containing $q$ is
\[
 \frac1{q-1}\biggl(1+2\sum_{\substack{p\in\Pset_{j,\nu}\\p\ne q}}
                       \frac1{p-1}\biggr)\ll\frac1q.
\]
The reciprocal sums over the other bands are bounded. Averaging the scalar coefficients preserves these bounds independently of $J$. By the Chinese remainder theorem,
\[
 \mathbb E_{b,\mathcal E}Z^{(b)}_{J}
 =\sum_Q\frac{c_Q}{\varphi(Q)}.
\]
For $q\in\mathcal E$, the condition $n\in\Omega_q$ means $q\mid n+1$. Hence $Q\mid n$ is incompatible with $n\in\Omega_r$ when $(Q,r)>1$. Otherwise the two conditions specify $|\Omega_r|$ classes modulo $Qr^2$. Therefore
\[
 \sum_{\substack{x/2<n+1\le x\\n\in\Omega_r}}\ind_{Q\mid n}
 =\frac{x}{2}g_b(r)\frac{\ind_{(Q,r)=1}}Q+O(|\Omega_r|).
\]
Using
\[
 \Bigl|\frac{\ind_{(Q,r)=1}}Q-\frac1{\varphi(Q)}\Bigr|
 \le\frac1{\varphi(Q)}
       \biggl(\sum_{q\mid(Q,r)}1+\sum_{q\mid Q}\frac1q\biggr),
\]
we sum over $Q$ in \eqref{mot:eq:Z-signed-expansion} and apply \eqref{mot:eq:coefficient-reciprocal-bounds}. The stated identity follows, with
\begin{equation}\label{mot:eq:weighted-local-pointwise-error}
 R(r)\ll_K xg_b(r)
       \biggl(\sum_{q\in\mathcal E}\frac1{q^2}
             +\sum_{\substack{q\mid r\\q\in\mathcal E}}\frac1q\biggr)
       +|\Omega_r|\exp(O_K(L^{0.2})).
\end{equation}
We sum \eqref{mot:eq:weighted-local-pointwise-error} over $d_1,d_2\le D_0$. Each squarefree $r=[d_1,d_2]\le D_0^2=x^{1/10}$ arises from at most $3^{\omega(r)}$ pairs. By Mertens' theorem \cite[Theorem~2.7(e) and Corollary~4.12(d)]{MV},
\[
 \sum_{r\mid P(z_0)}3^{\omega(r)}g_b(r)
 =\prod_{p\le z_0}(1+3g_b(p))\ll L^{9/2}.
\]
Moreover, $g_b(q)=1/q$ for $q\in\mathcal E$, so
\[
 \sum_{r\mid P(z_0)}3^{\omega(r)}g_b(r)
       \sum_{\substack{q\mid r\\q\in\mathcal E}}\frac1q
 \ll L^{9/2}\sum_{q\in\mathcal E}\frac1{q^2}.
\]
Since every $q\in\mathcal E$ exceeds $\exp(L^{0.1})$, the latter prime sum is $O(\exp(-L^{0.1}))$. For squarefree $r$, we have $3^{\omega(r)}=\sum_{m\mid r}d(m)$, whence $\sum_{r\le y}\mu^2(r)3^{\omega(r)}\le y\sum_{m\le y}d(m)/m\ll y(\log y)^2$. Together with $|\Omega_r|\le r^2$, this implies
\[
 \sum_{\substack{d_1,d_2\le D_0\\d_1,d_2\mid P(z_0)}}
       |R([d_1,d_2])|
 \ll_K xL^{9/2}\exp(-L^{0.1})
       +x^{3/10}L^2\exp(C_KL^{0.2}).
\]
Both terms are $O_{B,K}(xL^{-B})$ for every fixed $B>0$. The coefficient bounds above are preserved by the normalized averages over the $J$ families, so the estimate is uniform throughout \eqref{common:eq:family-range}.
\end{proof}

\begin{proof}[Proof of Lemma~\ref{mot:lem:sparse-upper}]
	Let
	\[
	 G=\sum_{\substack{d\le D_0\\d\mid P(z_0)}}
	      \prod_{p\mid d}\frac{g_b(p)}{1-g_b(p)}.
	\]
	By the Selberg sieve \cite[Section~6.5]{IK}, there are real weights $\xi_d$, supported on $d\le D_0$ with $d\mid P(z_0)$, such that $\xi_1=1$, $|\xi_d|\le1$, and
	\begin{equation}\label{mot:eq:selberg-quadratic-identity}
		\sum_{\substack{d_1,d_2\le D_0\\d_1,d_2\mid P(z_0)}}
		\xi_{d_1}\xi_{d_2}g_b([d_1,d_2])
		= G^{-1}.
	\end{equation}
	
	For each prime counted on the left below, $n=p-1$ avoids every $\Omega_q$ with $q\le z_0$, so the inner sum on the right equals $1$. Since $Z^{(b)}_{J}\ge0$, we have
	\begin{equation*}
		\sum_{x/2<p\le x} b(p-1)Z^{(b)}_{J}(p-1)
		\le \sum_{x/2<n+1\le x} Z^{(b)}_{J}(n)
		\bigg(\sum_{\substack{d\le D_0\\d\mid P(z_0)}}
		\xi_d\,\ind_{n\in\Omega_d}\bigg)^2.
	\end{equation*}
	The conditions $n\in\Omega_{d_1}$ and $n\in\Omega_{d_2}$ are together equivalent to $n\in\Omega_{[d_1,d_2]}$. We expand the square and apply Lemma~\ref{mot:lem:weighted-local-counting} with $B=10$, together with \eqref{mot:eq:selberg-quadratic-identity}. Since $|\xi_d|\le1$, we obtain
	\begin{equation}\label{mot:eq:post-local-count-sieve}
		\sum_{x/2<p\le x} b(p-1)Z^{(b)}_{J}(p-1)
		\le \frac{x}{2G}\mathbb E_{b,\mathcal E}\bigl(Z^{(b)}_{J}\bigr)
		+O_K(xL^{-10}).
	\end{equation}
	
	To bound $G$ from below, we use $1-\log d/\log D_0\le\ind_{d\le D_0}$ and obtain
	\[
	\begin{aligned}
	 G&\ge\sum_{d\mid P(z_0)}
	       \Bigl(1-\frac{\log d}{\log D_0}\Bigr)
	       \prod_{p\mid d}\frac{g_b(p)}{1-g_b(p)}\\
	  &=\prod_{p\le z_0}(1-g_b(p))^{-1}
	       \bigg(1-\frac1{\log D_0}\sum_{p\le z_0}g_b(p)\log p\bigg)
	   \gg L^{3/2}.
	\end{aligned}
	\]
	Here $\log D_0=5\log z_0$. By \eqref{mot:eq:terminal-local-densities} and Mertens' theorem, we have
	\[
	\begin{aligned}
	 \sum_{p\le z_0}g_b(p)\log p
	 &=\Bigl(\frac32+o(1)\Bigr)\log z_0,\\
	 \prod_{p\le z_0}(1-g_b(p))
	 &=\prod_{p\le z_0}(1-p^{-1})
	   \prod_{\substack{p\le z_0\\p\equiv3\modf}}(1-p^{-1})
	\asymp L^{-3/2}.
	\end{aligned}
	\]
	Thus $G^{-1}\ll L^{-3/2}$ with an absolute implied constant. The lemma follows from \eqref{mot:eq:post-local-count-sieve} and \eqref{mot:eq:variance}.
\end{proof}

\subsection{Proof of Theorem~\ref{mot:thm:main}}
\begin{proof}
Fix $K$ sufficiently large for Proposition~\ref{mot:prop:weighted}. Let $J$ satisfy \eqref{common:eq:family-range}, and let $\Phi\in C_c^\infty((0,\infty))$ be nonnegative and supported in $[1/2,1]$. Apply Proposition~\ref{mot:prop:weighted} to every family $\Pset^{(j)}$. The prime powers $n=p^r$ with $r\ge2$ contribute $O_K(\|\Phi\|_\infty x^{3/5})$, since $\W(n)\ll_Kd(n)\ll_Kn^{1/20}$. On the support of $\Phi(p/x)$, we have $\log p=L+O(1)$ and $\log p\ge L/2$. By positivity, the weighted prime sum is at most $2\Theta_{b,\W_j}(x,\Phi)/L$. Thus, after division by $L$, replacing $\log p/L$ by $1$ costs $O_K(xL^{-5/2}(\|\Phi\|_\infty+\mathcal N_8(\Phi)L^{-1/2}))$. Averaging over $j$, we obtain
\[
 \frac1J\sum_{j=1}^J\sum_p b(p-1)\W_j(p-1)\Phi(p/x)
 =\frac{x}{L^{3/2}}\biggl(\kappa\Spr\int_0^\infty\Phi(t)\,\dd t
       +O_K\bigl(\mathcal N_8(\Phi)L^{-1/2}\bigr)\biggr).
\]
By Cauchy's inequality, \eqref{eq:iwaniec-order} and Lemma~\ref{mot:lem:sparse-upper},
\begin{align}
 &\biggl|\sum_p b(p-1)\Phi(p/x)
       \biggl(\frac1J\sum_{j=1}^J\W_j(p-1)-1\biggr)\biggr|\notag\\
 &\quad\le\|\Phi\|_\infty
       \biggl(\sum_{x/2<p\le x}b(p-1)\biggr)^{1/2}
       \biggl(\sum_{x/2<p\le x}b(p-1)Z^{(b)}_J(p-1)\biggr)^{1/2}\notag\\
 &\quad\ll_K\|\Phi\|_\infty\frac{x}{L^{3/2}}
       \bigl(J^{-1/2}+L^{-17/4}\bigr).\label{mot:eq:remove-weights-CS}
\end{align}
Consequently, uniformly for $J$ in \eqref{common:eq:family-range},
\[
 \biggl|\frac{L^{3/2}}x\sum_p b(p-1)\Phi(p/x)
       -\kappa\Spr\int_0^\infty\Phi(t)\,\dd t\biggr|
 \ll_K\|\Phi\|_\infty J^{-1/2}
       +\mathcal N_8(\Phi)L^{-1/2}.
\]
For sufficiently large $x$, choose $J=\lfloor\log L/(20K\log2)\rfloor$. Then
\begin{equation}\label{mot:eq:smoothed-prime-asymp}
 \sum_p b(p-1)\Phi(p/x)
 =\frac{x}{L^{3/2}}\biggl(\kappa\Spr\int_0^\infty\Phi(t)\,\dd t
       +O_K\biggl(\frac{\|\Phi\|_\infty}{\sqrt{\log L}}
                  +\mathcal N_8(\Phi)L^{-1/2}\biggr)\biggr).
\end{equation}
In particular, for fixed $\Phi$ the error in parentheses is $O_{K,\Phi}((\log L)^{-1/2})$.

The estimate \eqref{mot:eq:smoothed-prime-asymp} extends to nonnegative smooth $\Phi$ supported in $[1/3,4/3]$. Indeed, use a fixed nonnegative smooth partition of unity whose pieces are supported in intervals $(c/2,c)$, with $c>0$ fixed, and apply \eqref{mot:eq:smoothed-prime-asymp} to each rescaled piece at the scale $cx$. At each scale the prime sets and $J$ are chosen anew. The number of pieces is fixed, their seminorms are bounded by a constant times $\mathcal N_8(\Phi)$, and $\log(cx)=L+O(1)$, so the resulting bound has the same form.

Put $\delta=(\log L)^{-1/2}$. Choose smooth functions $\Phi^-_\delta,\Phi^+_\delta$, supported in $[1/3,4/3]$, such that
\[
\begin{gathered}
 0\le\Phi^-_\delta\le\ind_{(1/2,1]}
       \le\Phi^+_\delta\le1,\\
 \int_0^\infty\Phi^\pm_\delta(t)\,\dd t
       =\frac12+O(\delta),\qquad
 \mathcal N_8(\Phi^\pm_\delta)\ll\delta^{-8}.
\end{gathered}
\]
Since $\delta^{-8}L^{-1/2}\ll\delta$, the upper and lower bounds from \eqref{mot:eq:smoothed-prime-asymp} imply
\begin{equation}\label{mot:eq:dyadic-prime-quant}
 \sum_{x/2<p\le x}b(p-1)
 =\frac{\kappa\Spr x}{2L^{3/2}}
       +O_K\biggl(\frac{x}{L^{3/2}\sqrt{\log L}}\biggr).
\end{equation}
Sum \eqref{mot:eq:dyadic-prime-quant} over the scales $y=x2^{-h}\ge\sqrt x$. The errors total $O_K(xL^{-3/2}(\log L)^{-1/2})$. The main terms sum to $\kappa\Spr xL^{-3/2}+O(xL^{-5/2}+\sqrt x)$, since $\sum_{h\ge0}2^{-h-1}=1$ and $\sum_{h\ge0}h2^{-h}\ll1$. The remaining primes contribute $O(\sqrt x)$, which is negligible. This proves Theorem~\ref{mot:thm:main}.
\end{proof}

\noindent\textbf{AI use.} This manuscript presents results obtained through an interactive collaboration between AI and human researchers. Convinced by the mathematical substance of the work, we share this paper to express our excitement about this methodology and its broader potential.

\raggedbottom


\begin{thebibliography}{10}

\bibitem{BanksSiegel}
\bgroup{}W.~D. Banks\egroup{}, {On the OpenAI proof that Landau--Siegel zeros do not exist}, arXiv:2610.15483v1, 2026.

\bibitem{BocaAngles}
\bgroup{}F.~P. Boca\egroup{}, {Distribution of angles between geodesic rays associated with hyperbolic lattice points},  \emph{Q. J. Math.} \textbf{58} (2007), 281--295.

\bibitem{BZdirections}
\bgroup{}F.~P. Boca\egroup{} and \bgroup{}A.~Zaharescu\egroup{}, {On the correlations of directions in the Euclidean plane},  \emph{Trans. Amer. Math. Soc.} \textbf{358} (2006), 1797--1825.

\bibitem{Bredikhin}
\bgroup{}B.~M. Bredikhin\egroup{}, {Binary additive problems of indeterminate type. II. Analogue of the problem of Hardy and Littlewood}, \emph{Izv. Akad. Nauk SSSR Ser. Mat.} \textbf{27} (1963), 577--612.

\bibitem{BTSD}
\bgroup{}R.~de~la Bret\`eche\egroup{} and \bgroup{}G.~Tenenbaum\egroup{}, {Remarks on the Selberg--Delange method},  \emph{Acta Arith.} \textbf{200} no.~4 (2021), 349--369, corrected version arXiv:2010.12929v8, 24 June 2025.

\bibitem{BH}
\bgroup{}A.~E. Brouwer\egroup{} and \bgroup{}W.~H. Haemers\egroup{}, \emph{{Spectra of Graphs}}, \emph{Universitext}, Springer, New York, 2012.

\bibitem{FT}
\bgroup{}{\'E}.~Fouvry\egroup{} and \bgroup{}G.~Tenenbaum\egroup{}, {Multiplicative functions in large arithmetic progressions and applications},  \emph{Trans. Amer. Math. Soc.} \textbf{375} no.~1 (2022), 245--299.

\bibitem{FI98}
\bgroup{}J.~Friedlander\egroup{} and \bgroup{}H.~Iwaniec\egroup{}, {The polynomial $X^2+Y^4$ captures its primes},  \emph{Ann. of Math. (2)} \textbf{148} (1998), 945--1040.

\bibitem{FI}
\bgroup{}J.~Friedlander\egroup{} and \bgroup{}H.~Iwaniec\egroup{}, \emph{{Opera de Cribro}}, \emph{Amer. Math. Soc. Colloq. Publ.} \textbf{57}, American Mathematical Society, Providence, RI, 2010.

\bibitem{FMPSeo}
\bgroup{}T.~Furuta\egroup{}, \bgroup{}J.~Mi\'ci\'c~Hot\egroup{}, \bgroup{}J.~Pe\v{c}ari\'{c}\egroup{}, and \bgroup{}Y.~Seo\egroup{}, \emph{{Mond--Pe\v cari\'c Method in Operator Inequalities}}, \emph{Monographs in Inequalities} \textbf{1}, Element, Zagreb, 2005.

\bibitem{GreenSiegel}
\bgroup{}B.~Green\egroup{}, {On Landau--Siegel zeros}, Hexagon:2610.00311v2, October~17, 2026.

\bibitem{HB02}
\bgroup{}D.~R. Heath-Brown\egroup{}, {The density of rational points on curves and surfaces}, with an appendix by J.-L.~Colliot-Th\'el\`ene, \emph{Ann. of Math. (2)} \textbf{155} (2002), 553--598.

\bibitem{HB12}
\bgroup{}D.~R. Heath-Brown\egroup{}, {Square-free values of $n^2+1$}, \emph{Acta Arith.} \textbf{155} (2012), 1--13.

\bibitem{Helfgott}
\bgroup{}H.~A. Helfgott\egroup{}, {Expansion, divisibility and parity: an explanation},  in \emph{Combinatorial and Additive Number Theory V}, \emph{Springer Proc. Math. Stat.} \textbf{395}, Springer, Cham, 2022, pp.~199--237.

\bibitem{Iwaniec72}
\bgroup{}H.~Iwaniec\egroup{}, {Primes of the type $\varphi(x,y)+A$, where $\varphi$ is a quadratic form},  \emph{Acta Arith.} \textbf{21} (1972), 203--234.

\bibitem{Iwaniec76}
\bgroup{}H.~Iwaniec\egroup{}, {The half dimensional sieve},  \emph{Acta Arith.} \textbf{29} (1976), 69--95.

\bibitem{IK}
\bgroup{}H.~Iwaniec\egroup{} and \bgroup{}E.~Kowalski\egroup{}, \emph{{Analytic Number Theory}}, \emph{Amer. Math. Soc. Colloq. Publ.} \textbf{53}, American Mathematical Society, Providence, RI, 2004.

\bibitem{KMT}
\bgroup{}O.~Klurman\egroup{}, \bgroup{}A.~P. Mangerel\egroup{}, and \bgroup{}J.~Ter\"av\"ainen\egroup{}, {Multiplicative functions in short arithmetic progressions},  \emph{Proc. Lond. Math. Soc. (3)} \textbf{127} no.~2 (2023), 366--446.

\bibitem{Koukou}
\bgroup{}D.~Koukoulopoulos\egroup{}, \emph{{The Distribution of Prime Numbers}}, \emph{Grad. Stud. Math.} \textbf{203}, American Mathematical Society, Providence, RI, 2019.

\bibitem{Landau}
\bgroup{}E.~Landau\egroup{}, {\"Uber die Einteilung der positiven ganzen Zahlen in vier Klassen nach der Mindestzahl der zu ihrer additiven Zusammensetzung erforderlichen Quadrate},  \emph{Arch. Math. Phys. (3)} \textbf{13} (1908), 305--312.

\bibitem{Linnik}
\bgroup{}{\relax Yu}.~V. Linnik\egroup{}, \emph{{The Dispersion Method in Binary Additive Problems}}, \emph{Transl. Math. Monogr.} \textbf{4}, American Mathematical Society, Providence, RI, 1963.

\bibitem{LWW}
\bgroup{}G.~L\"u\egroup{}, \bgroup{}Z.~Wang\egroup{}, and \bgroup{}J.~Wu\egroup{}, {On the Hardy--Littlewood--Chowla conjecture}, arXiv:2610.20242v1, 2026.

\bibitem{Mont71}
\bgroup{}H.~L. Montgomery\egroup{}, \emph{{Topics in Multiplicative Number Theory}}, \emph{Lecture Notes in Math.} \textbf{227}, Springer-Verlag, Berlin--New York, 1971.

\bibitem{MV73}
\bgroup{}H.~L. Montgomery\egroup{} and \bgroup{}R.~C. Vaughan\egroup{}, {The large sieve},  \emph{Mathematika} \textbf{20} (1973), 119--134.

\bibitem{MV}
\bgroup{}H.~L. Montgomery\egroup{} and \bgroup{}R.~C. Vaughan\egroup{}, \emph{{Multiplicative Number Theory I. Classical Theory}}, \emph{Cambridge Stud. Adv. Math.} \textbf{97}, Cambridge University Press, Cambridge, 2007.

\bibitem{MVII}
\bgroup{}H.~L. Montgomery\egroup{} and \bgroup{}R.~C. Vaughan\egroup{}, \emph{{Multiplicative Number Theory II. Primes and Sieves}}, \emph{Cambridge Stud. Adv. Math.} \textbf{218}, Cambridge University Press, Cambridge, 2026.

\bibitem{MVIII}
\bgroup{}H.~L. Montgomery\egroup{} and \bgroup{}R.~C. Vaughan\egroup{}, \emph{{Multiplicative Number Theory III}}, manuscript in preparation, \nolinkurl{https://personal.science.psu.edu/rcv4/Vol3/Vol3.pdf}.

\bibitem{Motohashi}
\bgroup{}Y.~Motohashi\egroup{}, {On the distribution of prime numbers which are of the form $x^2+y^2+1$},  \emph{Acta Arith.} \textbf{16} (1970), 351--363.

\bibitem{MotohashiII}
\bgroup{}Y.~Motohashi\egroup{}, {On the distribution of prime numbers which are of the form $x^2+y^2+1$. II},  \emph{Acta Math. Acad. Sci. Hungar.} \textbf{22} (1971/72), 207--210.

\bibitem{PD}
\begingroup\raggedright
\bgroup{}{OpenAI}\egroup{}, {The Poisson--Dirichlet law for prime predecessors}, preprint, September~24, 2026, 74 pp., \nolinkurl{https://github.com/openai/math/blob/main/preprints/The-Poisson-Dirichlet-Law-for-Prime-Predecessors-September-24-2026/paper.pdf}.
\par\endgroup

\bibitem{Pascadi}
\bgroup{}A.~Pascadi\egroup{}, {Non-abelian amplification and bilinear forms with Kloosterman sums},  \emph{Geom. Funct. Anal.} \textbf{36} (2026), 1174--1243.

\bibitem{Pilatte}
\bgroup{}C.~Pilatte\egroup{}, {Improved bounds for the two-point logarithmic Chowla conjecture},  \emph{J. Amer. Math. Soc.} (2026), published electronically September 10, 2026.

\bibitem{Rudnick18}
\bgroup{}Z.~Rudnick\egroup{}, {A metric theory of minimal gaps}, \emph{Mathematika} \textbf{64} (2018), 628--636.

\bibitem{Schur}
\bgroup{}I.~Schur\egroup{}, {Bemerkungen zur Theorie der beschr\"ankten Bilinearformen mit unendlich vielen Ver\"anderlichen},  \emph{J. Reine Angew. Math.} \textbf{140} (1911), 1--28.

\bibitem{Shiu}
\bgroup{}P.~Shiu\egroup{}, {A Brun--Titchmarsh theorem for multiplicative functions},  \emph{J. Reine Angew. Math.} \textbf{313} (1980), 161--170.

\bibitem{Wongkew}
\bgroup{}R.~Wongkew\egroup{}, {Volumes of tubular neighbourhoods of real algebraic varieties},  \emph{Pacific J. Math.} \textbf{159} (1993), 177--184.

\end{thebibliography}
\end{document}